\documentclass{siamart1116}

\usepackage{lipsum}
\usepackage{amsfonts}
\usepackage{amsopn}
\usepackage{subeqnar}
\usepackage{subfig}
\usepackage{overpic}
\usepackage{tikz}
\usetikzlibrary{calc,intersections,positioning}
\usepackage{graphicx}
\usepackage{epstopdf}
\usepackage{algorithm,algpseudocode,caption}

\DeclareCaptionFormat{custom}{%
\begin{minipage}[t]{0.2\linewidth}\vspace*{0pt}#1#2\end{minipage}%
\begin{minipage}[t]{0.8\linewidth}\vspace*{0pt}#3\end{minipage}%
}
\Crefname{ALC@unique}{Line}{Lines} 
\ifpdf
  \DeclareGraphicsExtensions{.eps,.pdf,.png,.jpg}
\else
  \DeclareGraphicsExtensions{.eps}
\fi

\numberwithin{theorem}{section}

\newcommand{\TheTitle}{Optimized Sampling for Multiscale Dynamics} 
\newcommand{\TheAuthors}{K. Manohar, E. Kaiser, S. L. Brunton, and J. N. Kutz}

\headers{\TheTitle}{\TheAuthors}

\title{{\TheTitle}\thanks{
{\bf Funding:} JNK acknowledges support from the Air Force Office of Scientific Research
      (AFOSR) grant FA9550-15-1-0385. SLB acknowledges support from the Army Research Office (W911NF-17-1-0306).  SLB and JNK acknowledge support from the Defense Advanced Research Projects Agency (DARPA contract HR011-16-C-0016).  KM and SLB acknowledge support from the Boeing Corporation (SSOW-BRT-W0714-0004).}}

\author{
  Krithika Manohar\thanks{Department of Applied Mathematics, University of Washington, Seattle, WA
    (\email{kmanohar@uw.edu}, \email{kutz@uw.edu}, \url{website}).}
  \and
  Eurika Kaiser\thanks{Department of Mechanical Engineering, University of Washington, Seattle, WA (\email{eurika@uw.edu}, \email{sbrunton@uw.edu})}
  \and  
  Steven L. Brunton\footnotemark[3]
  \and
  J. Nathan Kutz\footnotemark[2]
}

\usepackage{amsopn}

\ifpdf
\hypersetup{
  pdftitle={\TheTitle},
  pdfauthor={\TheAuthors}
}
\fi

\newcommand{\ba}{\mathbf{a}}
\newcommand{\bA}{\mathbf{A}}

\newcommand{\be}{\mathbf{e}}
\newcommand{\bx}{\mathbf{x}}
\newcommand{\by}{\mathbf{y}}

\newcommand{\bC}{\mathbf{C}}

\newcommand{\bR}{\mathbf{R}}
\newcommand{\bQ}{\mathbf{Q}}
\newcommand{\bX}{\mathbf{X}}
\newcommand{\bSigma}{\mathbf{\Sigma}}

\newcommand{\bgam}{\boldsymbol\gamma}
\newcommand{\Ind}{\boldsymbol\beta}
\newcommand{\Sub}{\boldsymbol\alpha}

\newcommand{\bphi}{\boldsymbol{\phi}}
\newcommand{\bPhi}{\boldsymbol{\Phi}}

\newcommand{\reals}{\mathbb{R}}
\newcommand{\ie}{i.e.~}

\DeclareMathOperator*{\argmax}{arg\rm{}max}
\DeclareMathOperator*{\argmin}{arg\rm{}min}

\DeclareMathOperator*{\Max}{maximize}

\DeclareMathOperator{\vol}{vol}

\graphicspath{{update/}}

\begin{document}

\maketitle

\begin{abstract}
The characterization of intermittent, multiscale and transient dynamics using data-driven analysis remains an open challenge. 
We demonstrate an application of the Dynamic Mode Decomposition (DMD) with sparse sampling for the diagnostic analysis of multiscale physics. The DMD method is an ideal spatiotemporal matrix decomposition that correlates spatial features of computational or experimental data to periodic temporal behavior.  DMD can be modified into a multiresolution analysis to separate complex dynamics into a hierarchy of multiresolution timescale components, where each level of the hierarchy divides dynamics into distinct background (slow) and foreground (fast) timescales. 
The multiresolution DMD is capable of characterizing nonlinear dynamical systems in an equation-free manner by recursively decomposing the state of the system into low-rank spatial modes and their temporal Fourier dynamics. 
Moreover, these multiresolution DMD modes can be used to determine sparse sampling locations which are nearly optimal for dynamic regime classification and full state reconstruction.
Specifically, optimized sensors are efficiently chosen using QR column pivots of the DMD library, thus avoiding an NP-hard selection process. 
We demonstrate the efficacy of the method on several examples, including global sea-surface temperature data, and show that only a small number of sensors are needed for accurate global reconstructions and classification of El Ni\~no events.
%
\end{abstract}




\section{Introduction}\label{Sec:Intro}

The accurate and efficient modeling and computation of multiscale, spatio-temporal phenomena remains a grand challenge across almost every physical, biological and engineering discipline.  Indeed, the dynamic interactions that persist across multiple timescales are typical of many complex systems in fluid dynamics, material science, atmospheric and ocean interactions, networked systems and neurobiology.  
Remarkably, the interactions that occur at various micro- and macro-scales generate phenomena that are inherently low-rank, i.e. many multiscale systems manifest dominant coherent patterns (attractors) of activity that can have disparate times-scales.  In such situations, we can leverage dimensionality reduction techniques to obtain interpretable models directly from data for downstream prediction and control.  However, complex interactions of spatiotemporal scales and limited measurement capabilities confound efforts to separate the various micro- and macro-scale physics.  To alleviate this critical limitation, we develop a multiresolution analysis (MRA)~\cite{wavelet,wavelet2} algorithm for the separation of multiscale, low-rank phenomena by its intrinsic timescale in order to determine optimal measurement locations.    We build on the multiresolution dynamic mode decomposition (mrDMD)~\cite{Kutz2016siads} which capitalizes on a wavelet-like decomposition for MRA and embeds the dynamics at each timescale using a low-rank feature extraction technique such as principal component analysis (PCA).  Our method further capitalizes on the extracted low-rank features in order to optimally sample the multiscale physics.

The measurement of multiscale phenomena is critical for characterizing complex systems for prediction and control.  Specifically, sensor placement is central for estimating the full state dynamics and coherent structures in such applications as ocean monitoring~\cite{Yildirim2009oceanmod}, fluid dynamics~\cite{nair2015network,Kaiser2018jcp,Bai2015springer}, neural stimulation in the brain~\cite{Brunton2016jns}, and control~\cite{Han1999sms}.  Given the cost and limitations of measurements in such systems, optimizing sensor placement is a central mathematical challenge.  Of specific concern in this work is the ability of sensors to respect the multiscale dynamics induced in the complex system being measured.  By taking advantage of the multiscale features extracted from our mrDMD, we can move to a  hyperreduction framework whereby the minimal placement of sensors is determined through optimization.  Thus our methodology simultaneously solves a sparse sampling problem by leveraging meaningful features in data, and also accounts for the different spatiotemporal scales in complex systems.

The workhorse algorithm behind mrDMD is the dynamic mode decomposition (DMD)~\cite{Schmid2010jfm,Rowley2009jfm,Tu2014jcd}, a matrix decomposition that identifies low-rank spatial structures and their corresponding temporal Fourier dynamics~\cite{Kutz2016book}. 
This separation is contrasted with variance or energy-based reductions such as proper orthogonal decomposition~\cite{Kutz:2013,benner2015survey} (POD, alternatively known as PCA, Karhunen-Lo\`{e}ve expansion, Empirical Orthogonal Functions (EOF), etc).
POD identifies spatial modal structures based on temporal correlations and is commonly used as an intermediate low-rank projection within DMD.   DMD can be considered a combination of the Discrete Fourier Transform and Proper Orthogonal Decomposition for preserving both temporal correlation and frequency information.  It further provides interpretable, low-rank features in data that can directly be used for prediction or for building dynamic models.  

The low-rank patterns of DMD can be further exploited for sparse (greedy) sampling of the dynamics.  Thus our secondary goal for multiscale characterization is in determining optimal sensor placement.  This requires a principled measurement strategy for spatially sampling the multiscale phenomena of interest. Reduced sensors are especially desirable for forecasting localized phenomena in high-dimensional data or in regimes where sensors are costly or limited.  However, optimizing discrete sensor locations for specific downstream tasks involves a combinatorial search over spatial gridpoints -- an intractable NP-hard computation for even moderately sized grids. Instead, we sample multiscale features using greedy matrix pivoting methods from reduced order modeling~\cite{benner2015survey}, hence bypassing the discrete combinatorial search. In the reduced order modeling (ROM) context this sampling is also known as {\em hyperreduction}, which is distinct from the feature reduction step.  We show that mrDMD correctly identifies dynamics occurring at different timescales in multiscale data. Secondly, multiscale sensors derived from QR matrix pivots of mrDMD modes are shown to spatially cover localized coherent structures exhibiting specific frequencies. We then leverage resulting mrDMD modes and multiscale sensors to estimate temporal behavior, for which multiscale reductions are demonstrably more accurate than POD-based reduction. 


\section{Background}\label{Sec:Background}

This section develops the background theory of dynamic mode decomposition and its multiresolution variant mrDMD. Our goal is to construct a spatiotemporal decomposition directly from time snapshots of high-dimensional system states $\bx\in\reals^n$, such that their evolution can be expressed as a linear combination of $r$ spatial modes $\bphi_k(\boldsymbol{\xi})$ governed by time dynamics $a(t)$
\begin{equation}\label{eq:at}
\bx(t) = \sum_{k=1}^r a_k(t)\boldsymbol{\phi}_k(\boldsymbol{\xi}).
\end{equation}
When $r\ll n$, this is considered a low-rank embedding of the dynamics.  For POD, the modes $\boldsymbol{\phi}_k(\boldsymbol{\xi})$ are determined from a singular value decomposition of a data matrix, and the dynamics are then projected into this space~\cite{Kutz:2013,benner2015survey}.
DMD and mrDMD instead enforce that each mode be further constrained so that $a_k(t)=b_k \exp(i\omega_k t)$.

\subsection{Dynamic mode decomposition}

The DMD originated as a spectral decomposition for identifying coherent structures in fluid dynamics~\cite{Schmid2010jfm,Rowley2009jfm}. Algorithmic refinements such as exact DMD~\cite{Tu2014jcd} and practical developments~\cite{Kutz2016book} quickly followed, establishing DMD as a rigorous data-driven framework for spatiotemporal analysis and prediction. Many additional variants have been developed that capitalize on sparse,compressive measurements~\cite{Jovanovic2014pof,Brunton2015jcd}, input-output systems and control~\cite{Proctor2016siads}, improved denoising~\cite{askham2017variable}, multiscale structure~\cite{Kutz2016siads}, randomized linear algebra~\cite{bistrian2017randomized,erichson2017randomized}, and streaming data~\cite{hemati2014dynamic,pendergrass2016streaming}.  

The DMD separates spatiotemporal dynamics into a linear decomposition of coherent structures (DMD eigenmodes) with fixed temporal behavior (DMD eigenvalues). In particular, the DMD has connections to the Koopman operator which acts as the forward operator on scalar observables of a system. Given a discrete-time dynamical system 
\begin{equation}
\label{eqn:nl_sys}
\bx_{k+1} = \mathbf{F}(\bx_k),
\end{equation}
the Koopman operator $\mathcal{K}$ is defined as the linear forward operator acting on all  scalar observables of the state~\cite{Koopman1931pnas,Koopman1932pnas}
\begin{equation}
\mathcal{K}g(\bx_k) = g(\mathbf{F}(\bx_k)) = g(\bx_{k+1}),
\end{equation}
thus linearizing nonlinear dynamics to completely characterize the system. The approximation of Koopman eigenfunctions and eigenvalues through carefully selected observables is an active area of research which has given rise to  Koopman mode decomposition~\cite{Budivsic2012chaos}, extended DMD~\cite{Williams2015jns}, and kernel DMD~\cite{Williams2014kernel}. DMD can be considered a special case of Koopman decomposition which restricts observables $g$ to be point measurements of state. Thus the inputs to the DMD algorithm are the following data matrices of the measured state trajectory
\begin{align}
\label{eqn:data_matrix}
\bX &= \begin{bmatrix} \bx_1 & \bx_2 &\dots &\bx_{m-1} \end{bmatrix} \nonumber \\ 
\bX' &= \begin{bmatrix} \bx_2 & \bx_3 &\dots &\bx_{m} \end{bmatrix}, 
\end{align}
assumed locally related by a linear forward evolution map $\bA$ 
\begin{equation}
\label{eqn:linear_sys}
\bX' = \bA\bX.
\end{equation}
Hence, DMD fits a linear approximation $\bA$ to the system~\eqref{eqn:nl_sys}, but without explicitly computing $\bA$ using an expensive pseudoinverse operation. Instead, exact DMD computes its eigendecomposition $\bA\bPhi = \bPhi\mathbf{\Lambda}$ using low-rank approximations of the data matrices. Rank truncation is achieved using the singular value decomposition to obtain $r\ll n$ eigenvectors and eigenvalues which evolve the system forward in time -- a more efficient and interpretable representation of dynamics than $\bA$ alone. Each state can be efficiently computed from a linear combination of DMD eigenvectors or eigenmodes (columns $\bphi_k$ of $\bPhi$), DMD eigenvalues ($\lambda_k=\mathbf{\Lambda}_{kk}$) and corresponding modal amplitudes $\mathbf{b}\in\reals^k$
\begin{equation}
\bx_k = \bPhi\mathbf{\Lambda}^k\mathbf{b}.
\end{equation}
The equivalent expression in the continuous time setting uses a convenient scaling of the eigenvalues $\omega_k = \frac{\log{\lambda_k}}{i\Delta t}$
\begin{equation}
\label{eqn:dmd_approx}
\bx (t) = \sum_{k=1}^r b_k\bphi_k(\boldsymbol{\xi})\exp(i\omega_k t).
\end{equation}
{\em Remark: The DMD approximation as stated introduces systemic bias in the eigenvalue computation of noisy data\cite{Dawson2016EF,Hemati2017tcfd}. To correct for this, we construct a weighted approximation of $\bA$ that incorporates both forward ($\bX'=\bA_f\bX$) and backward ($\bX = \bA_b\bX'$) time evolution, as formulated by Dawson et al~\cite{Dawson2016EF}
\begin{equation}
\bA = (\bA_f\bA_b^{-1})^{1/2}.
\end{equation}
This forward-backward DMD method is used throughout our results. }

DMD provides a compelling description of many physical systems with low-rank, periodic temporal behavior, including fluid dynamics~\cite{Tu2014ef,Rowley2009jfm}, ocean sciences~\cite{Kutz2016siads} and neuroscience~\cite{Brunton2016jns}.  The DMD approximation fails, however, when approximating intermittent and transient phenomenon, both of which violate the constrained temporal response in (\ref{eqn:dmd_approx}).  In such cases, Fourier modes in time are insufficient to capture a richer set of dynamics.

\subsection{Multiresolution DMD}

\begin{figure}[t]
\begin{center}
\vspace*{0.2in}
\begin{overpic}[width=0.45\textwidth]{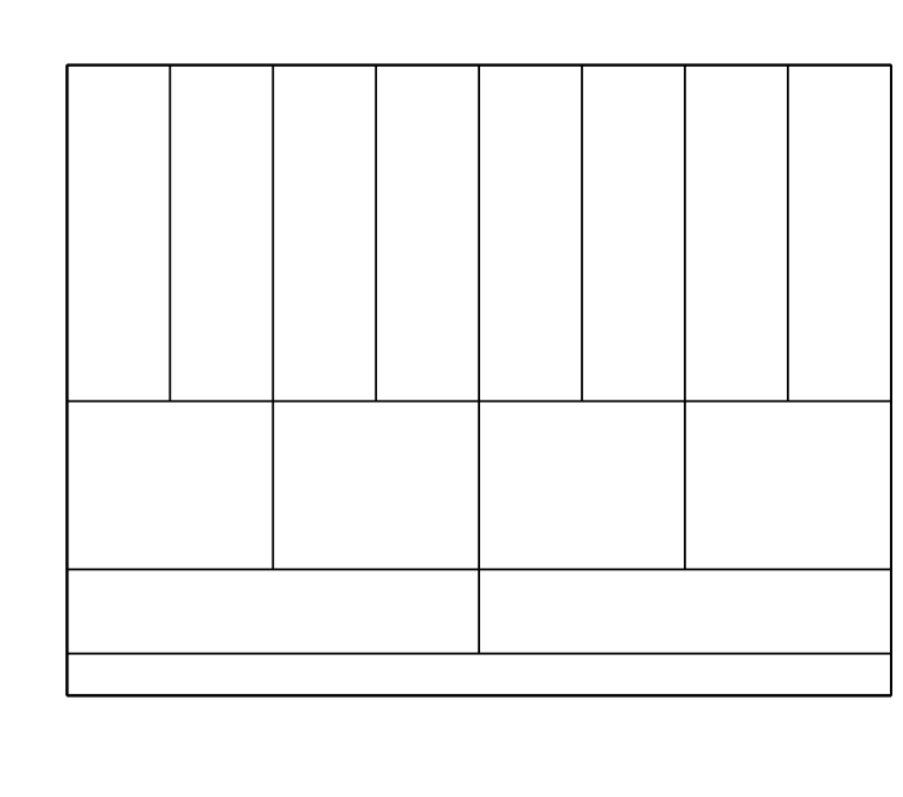}
\put(100,9){$\bphi_k^{(1,1)}$}
\put(30,18.25){$\bphi_k^{(2,1)}$}
\put(73,18.25){$\bphi_k^{(2,2)}$}
\put(15,32){$\bphi_k^{(3,1)}$}
\put(35,32){$\bphi_k^{(3,2)}$}
\put(60,32){$\bphi_k^{(3,3)}$}
\put(80,32){$\bphi_k^{(3,4)}$}
\put(7,90){$\bphi_k^{(4,1)} \,\,\, \cdots$}
\put(87,90){$\bphi_k^{(4,8)}$}
\put(10,87){\rotatebox{-90}{$\xrightarrow{\hspace*{0.5cm}}$}}
\put(90,87){\rotatebox{-90}{$\xrightarrow{\hspace*{0.5cm}}$}}
\put(90,15.5){\rotatebox{-180}{$\xrightarrow{\hspace*{0.7cm}}$}}
\put(40,5){time ($t$)}
\put(0,27){\rotatebox{90}{frequency ($\omega$)}}
\normalsize
\end{overpic}
\hspace{2em}
\begin{tikzpicture}
	\node (Phi) at (0,0) [draw,minimum width=2.45cm,minimum height=4cm] {$\bphi_k^{(1,1)} \dots \bphi_k^{(4,8)}$};
	\foreach \x in {-.75,-.3,.3,.75}{
  	\draw [help lines,dotted] (\x,-2) -- (\x,2);}
	\node [below=1.5mm of Phi] {$\bPhi$};
	\node [above=2mm of Phi] {DMD library};	
	\node (s1) [draw,thick,red,minimum width=2.7cm, minimum height=.25cm, below=-37mm of Phi] {};
	\node (s2) [draw,thick,red,minimum width=2.7cm, minimum height=.25cm,below=-32mm of Phi] {};
	\node (s3) [draw,thick,red,minimum width=2.7cm, minimum height=.25cm, below=-15mm of Phi] {};
	\node (s4) [draw,thick,red,minimum width=2.7cm, minimum height=.25cm,below=-5mm of Phi] {};
	\node[inner sep=0pt,draw=black,thick,right=5mm of Phi] (map) 
	{\includegraphics[width=.2\textwidth]{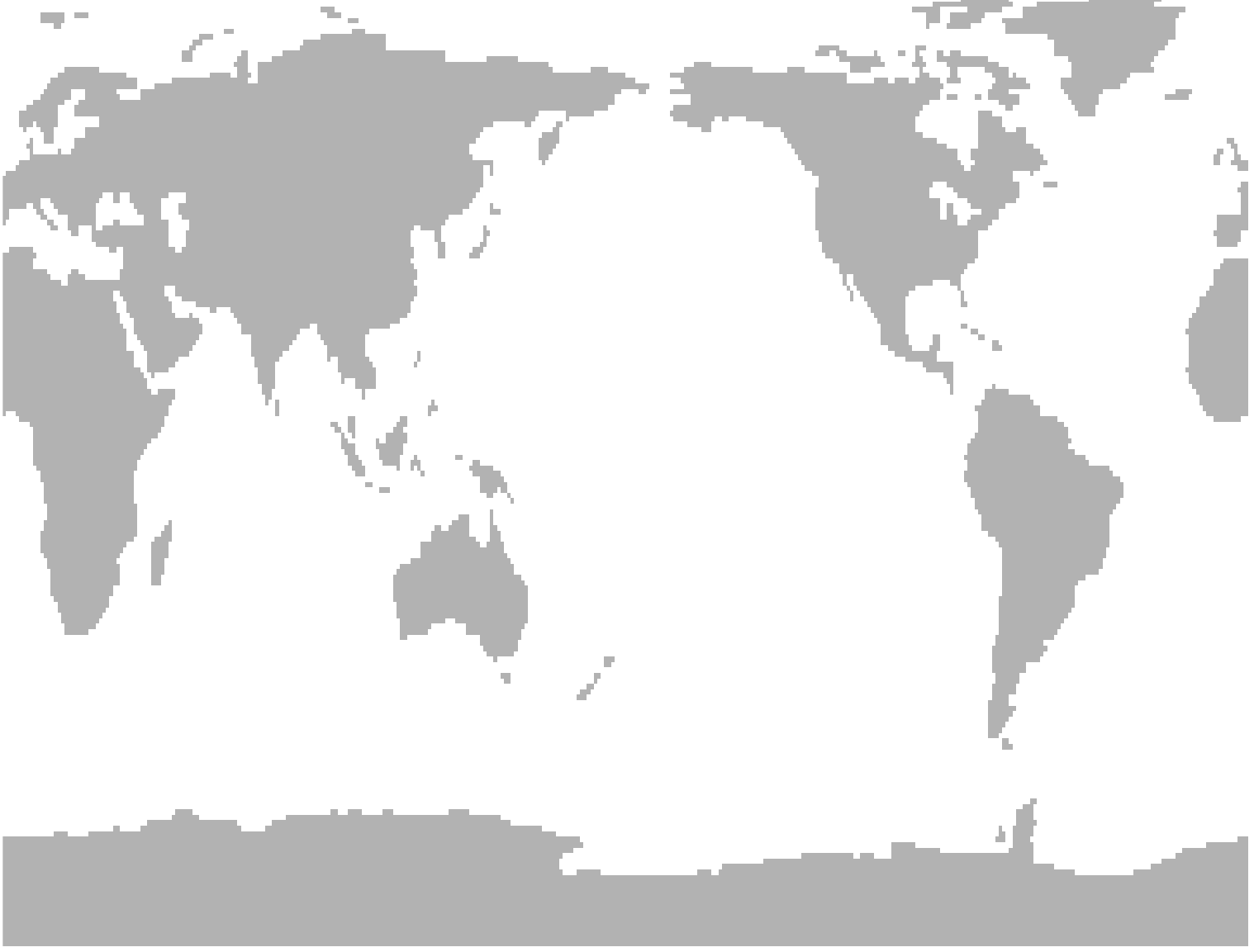}};
	\fill[orange]  (3,.5) circle (0.5ex); 
	\draw[->,thick,red] (s1.east) -- (3,.5) ;
	\fill[orange]  (2.7,.3) circle (0.5ex); 
	\draw[->,thick,red] (s2.east) -- (2.7,.3) ;
	\fill[orange]  (4,.2) circle (0.5ex); 
	\draw[->,thick,red] (s3.east) -- (4,.2) ;
	\fill[orange]  (3.6,-.05) circle (0.5ex); 
	\draw[->,thick,red] (s4.east) -- (3.6,-.05) ;
	\node [above=1.5mm of map] {Sparse sensors};
	\end{tikzpicture}
\end{center}
\caption{\label{fig:wavelet} Illustration of the multiresolution analysis and sampling framework.  On the left are the mrDMD modes $\bphi_k^{\ell,j}({\boldsymbol{\xi}})$ and their position in the hierarchy.  The triplet of integer values, $\ell, j$ and $k$, uniquely expresses the time level, bin and mode of
the decomposition. Depicted in the middle is the matrix library of all mrDMD modes, which are then leveraged to obtain optimal samples (spatial sensors, right) for downstream estimation and prediction tasks.  }
\end{figure}
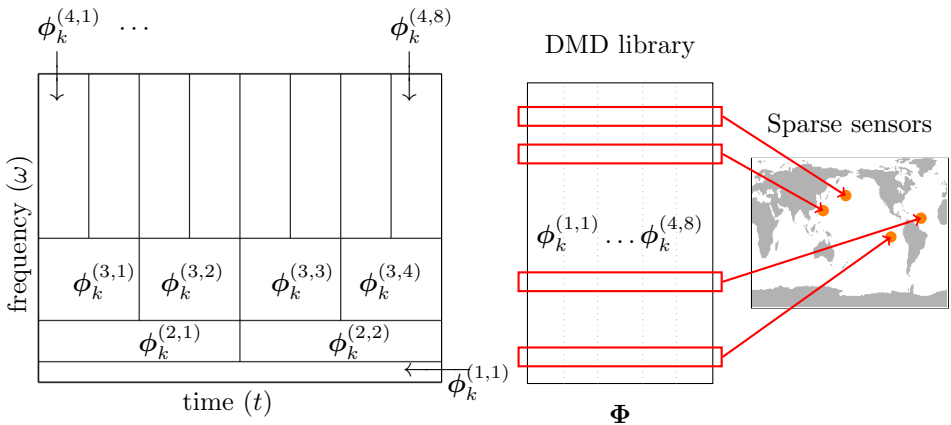

The mrDMD method~\cite{Kutz2016siads,Kutz2016book,Erichson2015mrDMD} circumvents some of the common shortcomings of standard DMD.  Specifically, mrDMD provides a hierarchical temporal sampling framework, much like a wavelet decomposition, whereby disparate timescale phenomena can be characterized recursively.
The recursive nature of mrDMD, illustrated in Fig.~\ref{fig:wavelet}, can readily characterize distinct time scales by removing micro- and macro-timescale modes at different hierarchical levels.
The mathematical structure of mrDMD accounts for the number of levels ($L$) of the decomposition, the number of time bins ($J$) for each level, and the number of modes retained at each level ($m_L$).  Thus the solution is parametrized by the following three indices:
\begin{subeqnarray*}
  &&  \ell = 1, 2, \cdots , L, \mbox{ where } L= \mbox{\# of decomposition levels } \\
  &&  j = 1, 2, \cdots, J  \,\,\, \mbox{\# time bins per level} \,\, (J=2^{(\ell-1)})  \\
  &&  k= 1, 2, \cdots, m_{L} \,\,\, \mbox{\# of modes extracted at level $L$.}
\end{subeqnarray*}
To formally define the series solution for ${\bf x}_{\mbox{\tiny mrDMD}}(t)$, we propose the following indicator function
\begin{equation}
  f_{\ell,j} (t) \!=\! \left\{ \!\!  \begin{array}{cl}  1 & t\in[t_j,t_{j+1}] \\ 0 & \mbox{elsewhere}  \end{array}    \right.
  \,\, \mbox{with} \,\,\, j=1,2,\cdots, 2^{(\ell-1)}
  \label{eq:indicator}
\end{equation}
which is only non-zero in the interval, or time bin, associated with the value of $j$.  The parameter
$\ell$ denotes the level of the decomposition.

The three indices and  indicator function (\ref{eq:indicator}) give the mrDMD solution expansion
\begin{equation}
  {\bf x}_{\mbox{\tiny mrDMD}}(t)  = \sum_{\ell=1}^{L} \sum_{j=1}^{J} \sum_{k=1}^{m_L}
  f_{\ell,j}(t)  b_k^{(\ell,j)} \bphi_k^{(\ell,j)}  \! ({\boldsymbol{\xi}}) \exp(i\omega_k^{(\ell,j)} t) 
  \, .
  \label{eq:dmd3}
\end{equation}
This is a concise definition of the mrDMD solution that includes the information on the level,
time bin location and number of modes extracted.
Figure~\ref{fig:wavelet} demonstrates the mrDMD decomposition in terms of the solution
(\ref{eq:dmd3}).  In particular, each mode is represented in its respective time bin and level.
An alternative interpretation of this solution is that it generalizes the linear mapping (\ref{eqn:linear_sys})
so that at each level $\ell$ of the decomposition one has
\begin{equation}
  {\bf x}^{(\ell,j)}_{j+1} = {\bf A}^{{(\ell,j)}} {\bf x}^{(\ell,j)}_j \, 
\end{equation}
where the matrix ${\bf A}^{(\ell,j)}$ captures the dynamics in a given time bin $j$ at level ${\ell}$.

The indicator function $f_{\ell,j}(t)$ acts as a sifting function for each time bin.  Interestingly,
this function acts as the Gab\'or window of a windowed Fourier transform~\cite{Kutz:2013}.  
Since our sampling bin has a hard cut-off of the time series, it may introduce some artificial
high-frequency oscillations.  Time-series analysis, and wavelets in particular, introduce
various functional forms that can be used in an advantageous way.  Thus thinking 
more broadly, one can imagine using wavelet functions for the sifting operation, allowing
the time function $f_{\ell,j}(t)$ to take the form of one of the many potential wavelet basis, i.e.
Haar, Daubechies, Mexican Hat, etc.


\section{Multiresolution analysis framework}\label{Sec:MSSP}

We propose a multiresolution analysis decomposition and spatial sampling framework for the estimation and prediction of multiscale dynamics. 
Localized time-frequency analysis presents difficulties for future state prediction as the active modes in a new temporal window are unknown. Thus additional information is required to estimate temporal coefficients $\ba(t)$. This additional information is commonly provided in the form of current state observations or sensor measurements. In this section we frame future state prediction as the reconstruction of new high-dimensional states from observations, or equivalently, estimating the correct temporal coefficients within the mrDMD library of modes. This is a common problem in engineering applications; the second half of this section describes how to optimize state measurement locations to minimize prediction error.

\subsection{Mathematical formulation}

The mrDMD yields a {\em library} of possible dynamics which can be arranged in matrix form 
\begin{equation}
\bPhi = \begin{bmatrix}
\bphi_k^{(1,1)} & \bphi_k^{(2,1)} & \bphi_k^{(2,2)} & \bphi_k^{(3,1)} & \dots & \bphi_k^{(3,4)} & \bphi_k^{(4,1)} & \dots & \bphi_k^{(4,8)} 
\end{bmatrix}.
\end{equation}
There may be more than mode per level indexed by $k$ in the library $\bPhi$,  which is now an $n\times M$ matrix.  
Thus any state can be approximated by a linear combination of the columns of $\bPhi$ with time-dependent coefficients
\begin{equation}
\label{eqn:full_state}
  \bx(t) = \bPhi\ba(t).
\end{equation}
The challenge with future state prediction is identifying the subset of $\ba$'s components that are active (nonzero) at a future time, because mrDMD modes are localized in the time-frequency domain and do not enforce globally periodic temporal behavior.

The necessary information is often provided by point observations of the unknown state vector -- an extremely prevalent construct in engineering and controls. Then there exists a linear relationship between observations, stored in a vector $\by$, and time coefficients. The components of $\by$ result from a linear observation map $\bC$ applied to $\bx$ and there is additive white noise $\eta\sim \mathcal{N}(0,\sigma^2)$
\begin{align*}
\by &= \bC\bx + \eta \\
 &= \bC\bPhi\ba + \eta.
\end{align*}
The following assumptions are made: $\by\in\reals^p$, where $p\ll n$, and the measurement matrix $\bC\in\reals^{p\times n}$ consists of $p$ rows of the $n\times n$ identity matrix which index measurements of the full state.
This linear inverse problem for coefficient estimation is generally ill-posed (underdetermined) without additional regularization. We impose two forms of regularization:
\begin{enumerate}
\item {\bf Coefficient sparsity.} To reflect the physical constraint that only a few modes are active at a given time, we impose a sparse structure on $\ba$ which minimizes the number of nonzero entries (see below). 
\item {\bf Designing $\bC$.} Point observations of states will be optimized to leverage mrDMD features of interest and eliminate redundancy between point samples.  In engineering settings this is known as sensor placement.
\end{enumerate}
At first glance both pose combinatorial subset selection problems, but as we shall see, there exist well-known methods to approximate both selections.

\subsection{Online estimation and prediction}\label{Sec:prediction}
Assume a fixed $\bC$ that is already optimized for estimation.
The sparsest coefficients $\ba_\star$ can be approximated using a convex $l_1$ constrained minimization
\begin{equation}
\label{eqn:sparse_pred}
\ba_\star = \argmin_\ba \|\ba\|_1 \mbox{ subject to } \by = \bC\bPhi\ba.
\end{equation}
The convex optimization posed in \eqref{eqn:sparse_pred} can be efficiently computed using basis pursuit~\cite{Tropp2007ieeetit}, which is available in $l_1$ optimization tools including CVX~\cite{Grant2008cvx}, SPGl1~\cite{spgl12007} and CoSaMP~\cite{Needell:2010}.
The full state is subsequently estimated from the active (nonzero) modal dynamics using $\hat\bx \approx \bPhi\ba_\star$ (see \eqref{eqn:full_state}).
A more accurate reconstruction can be achieved with traditional least-squares approximation using only the active modes. Explicitly, we first determine which library modes are active by indexing the nonzero solution coefficients
\begin{equation}
\Ind := \{ i\in\{1,\dots,r\} \mid a_{\star_i} \neq 0 \}.
\end{equation}
Once $\Ind$ is known, then state prediction reduces to a least-squares approximation using the subset of active modes within the library
\begin{align*}
\label{eqn:estimation_ls}
\hat\bx &\approx \bPhi_{\Ind}\ba_{\Ind} \\
&= \bPhi_{\Ind} (\bC\bPhi_{\Ind})^\dagger \by.
\end{align*}
When $\bPhi$ are POD modes, this is widely known as gappy POD, introduced in the early 90s for image reconstruction using eigenfaces~\cite{Everson1995gappy} and subsequently used for flow field reconstruction in ROMs~\cite{Willcox2006compfl}. 
Sparse approximation methods similar to the above description have applied compressed random measurements in POD mode libraries to classify different parameter regimes~\cite{Brunton2014siads,Bright2013pof}. In recent work~\cite{Kramer2017siads} classification is accomplished with sparse measurements in a DMD mode library.
However, to our best knowledge, this is the first known application of sparse estimation using a DMD mode library samples with optimized measurements. 





\subsection{Offline measurement design}

We now turn to the optimal design of the measurement operator to obtain a minimal set of non-redundant observations of the state. The concept of redundancy is implicitly quantified in the multiscale dynamical library, however, the repetition of modes across multiple time-frequency bins introduces redundancy. Optimal measurement design from this library would then result in multiple measurement locations with the same dynamics, which is undesirable in practice. This is remedied by first preprocessing the library and constructing an index set to filter out redundant modes (this is distinct from the local active set $\Ind$). For example, this subset can be selected from modes with dominant spatiotemporal signatures, by identifying mrDMD modal amplitudes which exceed some threshold $T$ 
\begin{equation}
\Sub := \{ i\in\{1,\dots,m\} \mid |b_{i}| > T \}.
\end{equation}
Alternatively, depending on the application, modes corresponding to specific time or frequency bins of the decomposition can be extracted in this step for  downstream sensor selection. 
The mathematical aim is to construct the measurements $\bC$ to minimize the least-squares approximation error, $\ba-\hat\ba$. Similar measurement design methods using the matrix QR factorization have been developed for POD modal bases~\cite{Drmac2016siam,Manohar2016jfs} and reformulated in the estimation setting~\cite{Manohar2017data}. For this step the library $\bPhi_{\Sub}$ and noise levels are assumed predetermined, but point measurements can be chosen to minimize some scalar measure of the ``size" of the error covariance ($\bSigma$):
\begin{equation}
\label{eqn:ynoise}
\bSigma = \mbox{Var}(\ba-\hat\ba) = \sigma^2[(\bC\bPhi_{\Sub})^T\bC\bPhi_{\Sub}]^{-1}.
\end{equation} 
The error covariance $\bSigma$ characterizes the minimum volume {\em $\rho$-confidence ellipsoid}, $\varepsilon_\rho$, that contains the least-squares estimation error, $\ba-\hat\ba$, with probability $\rho$. The D-optimal design criterion minimizes the volume of $\varepsilon_\rho$
\begin{equation}
\label{eqn:ellipsoidvol}
\vol(\varepsilon_\rho) = \delta_{\rho,r}\det\bSigma^{\frac{1}{2}},
\end{equation}
where $\delta_{\rho,r}$ only depends on $\rho$ and $r$, by minimizing the determinant, which is the only sensor dependent term. This optimization is equivalent to maximizing the determinant of the inverse
\begin{equation}
\label{eqn:dopt}
\Max_\bC ~~ \log \det \left[(\bC\bPhi_{\Sub})^T\bC\bPhi_{\Sub}\right].
\end{equation}
Here, we are imposing the following structure on the measurement matrix $\bC\in\reals^{p\times n}$
\begin{equation}
	\bC = [\be_{\gamma_1} ~|~ \be_{\gamma_2} ~|~ \dots ~|~ \be_{\gamma_p}]^T,
\end{equation}
where $\gamma_i\in\{1,\dots,n\}$ indexes the high-dimensional measurement space, and $\be_{\gamma_i}$ are the canonical unit vectors consisting of all zeros except for a unit entry at $\gamma_i$. 
This guarantees that each row of $\bC$ only measures from a single spatial location, corresponding to a point sensor.

The optimal design shrinks all components of the reconstruction error via the associated ellipsoid {\em volume}, \ie the choice that minimizes the determinant of $\bSigma$. This subset optimization is a combinatorial search over $n\choose p$ possibilities which quickly becomes computationally intractable even for moderately large $n$ and $p$. Fortunately, an extremely efficient greedy optimization method for~\eqref{eqn:dopt} is given by the pivoted matrix QR factorization of $\bPhi_{\Sub}^T$. QR pivoting has been extensively used to compute optimal quadrature and interpolation points~\cite{Businger1965nm,Sommariva2009qr,Drmac2016siam,Seshadri2017siam}, which are optimal samples in polynomial basis modes. QR pivoting is shown to outperform related methods for optimal selection in either computational accuracy or runtime, sometimes both~\cite{Manohar2017data}. Such methods include empirical interpolation methods~\cite{Barrault2004crm,Chaturantabut2010siamjsc} from ROMs, information theoretic criteria~\cite{Krause2008jmlr} and convex optimization~\cite{Joshi2009ieee,Boyd2004}.

The reduced matrix QR factorization with column pivoting decomposes a matrix $\bA\in\reals^{m\times n}$ into a unitary matrix $\bQ$, an upper-triangular matrix $\bR$ and a column permutation matrix $\bC$ such that $\bA\bC^T = \bQ\bR$. Recall that the determinant of a matrix, when expressed as a product of a unitary factor and an upper-triangular factor, is the product of the diagonal entries in the upper-triangular factor:
\begin{equation}
\left| \det \bA\bC^T \right| = |\det\bQ||\det\bR| = \prod_i |r_{ii}|,
\end{equation}
The pivoted QR permutes the matrix $\bPhi_{\Sub}^T$ with $\bC^T$ to enforce the following diagonal dominance structure on the diagonal entries of its $\bR$ factor~\cite{Drmac2016siam}:
	\begin{equation}
	\sigma_i^2 = |r_{ii}|^2 \ge \sum_{j=i}^k |r_{jk}|^2; \quad 1\le i \le k \le m.
	\end{equation}
	Column pivoting iteratively increments the volume of the pivoted submatrix by selecting a new pivot column with maximal 2-norm, then subtracting from every other column its orthogonal projection onto the pivot column (see Algorithm \ref{alg:qrpivot}). 	
	In this manner the QR factorization with column pivoting yields $p$ point measurement indices (pivots) that best characterize dominant dynamical modes $\bPhi_{\Sub}$
	\begin{equation}
	\bPhi_{\Sub}^T\bC^T = \bQ\bR.
	\end{equation}
	\begin{algorithm}[tb!]
		\caption{QR factorization with column pivoting of $\bA\in\reals^{n\times m}$ \newline Greedy optimization for placing $p$ sensors $\bgam$ from multiscale modes. \newline Usage: {\scshape QrPivot}( $\bPhi_{\Sub}^T$ , $p$ ) (if $p=m$) \newline {\color{white}Usage:} {\scshape QrPivot}( $\bPhi_{\Sub}\bPhi_{\Sub}^T$ , $p$ ) (if $p>m$)\label{alg:qrpivot}}
		\begin{algorithmic}[1]		
			\Procedure{qrPivot}{ $\bA,~p$ }
			
			\State $\bgam \gets [~~]$
			\For{$k=1,\dots,p$}
			\State $\gamma_k = \argmax_{j\notin \bgam} \|\mathbf{a}_j\|_2 $ 
			\State Find Householder $\tilde{\bQ}$ such that $\tilde{\bQ} \cdot \begin{bmatrix} a_{kk} \\ ~\\ \vdots \\ a_{nk} \end{bmatrix} = \begin{bmatrix}
			\star \\ 0 \\ \vdots \\ 0
			\end{bmatrix}$  \Comment{\parbox[t][][t]{.21\linewidth}{\itshape $\star$'s represent nonzero diagonal entries in $\bR$}}
			\State{$\bA \gets \mbox{diag}(I_{k-1},\tilde{\bQ}) \cdot \bA$ } \Comment{\parbox[t]{.5\linewidth}{\itshape Remove from all columns the orthogonal projection onto $\mathbf{a}_{\gamma_k}$}}
			\State $\bgam\gets [\bgam,~\gamma_k]$
			\EndFor \\
			\Return $\bgam$
			\EndProcedure
		\end{algorithmic}		
	\end{algorithm}
The QR pivoting algorithm is summarized in Algorithm~\ref{alg:qrpivot}. 
The standard pivoting formulation yields only as many pivots as there are columns of $\bPhi_{\alpha}$ (modes). However, {\em oversampling}, which refers to sampling more observations than modes, promotes robustness to measurement noise, thus regularizing the state estimation problem~\eqref{eqn:estimation_ls}. We include here a method first introduced in~\cite{Manohar2017data} for obtaining $p>m$ samples. These samples result from the pivoted QR factorization of $\bPhi_{\Sub}\bPhi_{\Sub}^T$, based on the observation that the singular values of the iverted error covariance $(\bC\bPhi_{\Sub})^T\bC\bPhi_{\Sub}$ coincide with the first $r$ singular values of its transpose and the diagonal entries of its $\bR$ factor. Thus, the pivoting formulation will automatically condition the desired determinant.

\section{Applications}\label{Sec:Applications}

This section demonstrates the application of our multiresolution analysis and sampling methods on data presenting multiscale temporal scales. The first example, a manufactured video, presents three spatial modes independently oscillating at different frequencies in overlapping time windows. We empirically quantify the accuracy of our analysis and estimation from optimized samples by comparing with the known dynamics.
By contrast, the second example seeks to identify intermittent warming phenomena from real global satellite data of sea surface temperatures. Here the accuracy metric is the correct identification of warming events in the validation window based on the same events identified in a separate training window. For both examples, results from the multiresolution analysis are contextualized with appropriate comparisons to POD and DMD.

\subsection{Multiscale video example}

The video data is generated from three spatial modes oscillating at different frequencies in overlapping time intervals, effectively appearing mixed in time.
Mathematically each snapshot is constructed as a linear combination of Gaussian  spatial modes
\begin{equation}
\bx(t) = \sum_{i=1}^3 a_i(t)\exp\left(-\frac{(\boldsymbol\xi-\boldsymbol\xi_{0_i})^2}{w_i}\right),
\end{equation}
where $\boldsymbol\xi\in\reals^{n}$ is the vectorized planar spatial grid ($n=nx\times ny$) and $w$ is a scalar width parameter.
The challenge with mrDMD is to identify the true modes $\bphi_i = \exp\left(-\frac{(\boldsymbol\xi-\boldsymbol\xi_{0_i})^2}{w_i}\right)$ within their respective temporal windows. 
We list explicitly the modal components with their time dynamics in \cref{tab:toy_dynamics} and \cref{fig:toy_true}.

\begin{table}[tbhp]
\caption{Multiscale video dynamics}
\label{tab:toy_dynamics}
\centering
\begin{tabular}{|l|| l |l|}
\hline
Mode & Time dynamics & frequency ($\omega_i$)\\ \hline
$\bphi_1$ & $a_1(t)  = \begin{cases} \exp(2\pi i\omega_1 t) & t\in[0,5]\\ 0&\text{elsewhere} \end{cases} $ & $\omega_1=5.55$ \\ \hline
$\bphi_2$ & $a_2(t)  = \begin{cases} \exp(2\pi i\omega_2 t) & t\in[2.5,7.5]\\ 0&\text{elsewhere}\end{cases}$ & $\omega_2=0.9$ \\ \hline
$\bphi_3$ & $a_3(t)  = \begin{cases} \exp(2\pi i\omega_3 t) & t\in[0,10]\\ 0&\text{elsewhere}\end{cases}$ & $\omega_3=0.15$ \\ \hline
\end{tabular}
\end{table}
\begin{figure}[tbhp]
	\centering
	\begin{tikzpicture}
		\node (map) at (9,1) {\includegraphics[scale=.325]{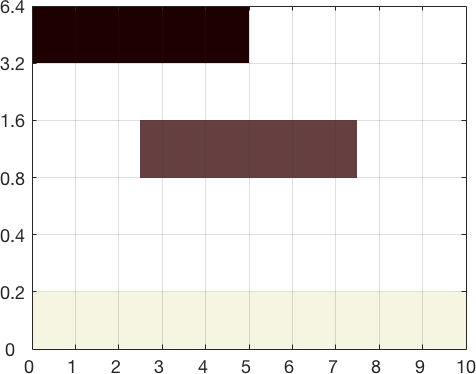}};	
    	\node at (0,3.2) {\frame{\begin{overpic}[scale=.15]{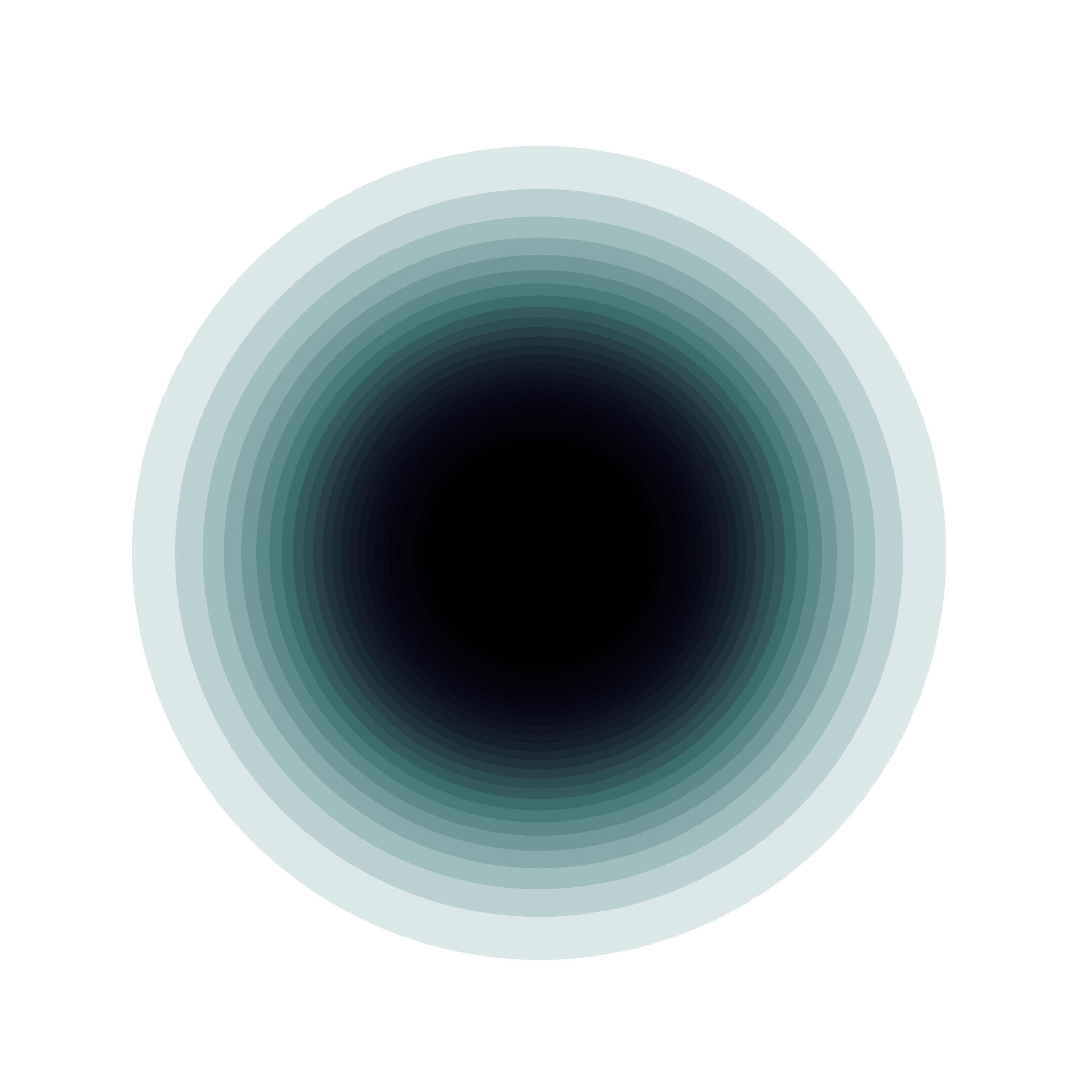}\put(5,82){$\bphi_1$}\end{overpic}}};
    	\node (a1) at (3.4,3.1) {\begin{overpic}[scale=0.3]{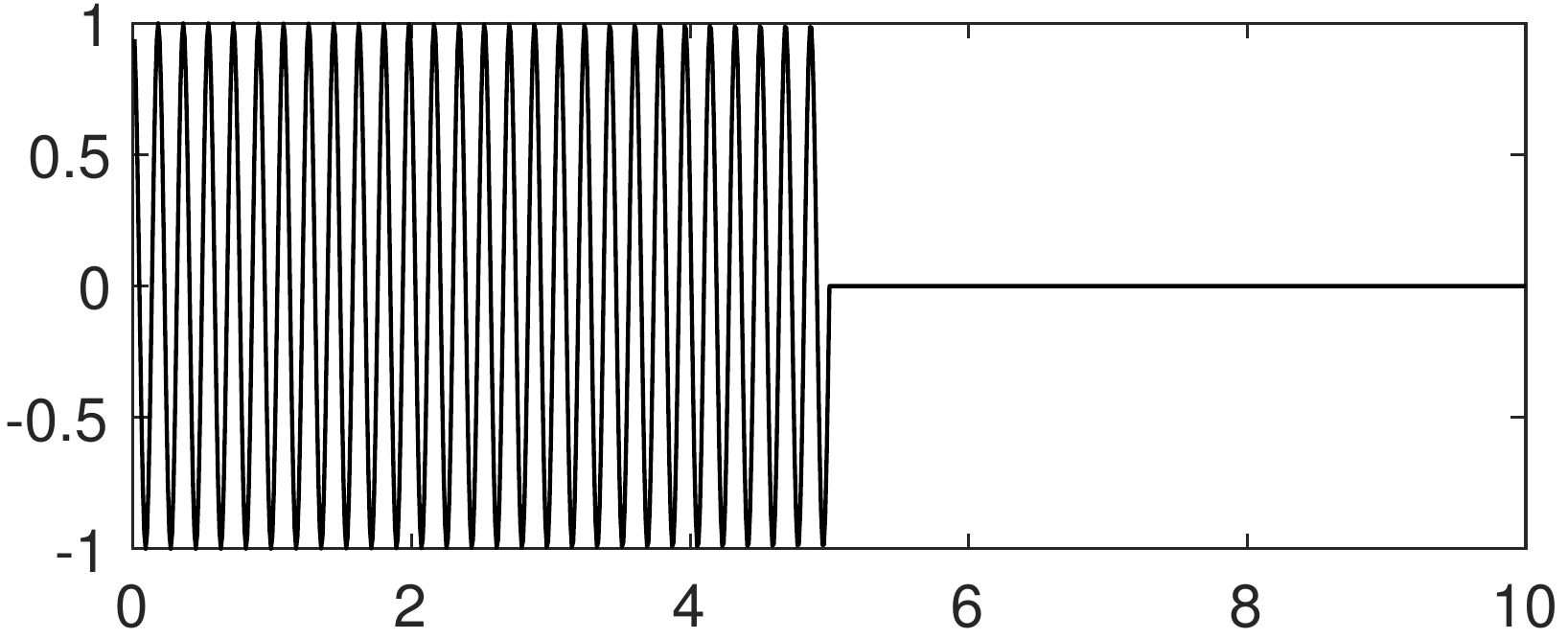}\put(80,32){$\ba_1(t)$}\end{overpic}};
    	\node at (0,1.1) {\frame{\begin{overpic}[scale=.155]{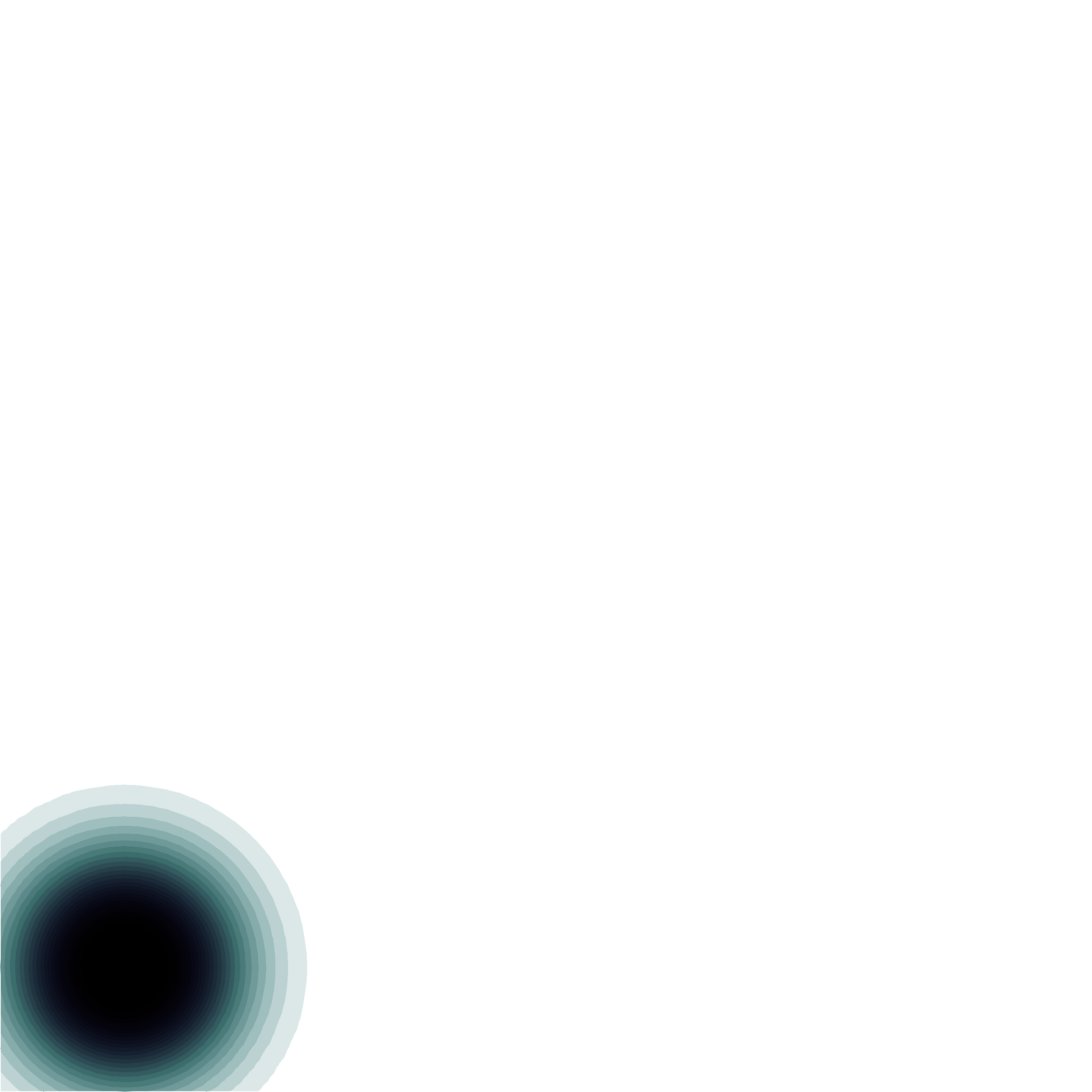}\put(5,82){$\bphi_2$}\end{overpic}}};
    	\node (a2) at (3.45,1.1) {\begin{overpic}[scale=0.3]{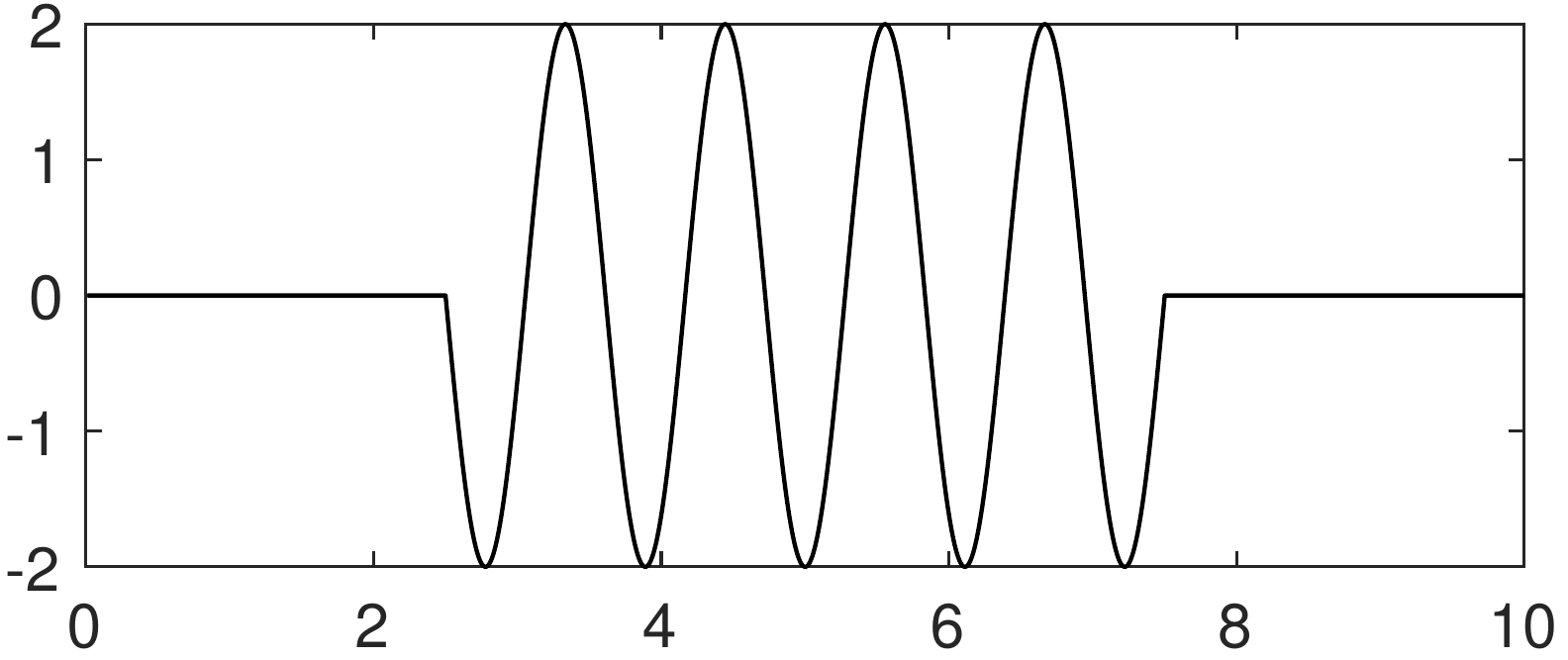}\put(80,32){$\ba_2(t)$}\end{overpic}};
    	\node at (0,-0.9) {\frame{\begin{overpic}[scale=.15]{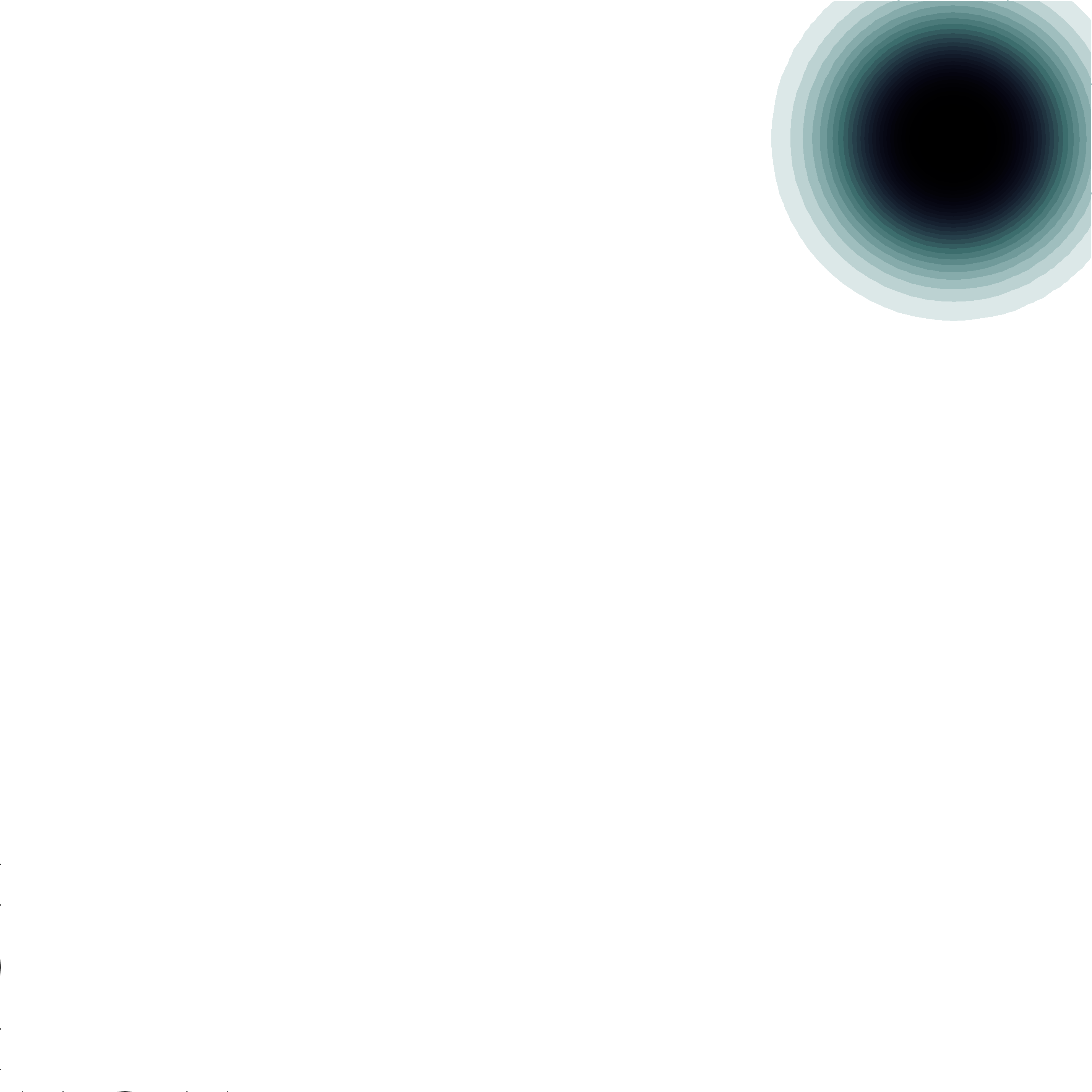}\put(5,82){$\bphi_3$}\end{overpic}}};
    	\node (a3) at (3.4,-1) {\begin{overpic}[scale=0.3]{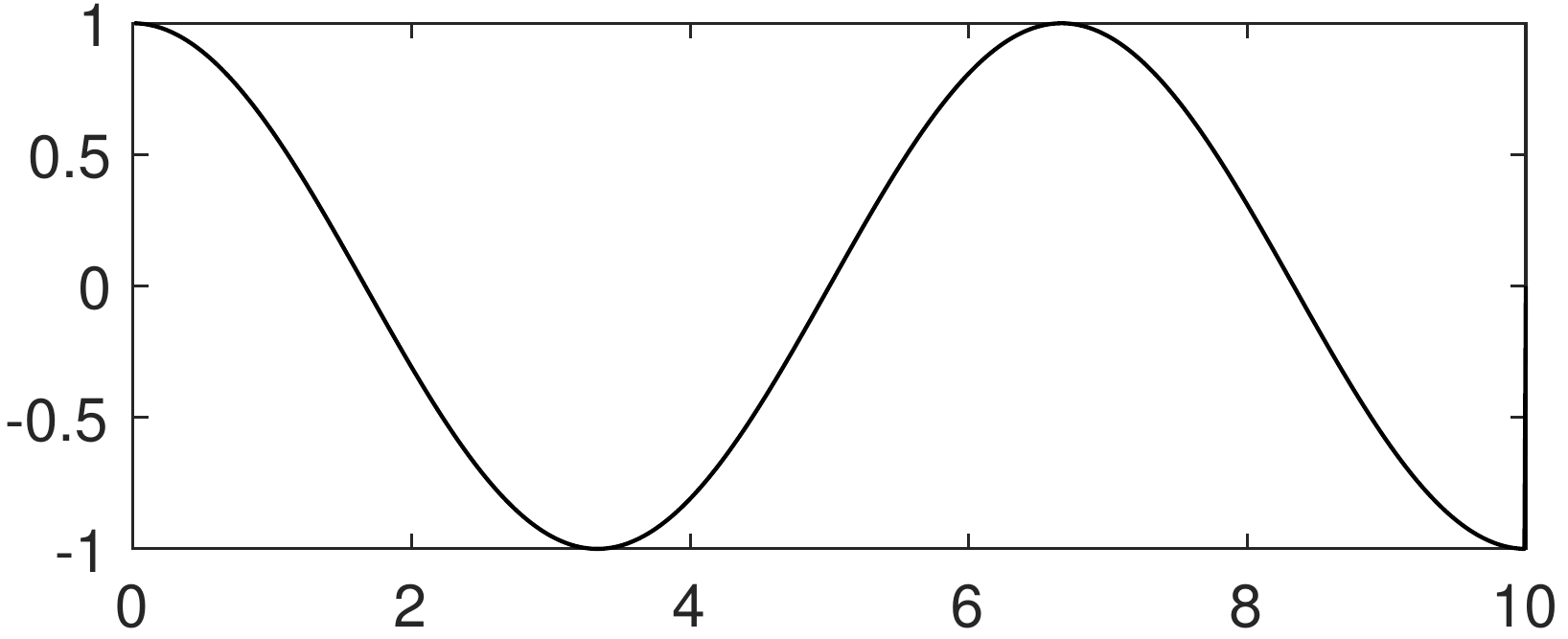}\put(80,32){$\ba_3(t)$}\end{overpic}};
	\node at (map.south) [align=left] {\\ mrDMD amplitude map};
	\node at (a3.south) [align=left] {\\ True dynamics};
	\draw[->,ultra thick] (5.75,4) -- (7.5,4) -- (7.5,3.25);
	\draw[->,ultra thick] (5.75,1.5) -- (7.75,1.5);
	\draw[->,ultra thick] (5.75,-0.75) -- (6.4,-0.75);
	\end{tikzpicture}
	\caption{{\bf Multiscale video}. The mrDMD analysis accurately identifies the active frequencies in each temporal window of the decomposition. The plot on the right colors each time-frequency bin by mrDMD mode amplitude.\label{fig:toy_true}}
\end{figure}
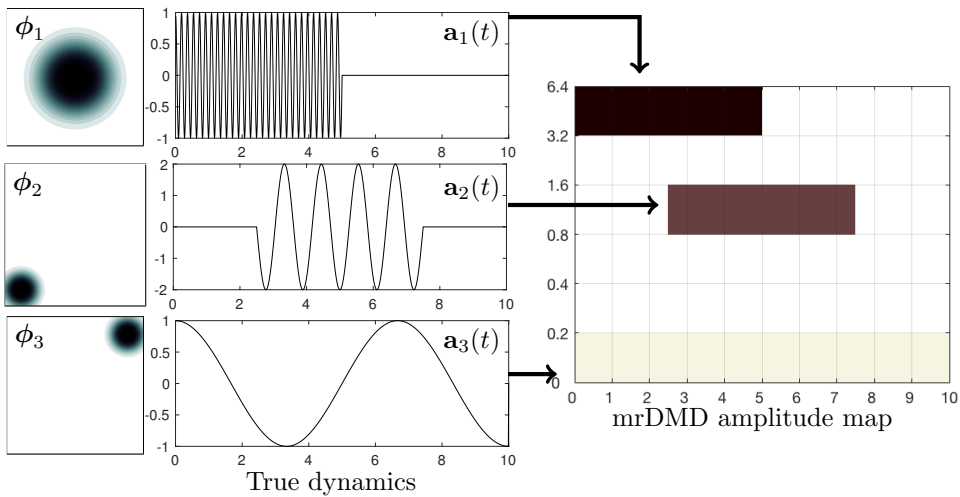
\begin{figure*}[t]
	\centering
	\subfloat[True modes 1-3]{\frame{\includegraphics[width=.15\textwidth]{toy_mode1.png}}~\frame{\includegraphics[width=.15\textwidth]{toy_mode2.png}}~\frame{\includegraphics[width=.15\textwidth]{toy_mode3.png}}
		\label{fig:toy_modes}}
~ \subfloat[mrDMD modes 1-3]{\frame{\includegraphics[width=.15\textwidth]{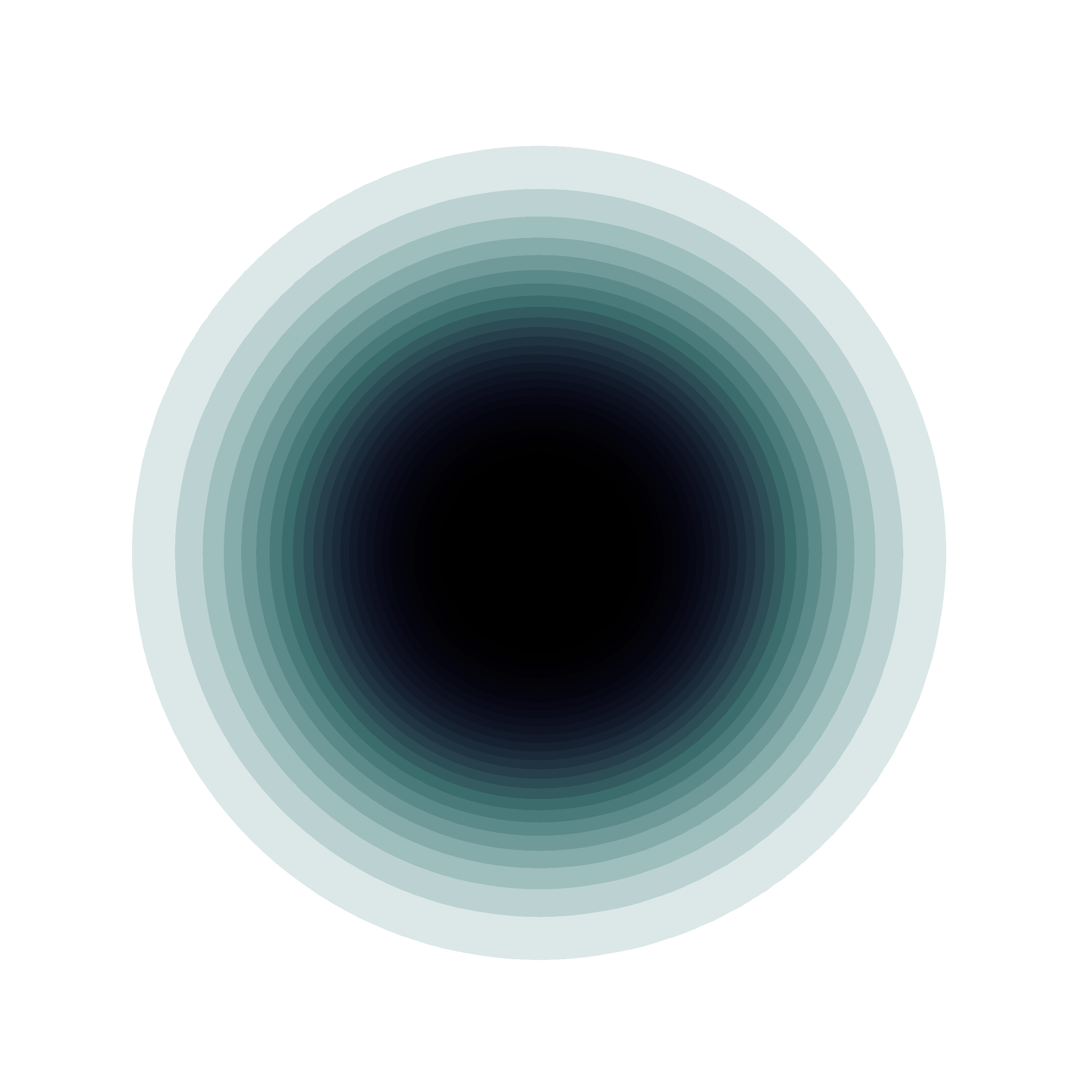}}~\frame{\includegraphics[width=.15\textwidth]{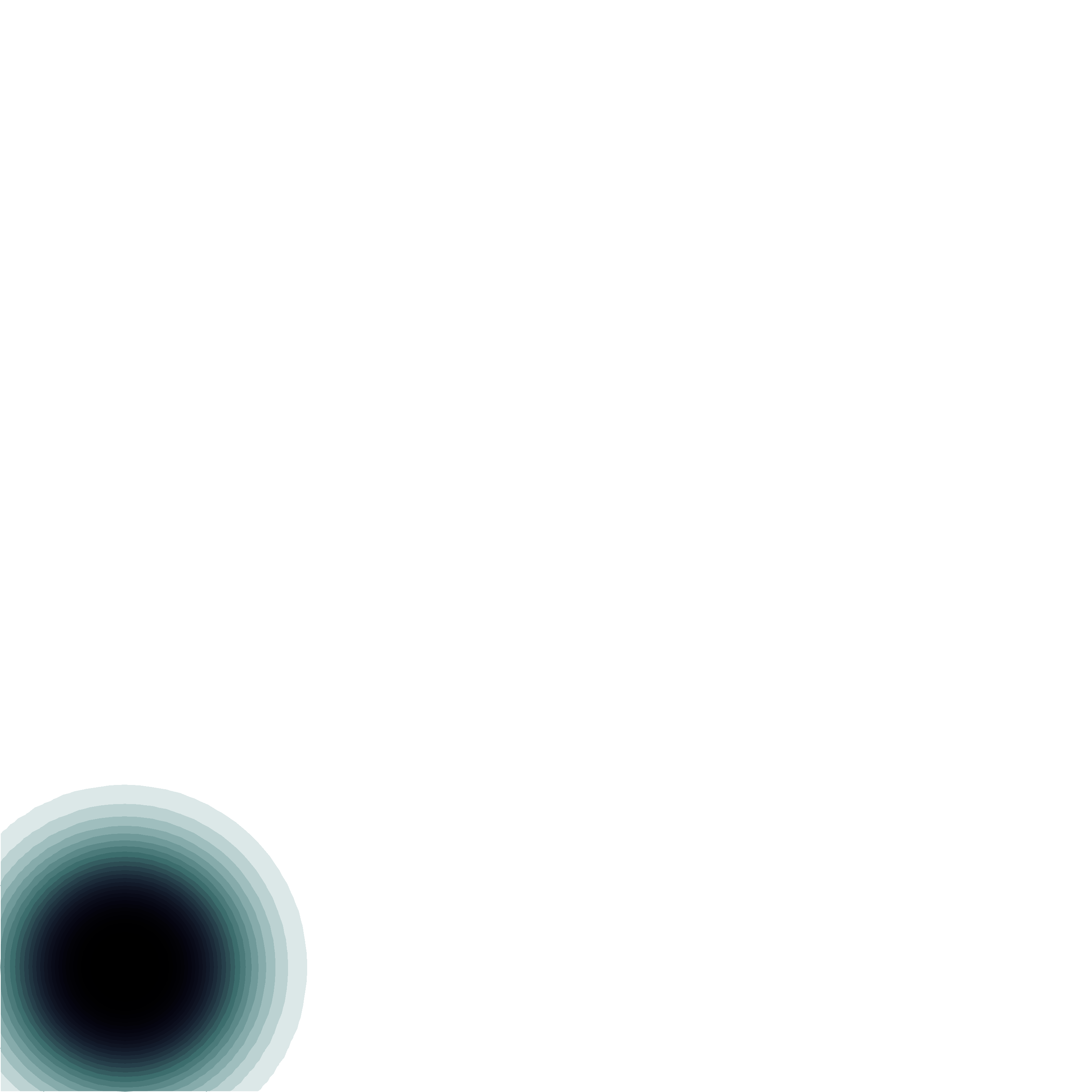}}~\frame{\includegraphics[width=.15\textwidth]{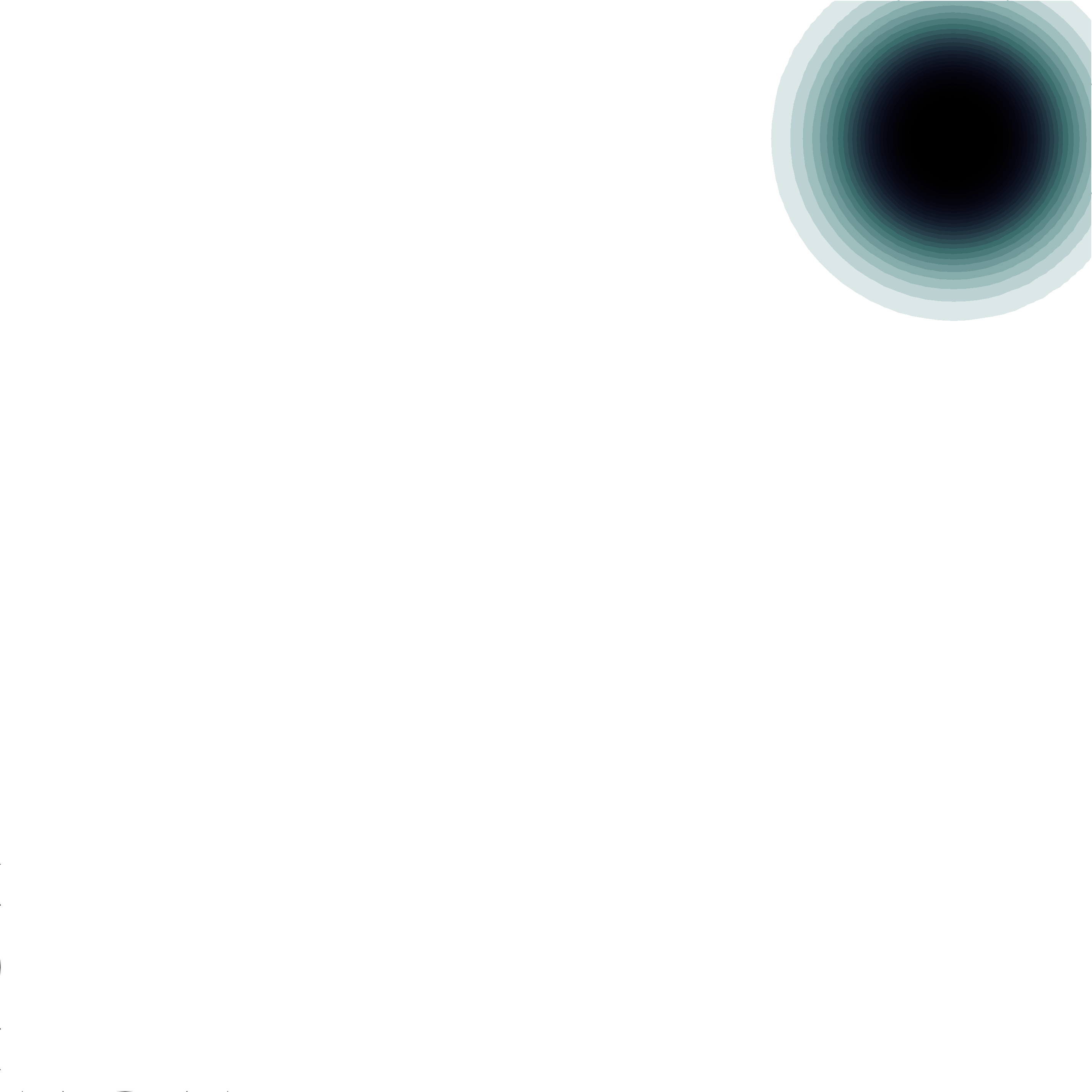}}
		\label{fig:toy_mr_modes}}
\\	\subfloat[DMD modes 1-3]{\frame{\includegraphics[width=.15\textwidth]{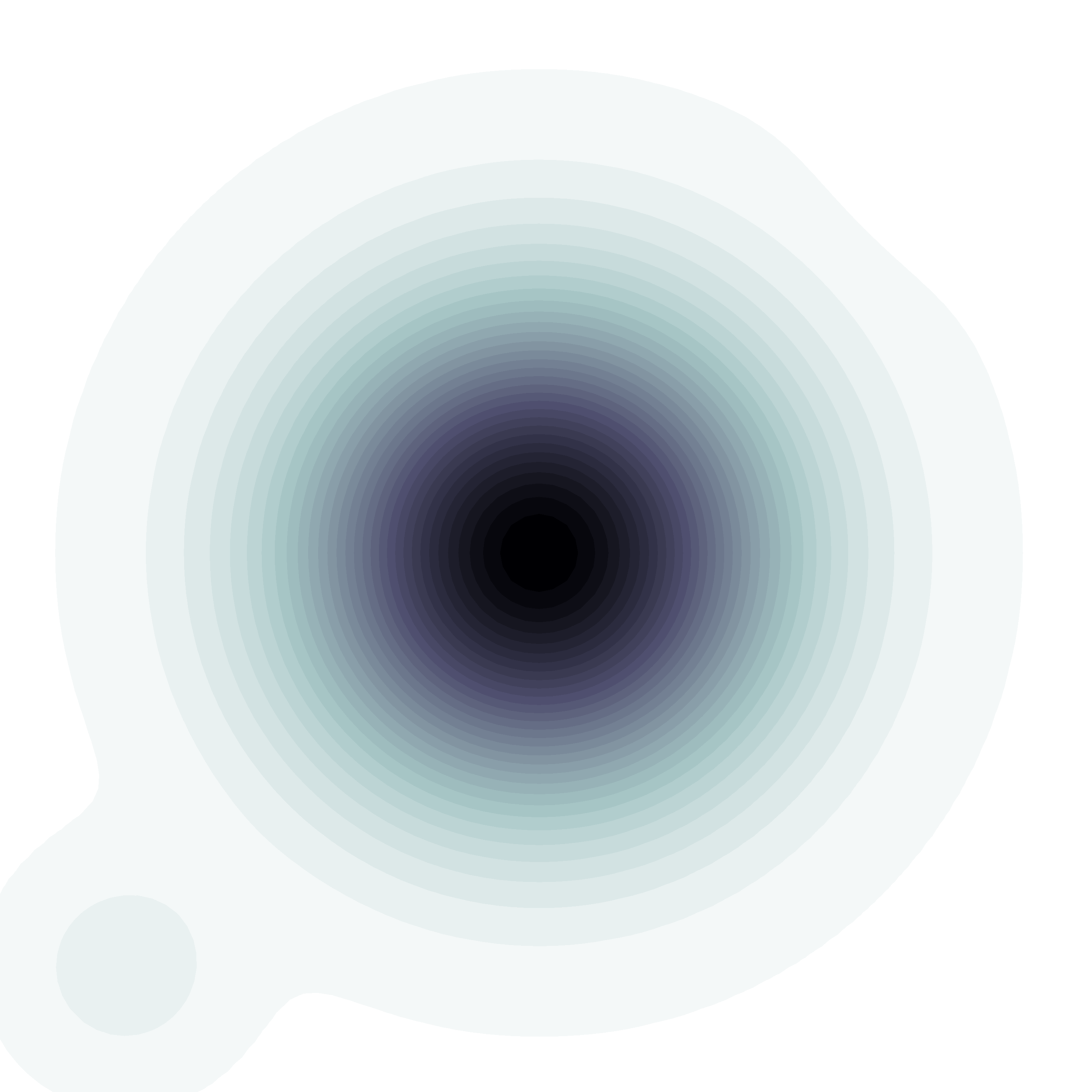}}~\frame{\includegraphics[width=.15\textwidth]{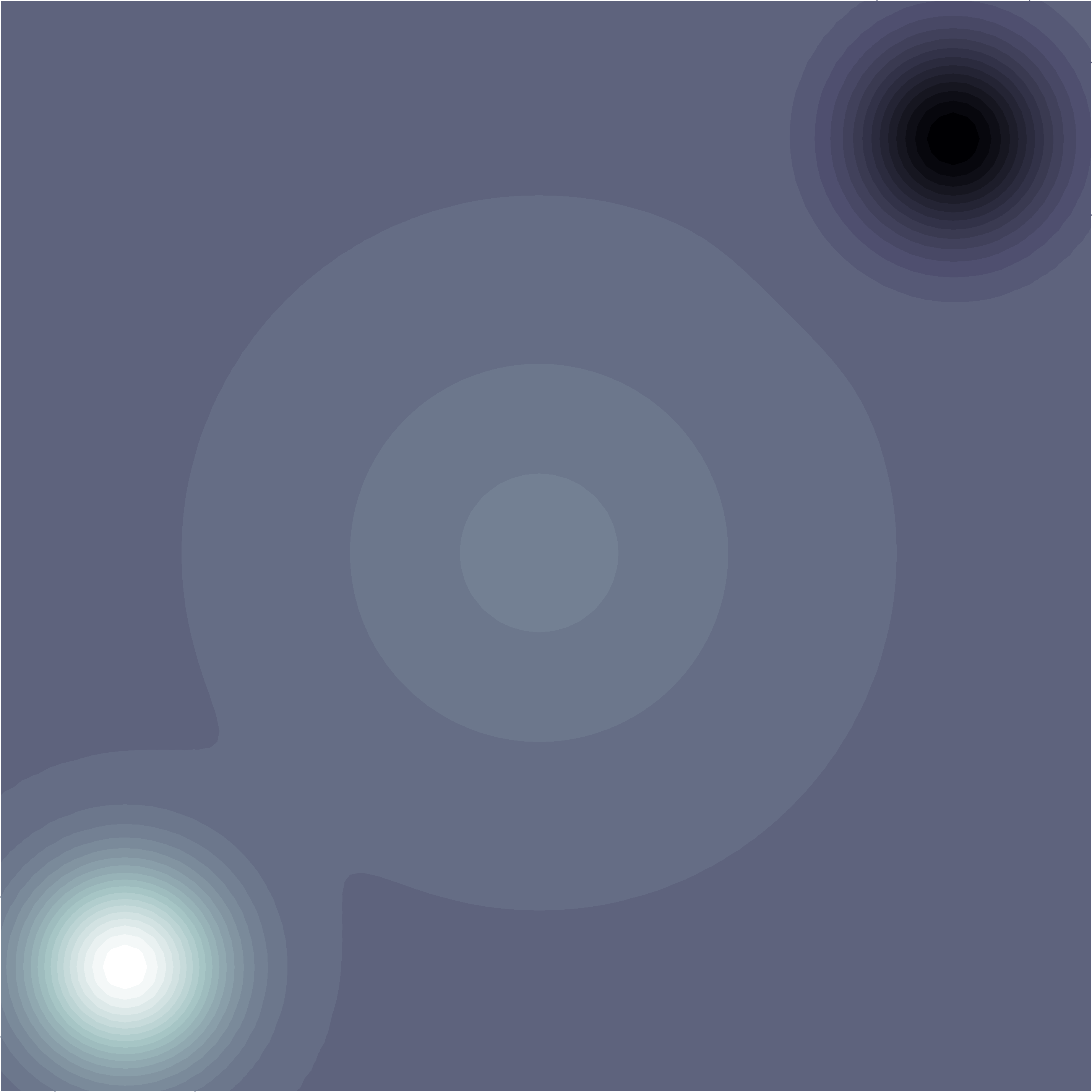}}~\frame{\includegraphics[width=.15\textwidth]{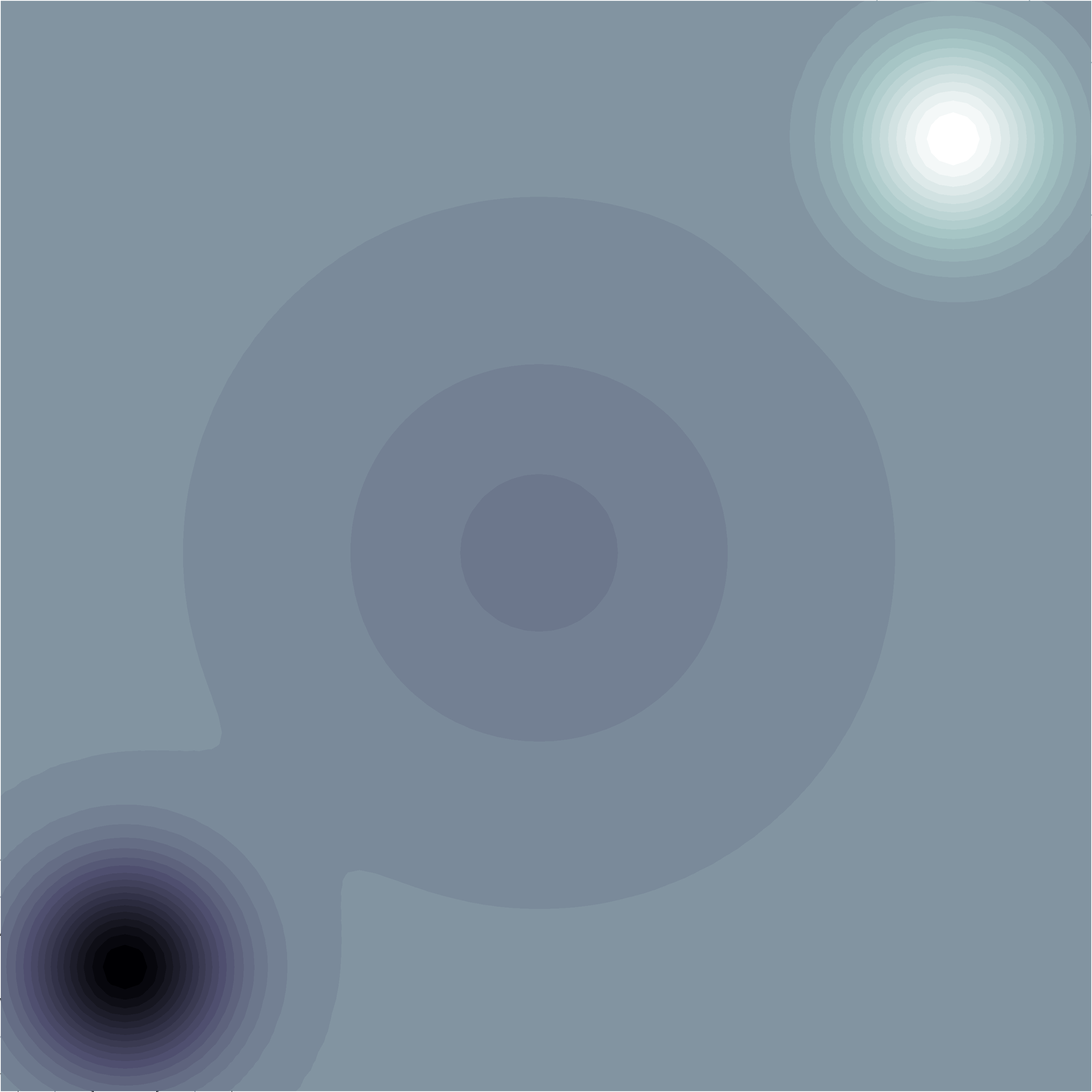}}
		\label{fig:toy_dmd_modes}}
~	\subfloat[POD modes 1-3]{\frame{\includegraphics[width=.15\textwidth]{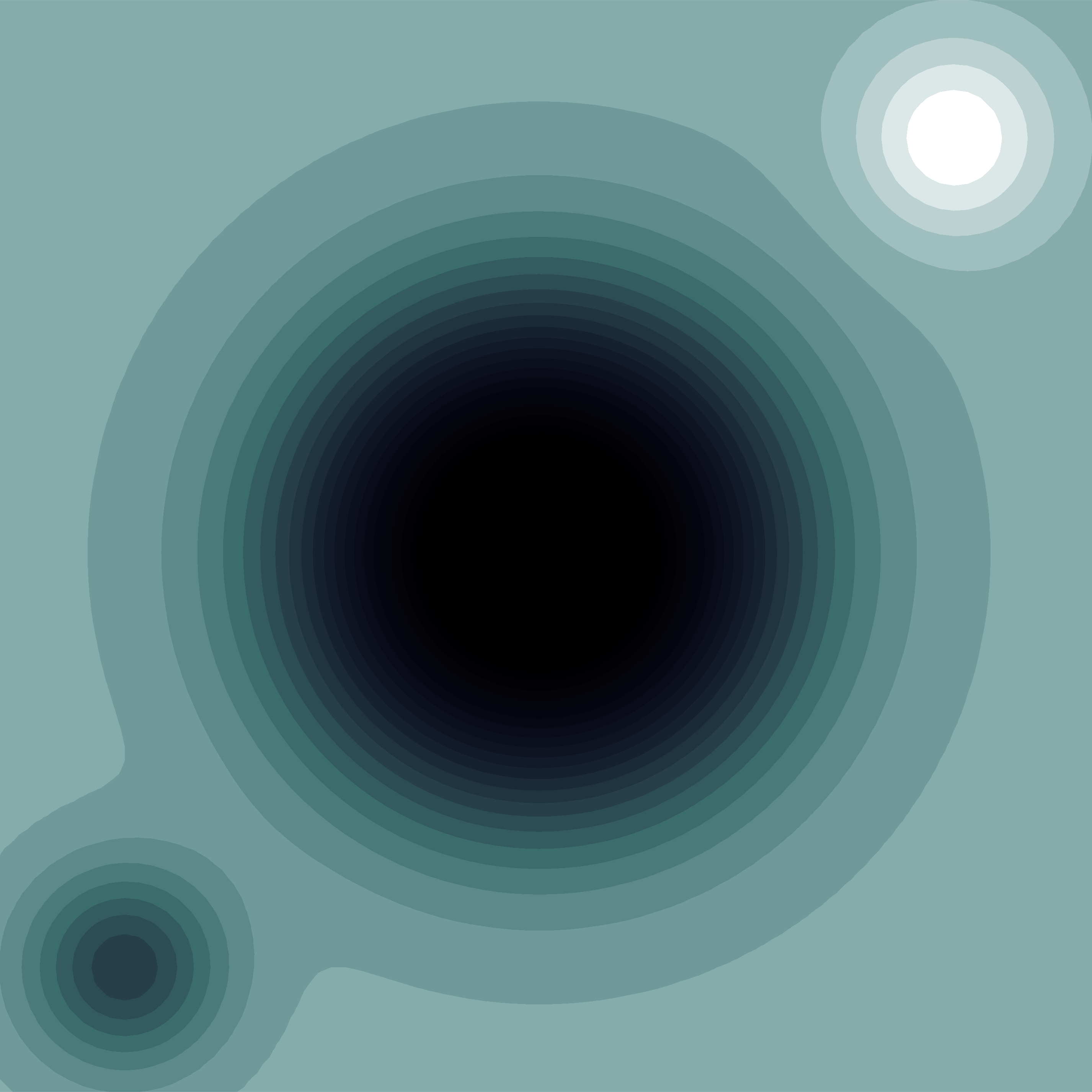}}~\frame{\includegraphics[width=.15\textwidth]{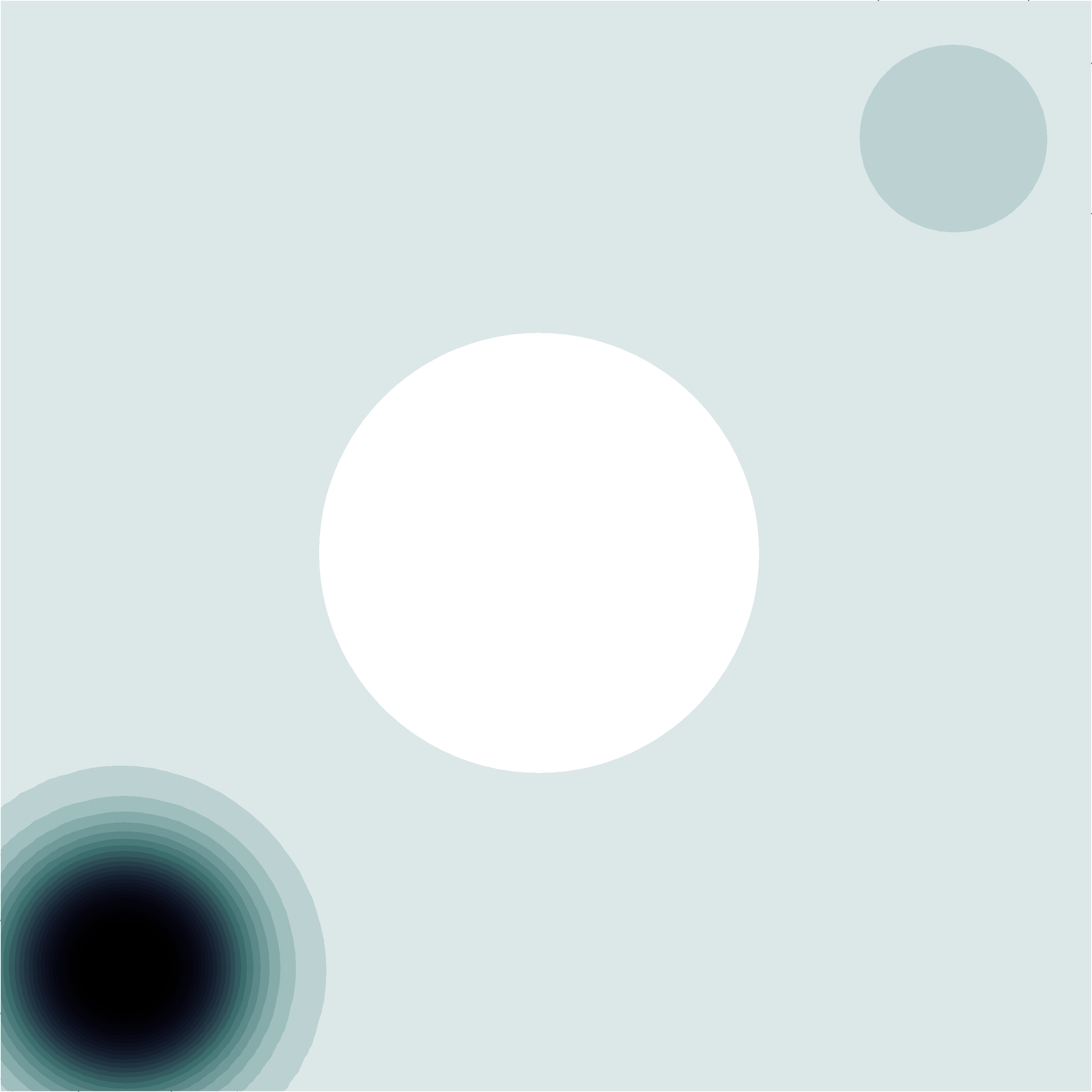}}~\frame{\includegraphics[width=.15\textwidth]{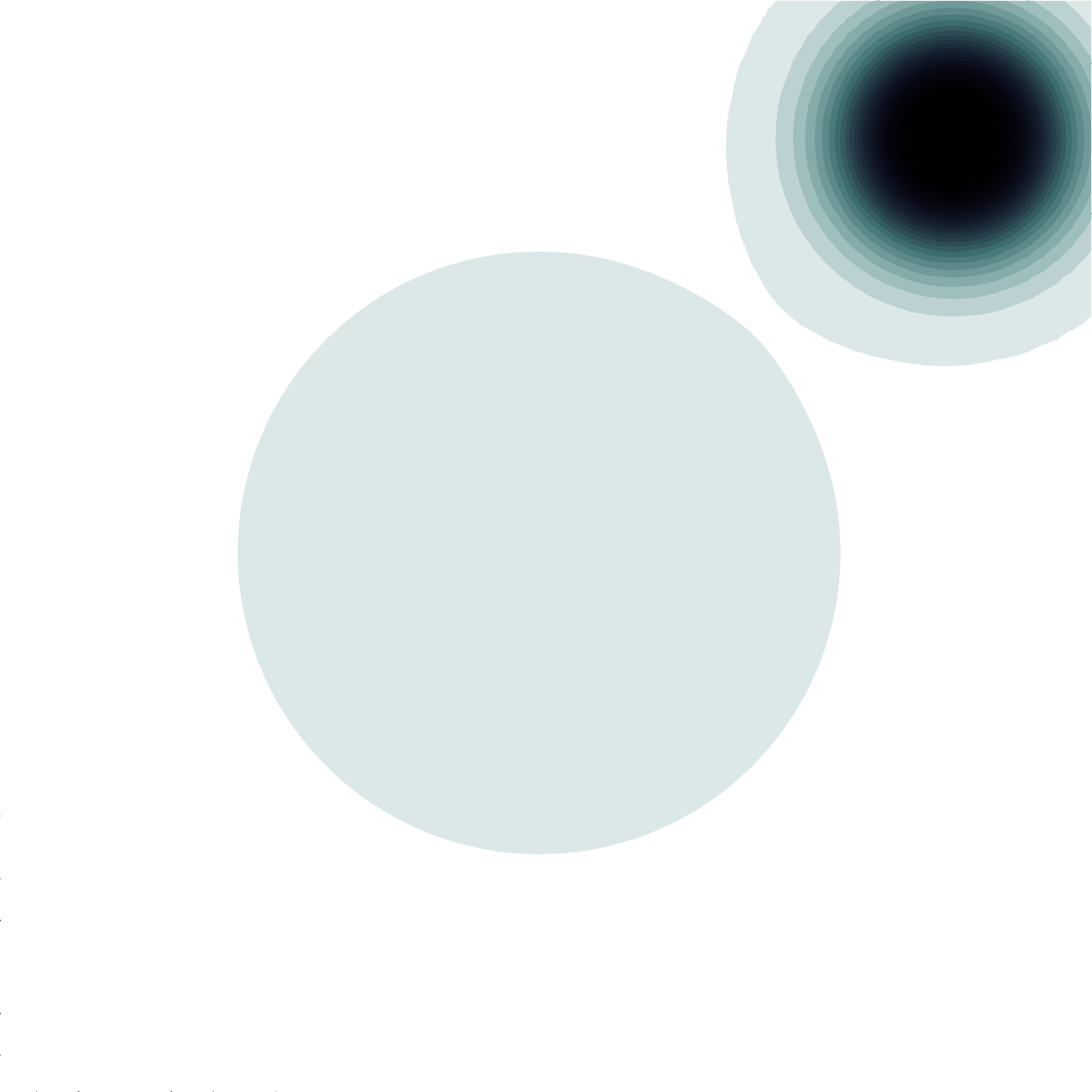}}
		\label{fig:toy_pod_modes}}
		\quad
		\\

\subfloat[QR sensors, left to right: true, mrDMD, DMD, POD]{\frame{\includegraphics[width=.15\textwidth]{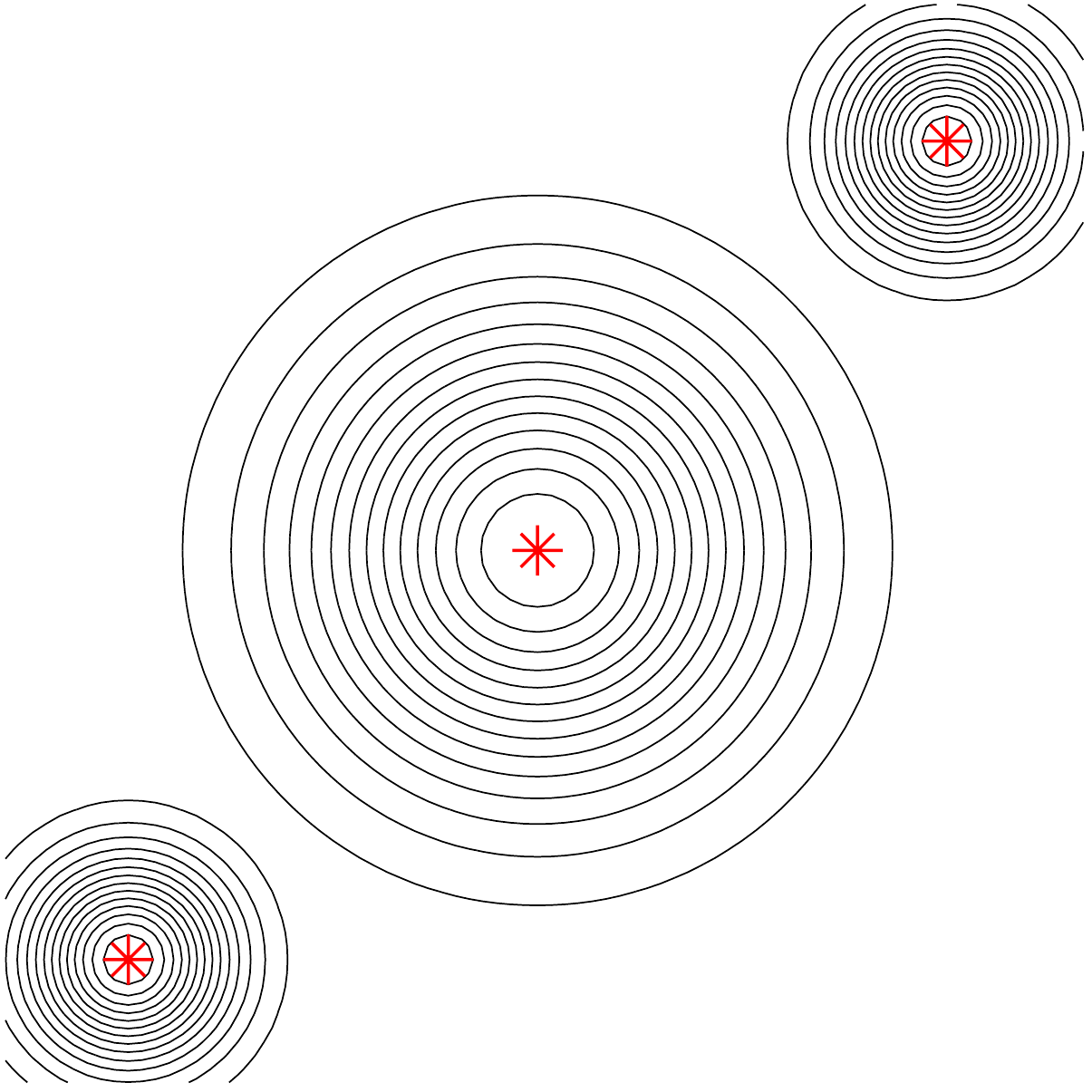}}~\frame{\includegraphics[width=.15\textwidth]{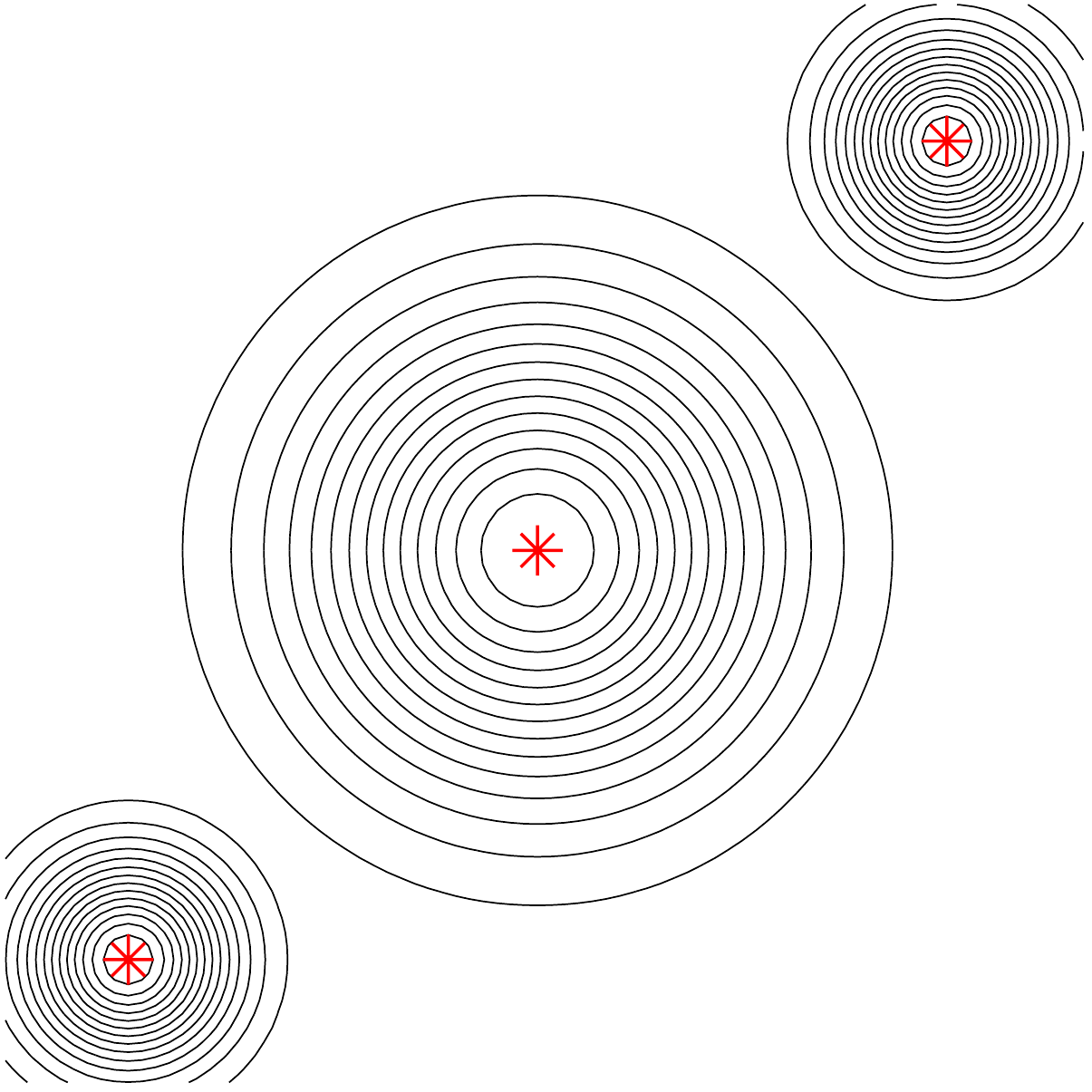}}~\frame{\includegraphics[width=.15\textwidth]{toy_qdmd_sensors.pdf}}~\frame{\includegraphics[width=.15\textwidth]{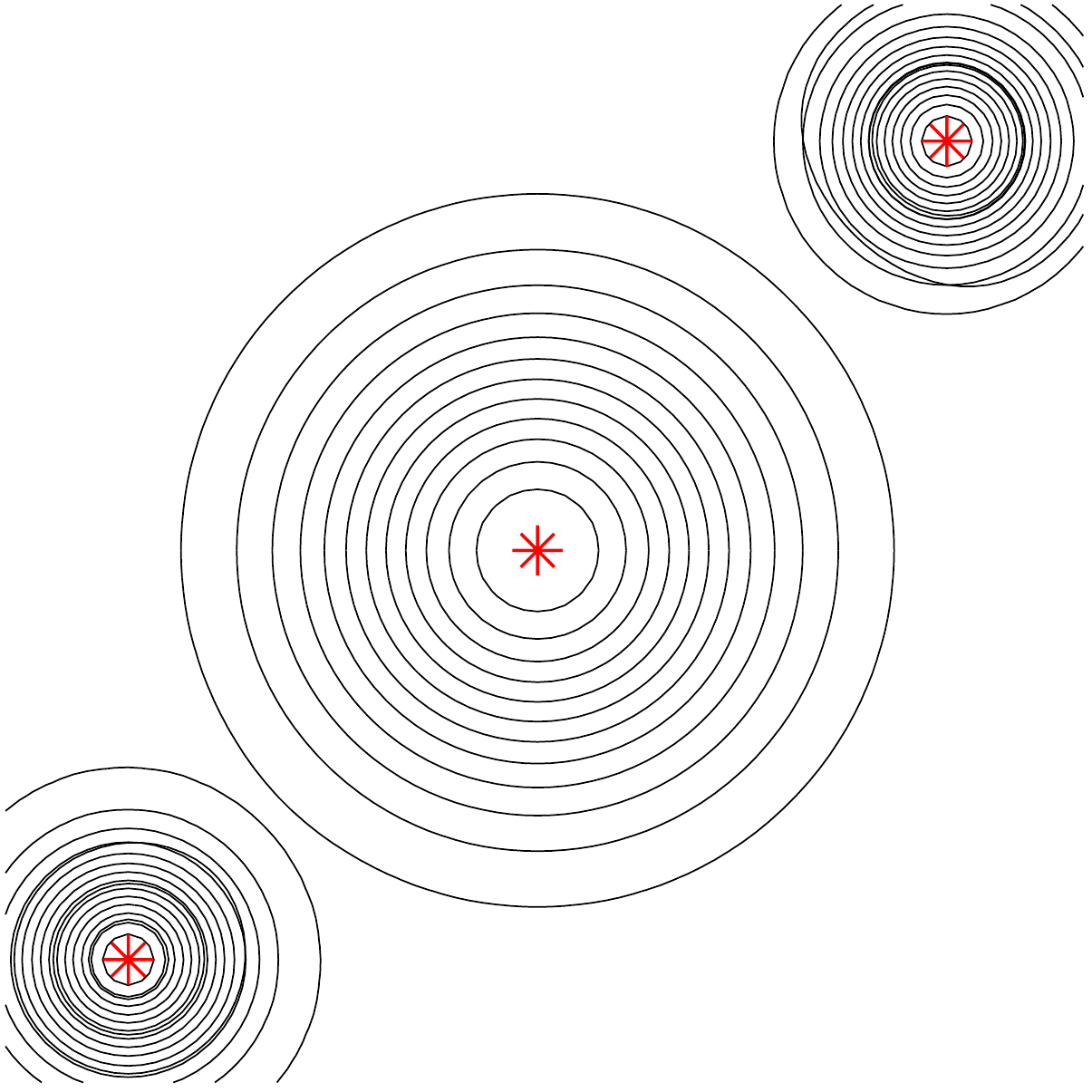}}
		\label{fig:toy_sensors}}
\caption{{\bf Comparison of spatial modes and sensors} shows that mrDMD (b) is ideal for isolating the true windowed dynamics (a), unlike standard DMD (c) or POD (d).\label{fig:toy_comp} Despite varying degrees of spatial mode separation, QR sensors extracted from the different basis modes (e) closely agree.}
\end{figure*}

The identification of spatial modes clearly depends on separation by temporal frequency. \cref{fig:toy_comp} compares mrDMD results with proper orthogonal decomposition (POD), the primary dimensionality reduction tool used in ROMs. The mrDMD time-frequency analysis successfully isolates the true modes of the system, while POD fails at this task since it is a variance-based decomposition in which modes are eigenvectors of data covariances. Hence, the POD spatial modes appear mixed. In addition, mrDMD outputs the correct frequencies associated with each of these modes. We do not give results for standard DMD since it will fail by fitting exponentials across the global time window.
Within the multiresolution analysis, active modes are repeated across the time-frequency bins, so we include only the mode with highest amplitude. These selected modes describe the feature selection or formation of the index set $\Sub$ used to train sensors. Similarly, we select the first three POD modes, which explain 99.9\% of the variance in the training data, as quantified by the POD eigenvalues.
Interestingly, both decompositions yield the same optimal samples for this example. This is because the three point sources of diffusion correspond to spatial extrema, which are indicated by the overlapping contour plots of the first three true, mrDMD and POD modes respectively. 
\begin{figure}[t!]
	\centering
	\begin{tikzpicture}		
	\node at (0,4) {\begin{overpic}[scale=.4]{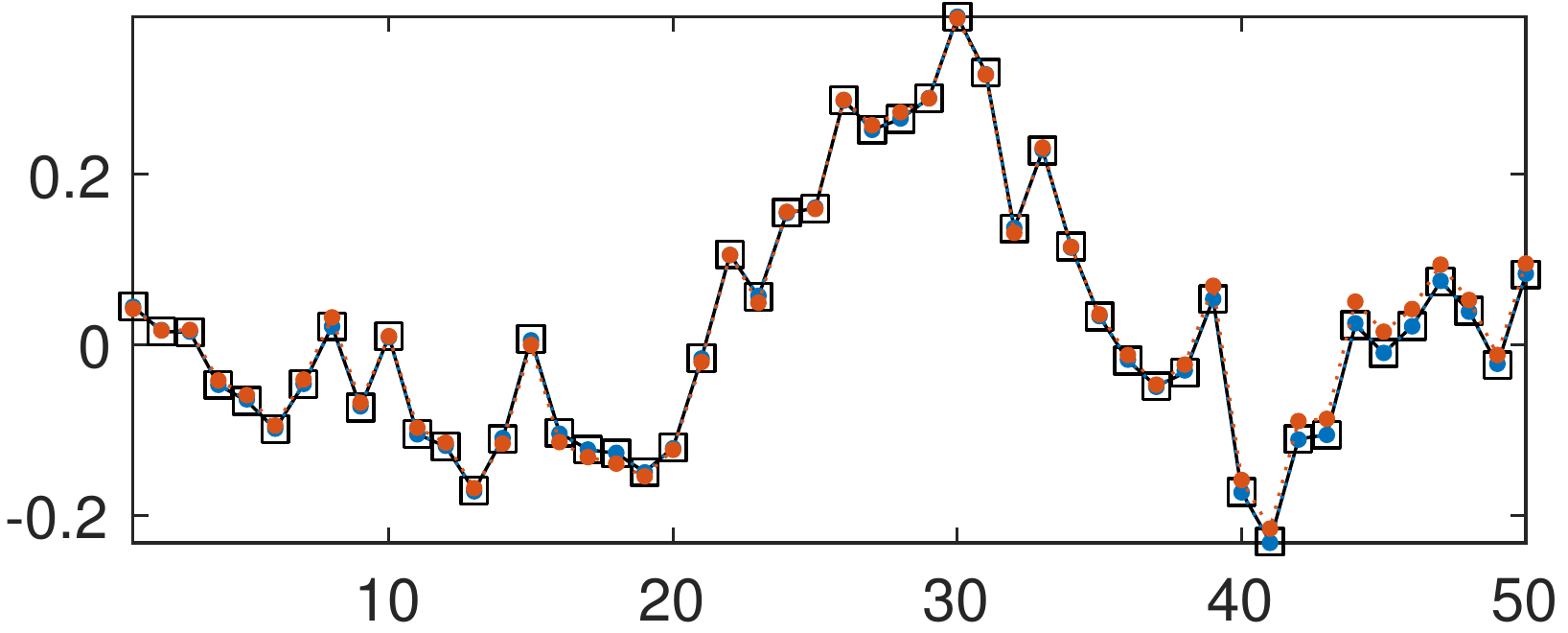}\put(-8,20){$\hat\ba_1(t)$}\end{overpic}};
	\node at (2,2) {\begin{overpic}[scale=.4]{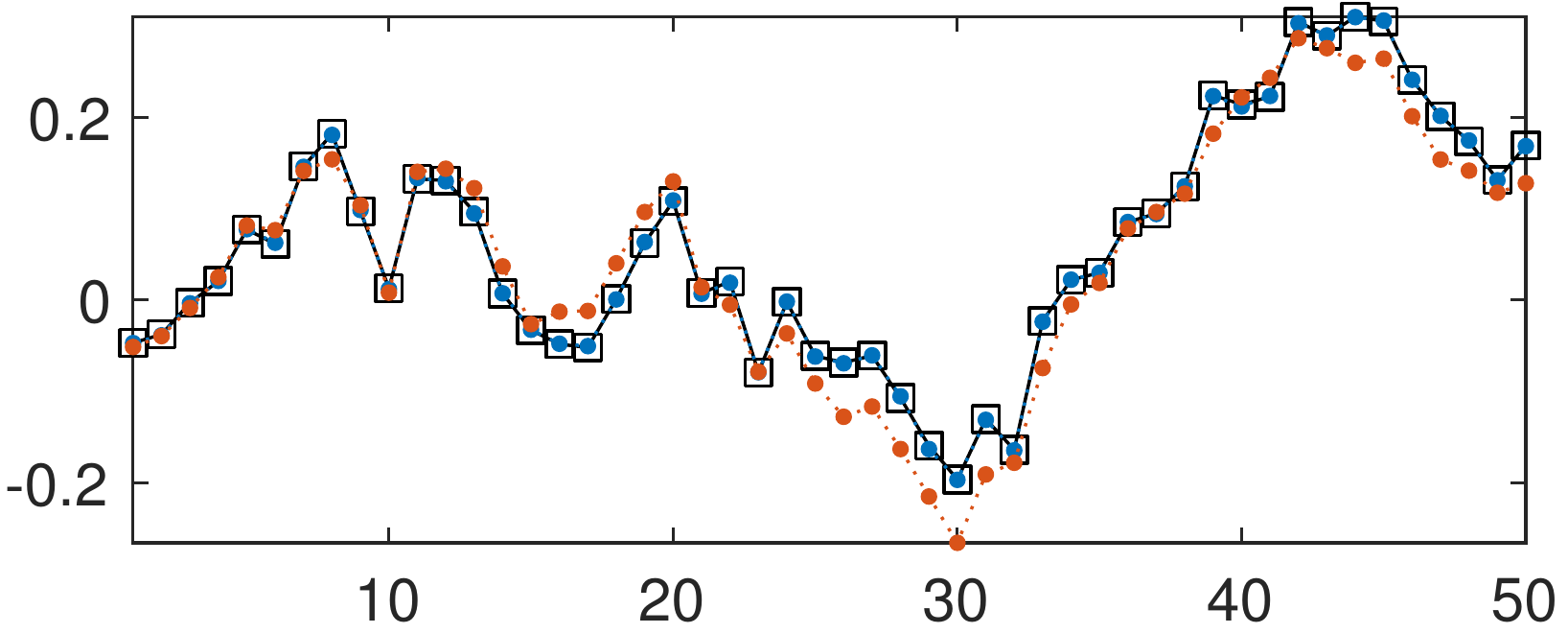}\put(-8,20){$\hat\ba_2(t)$}\end{overpic}};
	\node at (4,0) {\begin{overpic}[scale=.4]{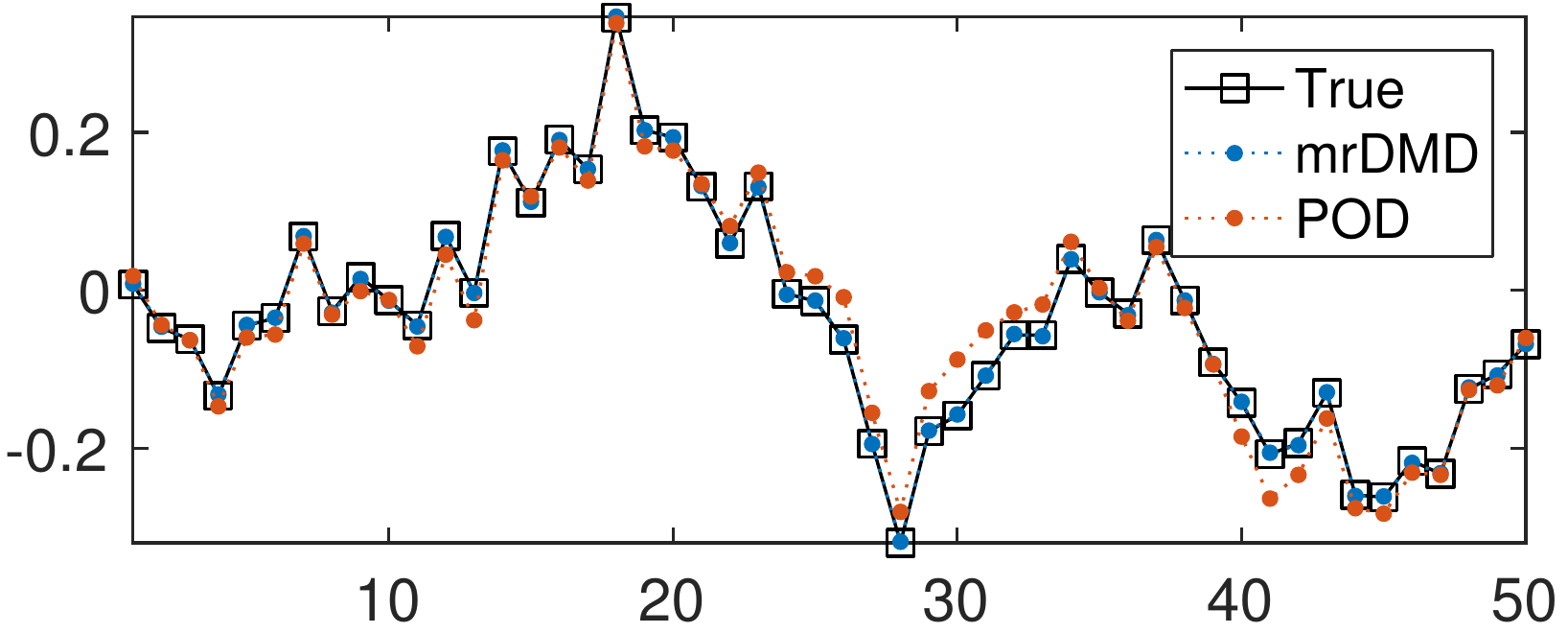}\put(-8,20){$\hat\ba_3(t)$}\put(47,-2){time}\end{overpic}};
	\end{tikzpicture}
\caption{{\bf Comparison of estimated coefficients}. \label{fig:toy_time} The estimation of time coefficients from mrDMD modes and QR-optimized sensors is more accurate than with POD modes. True coefficients in this test window are generated by a Gaussian random process.}
\end{figure}
\begin{figure*}
	\centering
	\begin{overpic}[width=.8\textwidth]{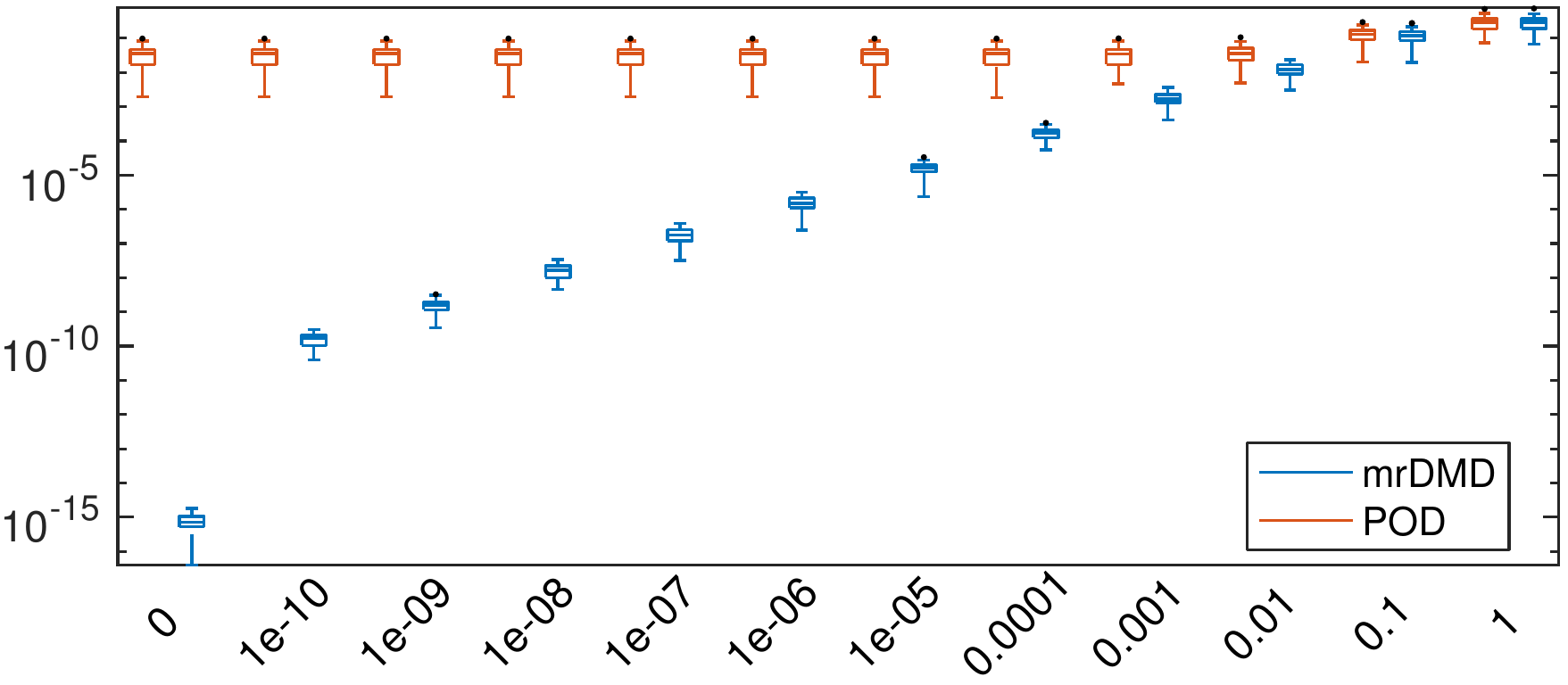}
	\put(-5,20){\rotatebox{90}{$l_2$ error}}
	\put(40,-2.5){Noise variance $\sigma$}
	\end{overpic}
	\caption{{\bf Noise corrupted measurements}. Approximation error with mrDMD modes and sensors trends several orders of magnitude lower than with POD modes and sensors for small levels of sensor white noise. \label{fig:toy_noise_study}}
\end{figure*}

Even though the same optimal samples are identified, feature selection greatly affects the quality of reconstruction for future state prediction. To test this, we generate test data using time coefficients for each mode from a  random process in order to estimate the true coefficients from the three optimal sampling locations. Further, white noise of variance $\sigma^2=0.001$ is added to the three sampled observations. Finally, coefficients estimated from the subselected basis $\bPhi_{\Sub}$ consisting of mrDMD modes~(\cref{fig:toy_mr_modes}) are compared to a basis of POD modes~(\cref{fig:toy_pod_modes}). Here least-squares estimation is used since the problem is well-posed (three equations and three unknowns). Results are plotted in~\cref{fig:toy_time}, where it can be seen that mrDMD-based estimation (blue) best approximates the true random process (black). By contrast, POD-based estimation (red) deviates from the true coefficients, particularly in the second and third modal coefficients.

We quantify this further in a noise study whose results are plotted in~\cref{fig:toy_noise_study}. White noise levels $\sigma^2$ are increased on a logarithmic scale to examine the $l_2$ reconstruction error trend of the coefficient vectors. Estimation with mrDMD points to an exponential increase in error until the signal appears saturated by noise at $\sigma^2=1$. POD-based estimation, however, appears to saturate much earlier at small noise levels, and is thus not as suitable for multiscale estimation and prediction. Specifically, the full state can now be reconstructed from noisy point observations in any validation window using estimated coefficients~\cref{eqn:full_state}.

In this video example, the true spatiotemporal modes were selected based on three different identified frequencies. For estimation, we leverage these three modes and associated sensors to reconstruct unseen random dynamics at future time instances. In many systems, it is not immediately obvious which of the identified dynamics are active in a new time window. This lack of knowledge of true dynamics typifies many real-world systems, in which nonlinear, chaotic dynamics and regime switching limit our ability to predict future states. However, we can harness observations from optimally sampled observations (sensors) to classify active dynamics in new time windows, which is illustrated here for sparse identification of El Ni\~no warming periods from global ocean temperature data.

\subsection{NOAA ocean surface temperature}

We apply multiresolution analysis to real satellite data of weekly ocean surface temperatures from 1990-2017. Here, the measure of success is accurate identification of intermittent warming and cooling events that are famously implicated in global weather patterns and climate change - El Ni\~no and La Ni\~na events. Weekly temperature means are collected on the entire 360x180 spatial grid, with data omitted over continents and land. The dimension of each weekly snapshot then reduces from $360\times 180=64800$ to $n=44219$ spatial gridpoints.  We run mrDMD in training windows of 16 years (multiples of 2 facilitate clean separation of annual scales) up to four decomposition levels so that the fastest identifiable frequency is biennial, effectively discarding the dominant annual scale from the analysis. This is done to parse out the El Ni\~no Southern Oscillation (ENSO), defined as any sustained temperature anomaly above running mean temperature with a duration of 9 to 24 months.

\begin{figure*}[tbhp]
	\centering
	\begin{tikzpicture}
	\node[inner sep=0pt] (amps) at (0,0)
	{\begin{overpic}[width=.475\textwidth]{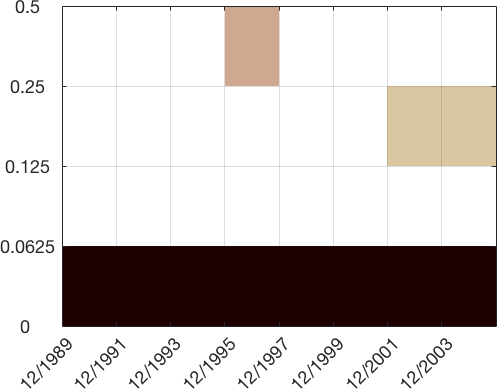}\put(-6,20){\rotatebox{90}{yearly frequency ($\omega$)}}\end{overpic}};
	\node[inner sep=0pt,draw=black,thick] (enso) at (0.1,4)
	{\includegraphics[width=.15\textwidth]{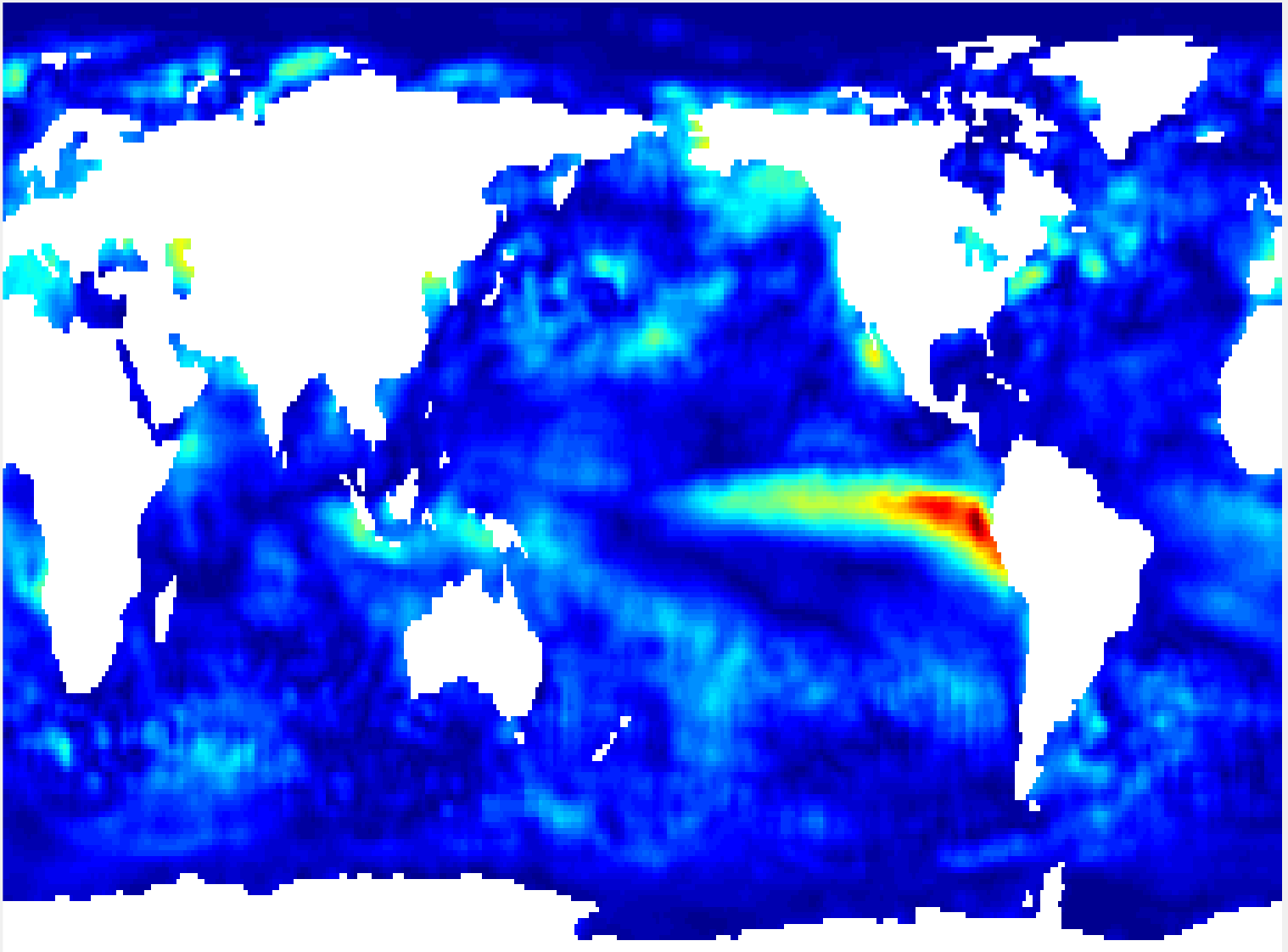}};
	\node[inner sep=0pt,draw=black,thick] (bso) at (2.3,4)
	{\includegraphics[width=.15\textwidth]{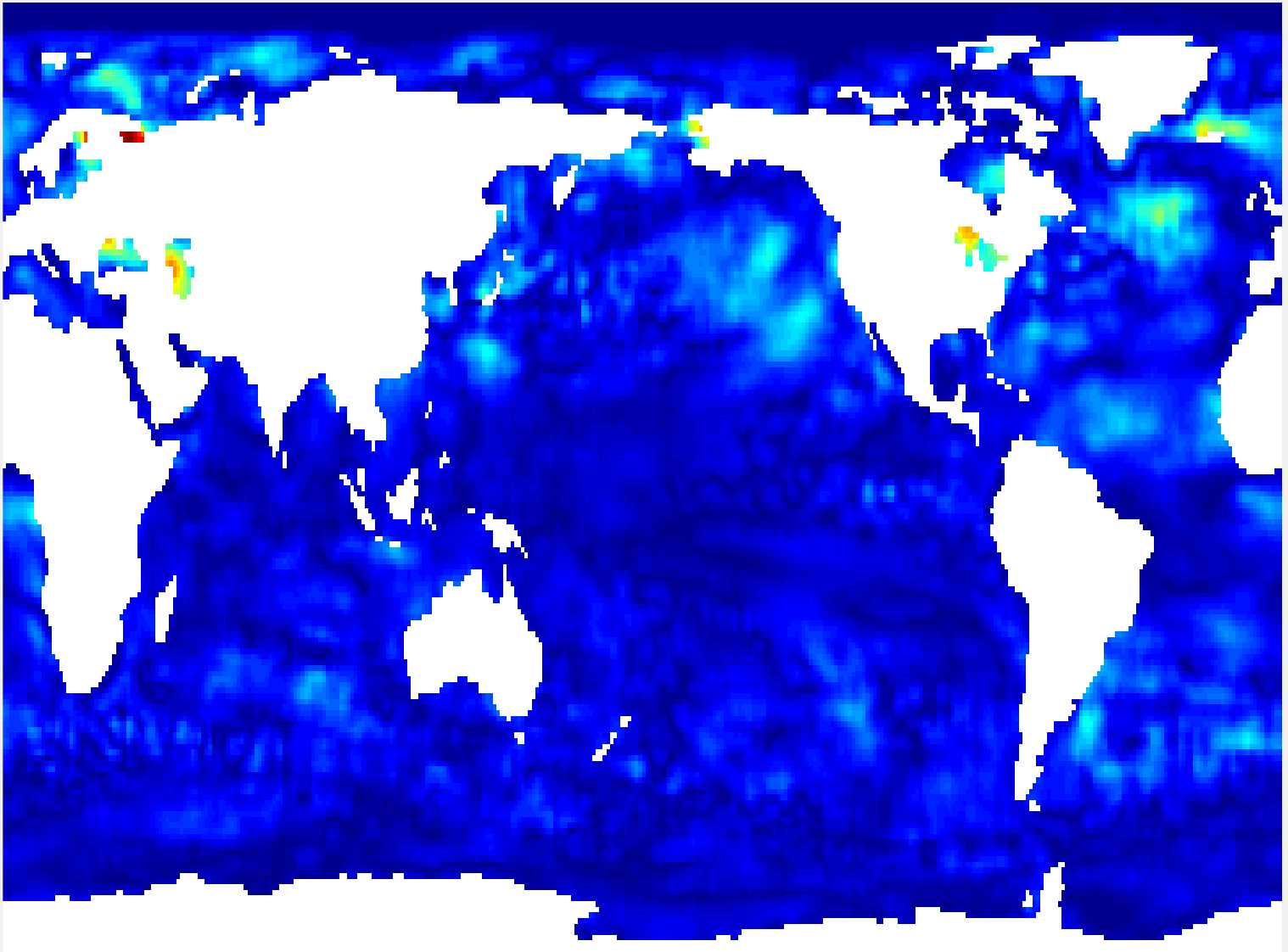}};
	\node[inner sep=0pt,draw=black,thick] (bg) at (.3,-1)
	{\includegraphics[width=.15\textwidth]{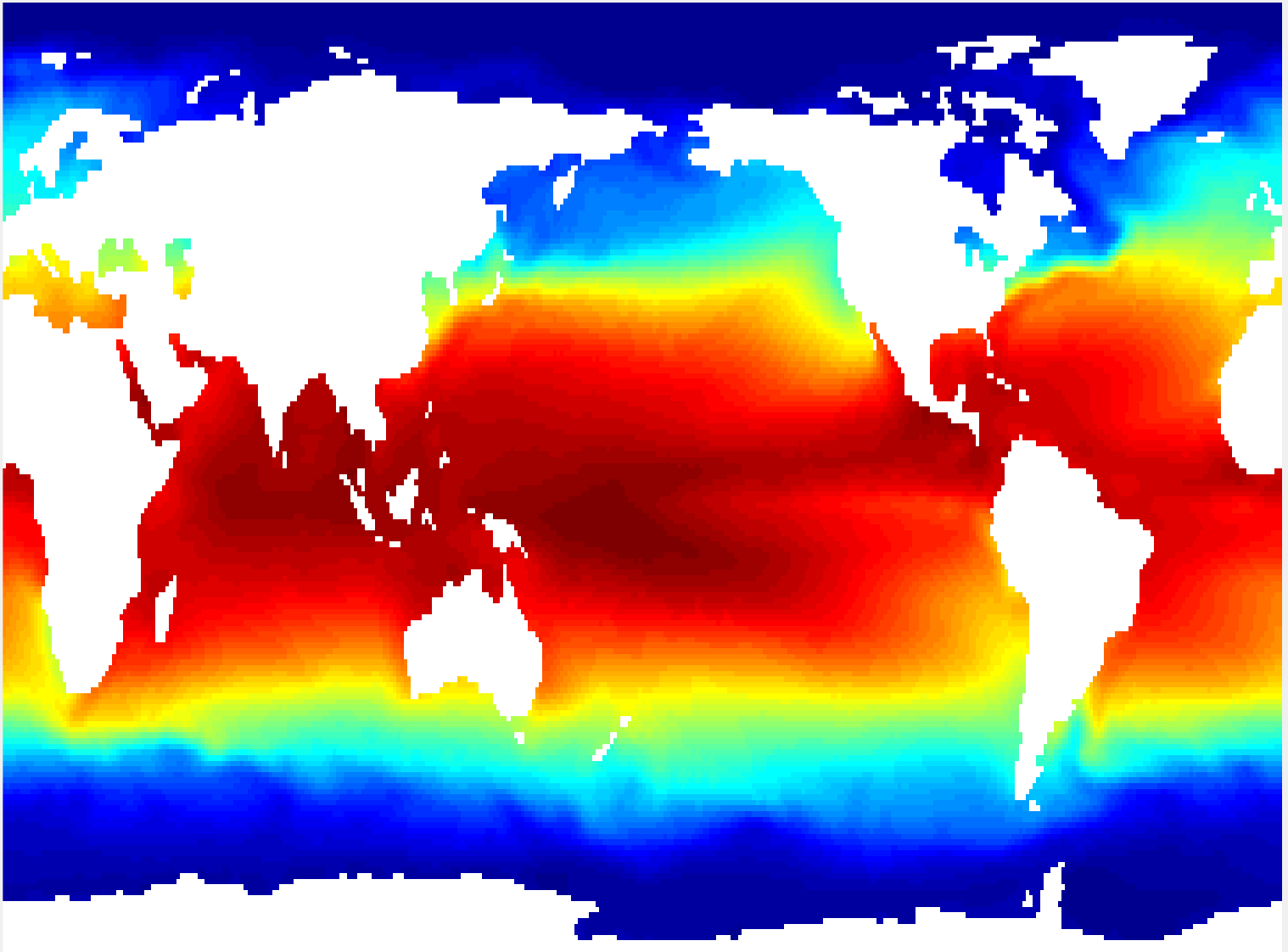}};
	\draw[->,thick] (enso.south) -- (0.1,2.5)
	node[near end,above] {ENSO};
	\draw[->,thick] (bso.south) -- (2.3,1)
	node[near end,above] {BSO};
	\end{tikzpicture}
	\begin{tikzpicture}
	\node[inner sep=0pt] (amps) at (0,0)
	{\includegraphics[width=.475\textwidth]{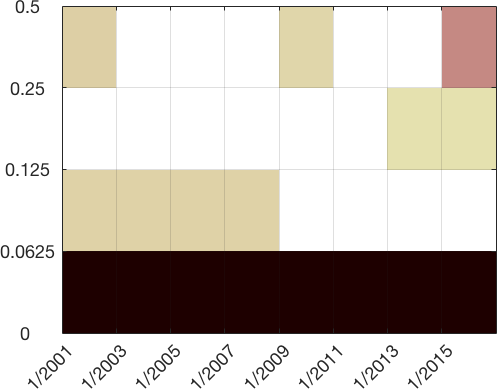}};
	\node[inner sep=0pt,draw=black,thick] (wkenso1) at (-1.9,4)
	{\includegraphics[width=.15\textwidth]{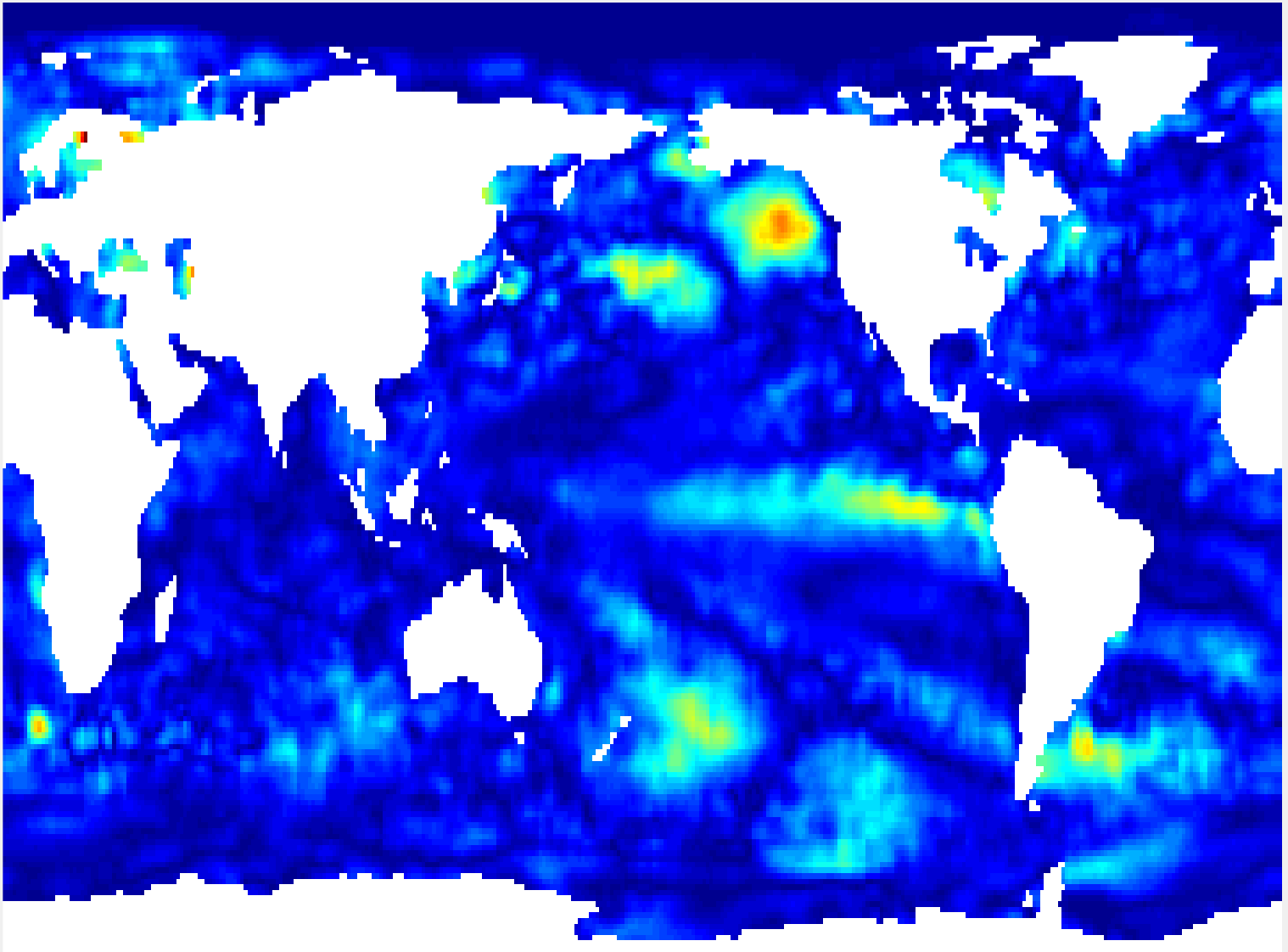}};
	\node[inner sep=0pt,draw=black,thick] (wkenso) at (.3,4)
	{\includegraphics[width=.15\textwidth]{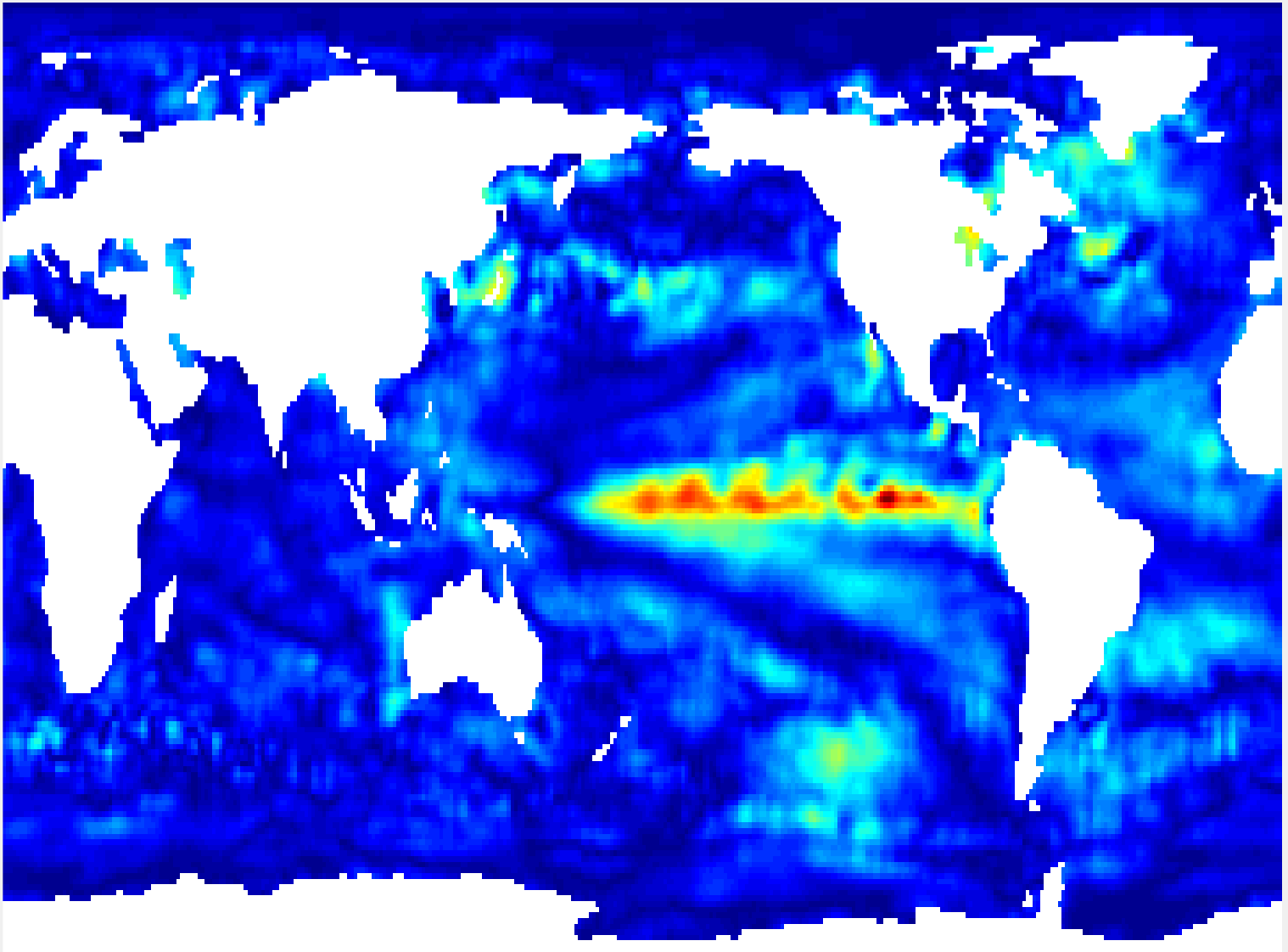}};
	\node[inner sep=0pt,draw=black,thick] (enso) at (2.3,4)
	{\includegraphics[width=.15\textwidth]{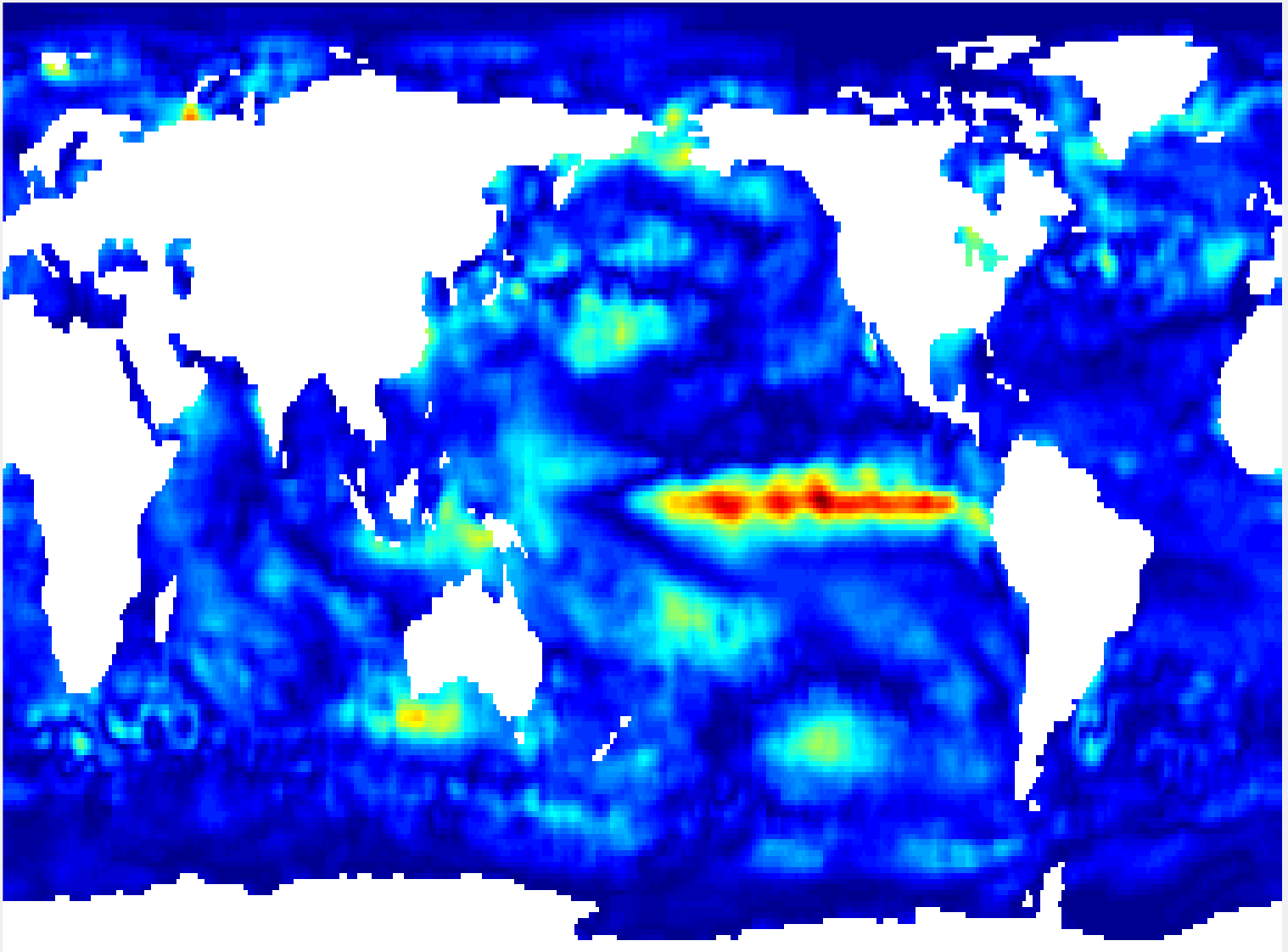}};
	\node[inner sep=0pt,draw=black,thick] (bg) at (.3,-1)
	{\includegraphics[width=.15\textwidth]{Phi_sst6}};
	\draw[->,thick] (wkenso1.south) -- (-1.9,2.5)
	node[near end,above] {weak ENSO};
	\draw[->,thick] (wkenso.south) -- (.5,2.5)
	node[near end,above] {weak ENSO};
	\draw[->,thick] (enso.south) -- (2.5,2.5)
	node[near end,above] {ENSO};
	\end{tikzpicture}	
	\caption{{\bf MrDMD map of modal amplitudes by level and time window}. The method captures several dynamically significant oscillations occurring slower than annual dynamics. The identified El Ni\~no (ENSO) and Barents Sea oscillations (BSO) are both rare events known to exert significant influence on global weather patterns.\label{fig:mrdmd_map_sst}}
\end{figure*}

\Cref{fig:mrdmd_map_sst} visualizes the resulting modal amplitude maps in the time-frequency domain, along with oscillatory dynamical modes identified by the decomposition. Each bin is colored by the average modal amplitude of identified dynamics. An energetic background mode is identified in both training windows (1990-2006, 2001-2016) corresponding to the mean temperature distribution. More importantly, lower-energy dynamics identified at the biennial scale correspond closely to well-documented ENSO events, the strongest of which occur in 1997-1999 and 2014-2016. Visual inspection clearly identifies these as ENSO events which cause a characteristic band of warm water across the South Pacific near coastal Peru.

\begin{figure}[t]
	\centering
	\setlength\fboxsep{0pt}
	\subfloat[mrDMD frequencies $\omega_k$]{\begin{overpic}[width=.32\textwidth]{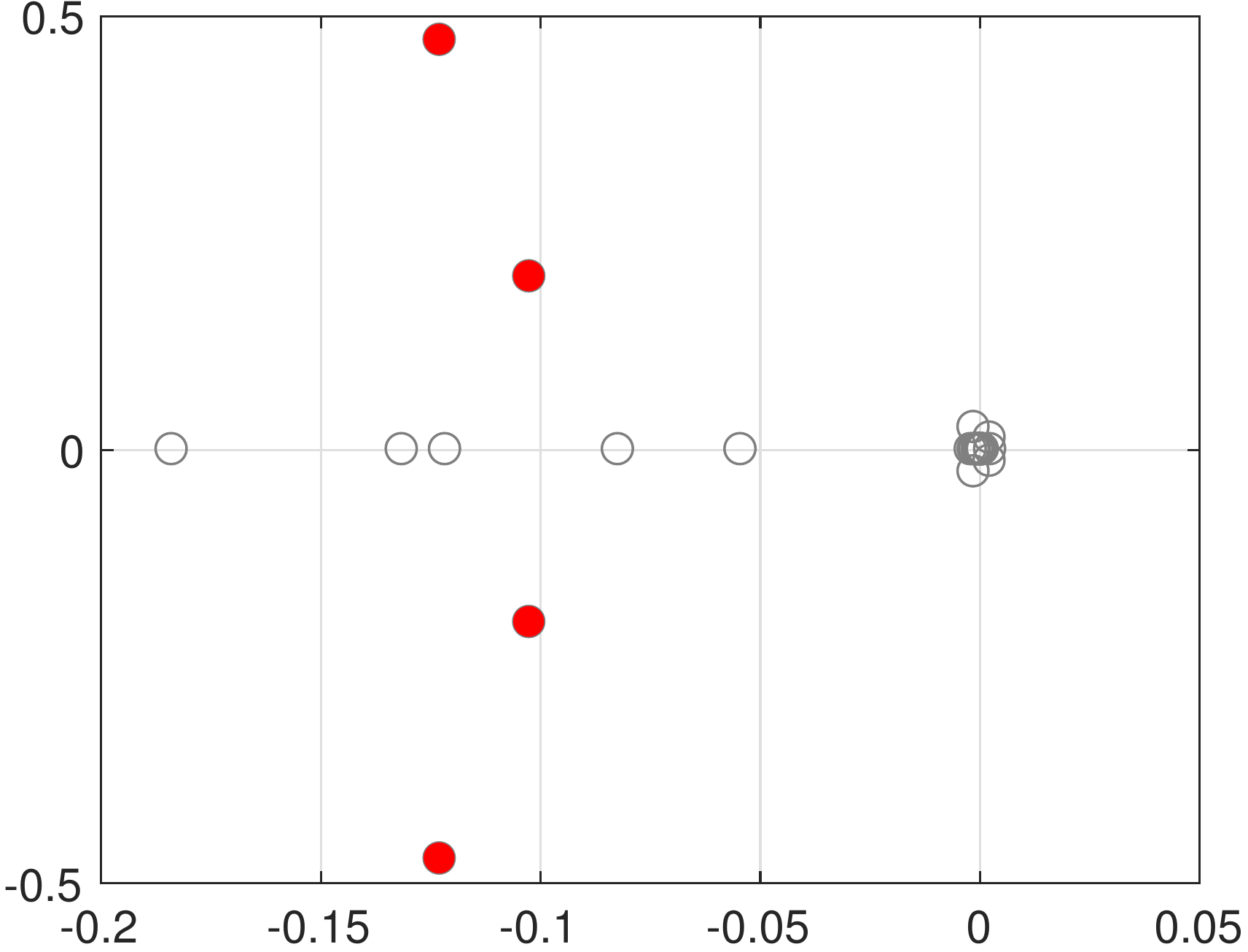}\put(85,65){$\mathbb{C}$}\end{overpic}
		\label{fig:mrdmd_freqs_sst}}
	~
	\subfloat[mrDMD sensors from QR]{\fbox{\includegraphics[width=.315\textwidth]{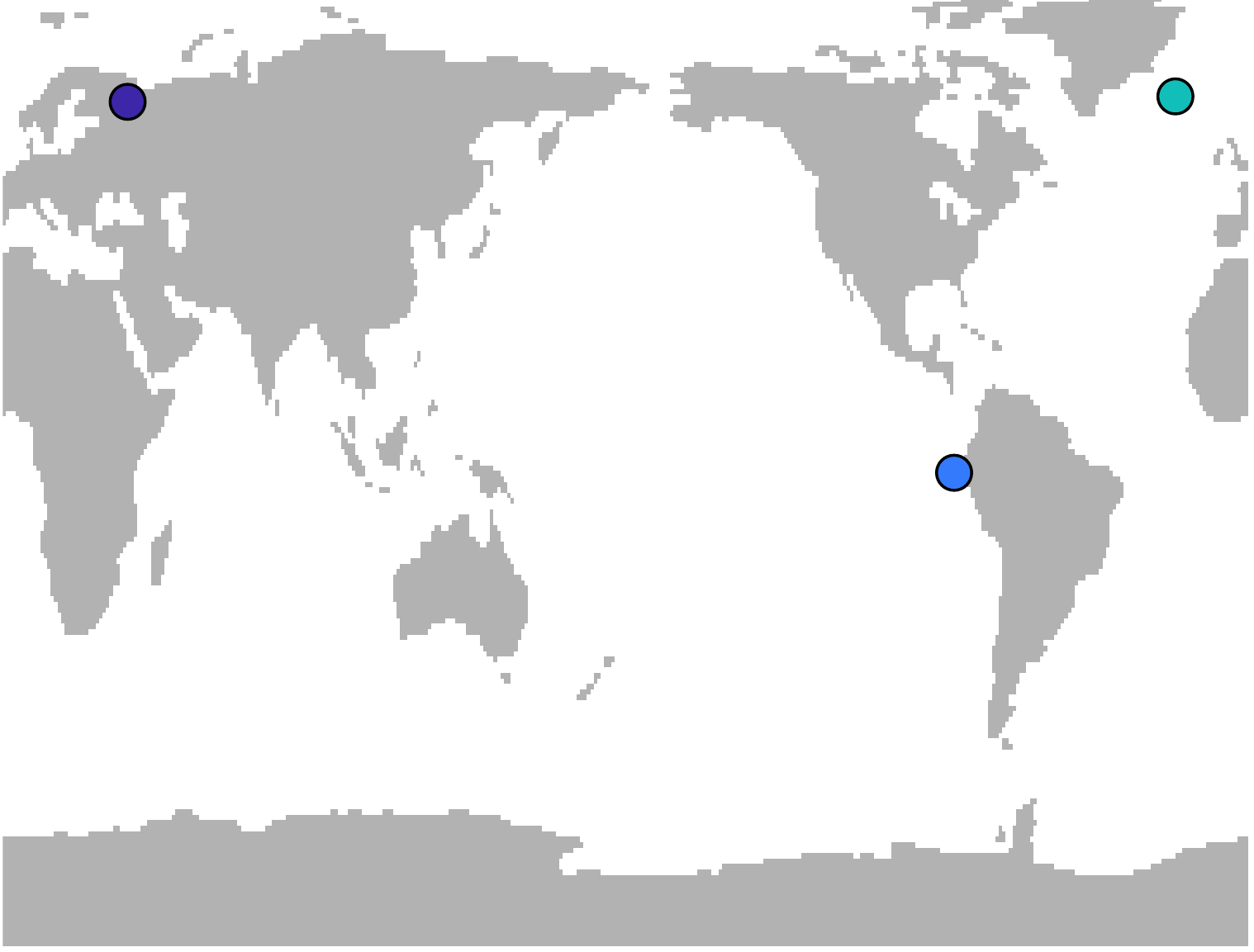}\label{fig:mssp_sens_sst}}}
	~
	\subfloat[Sensor time series]{\includegraphics[width=.3\textwidth]{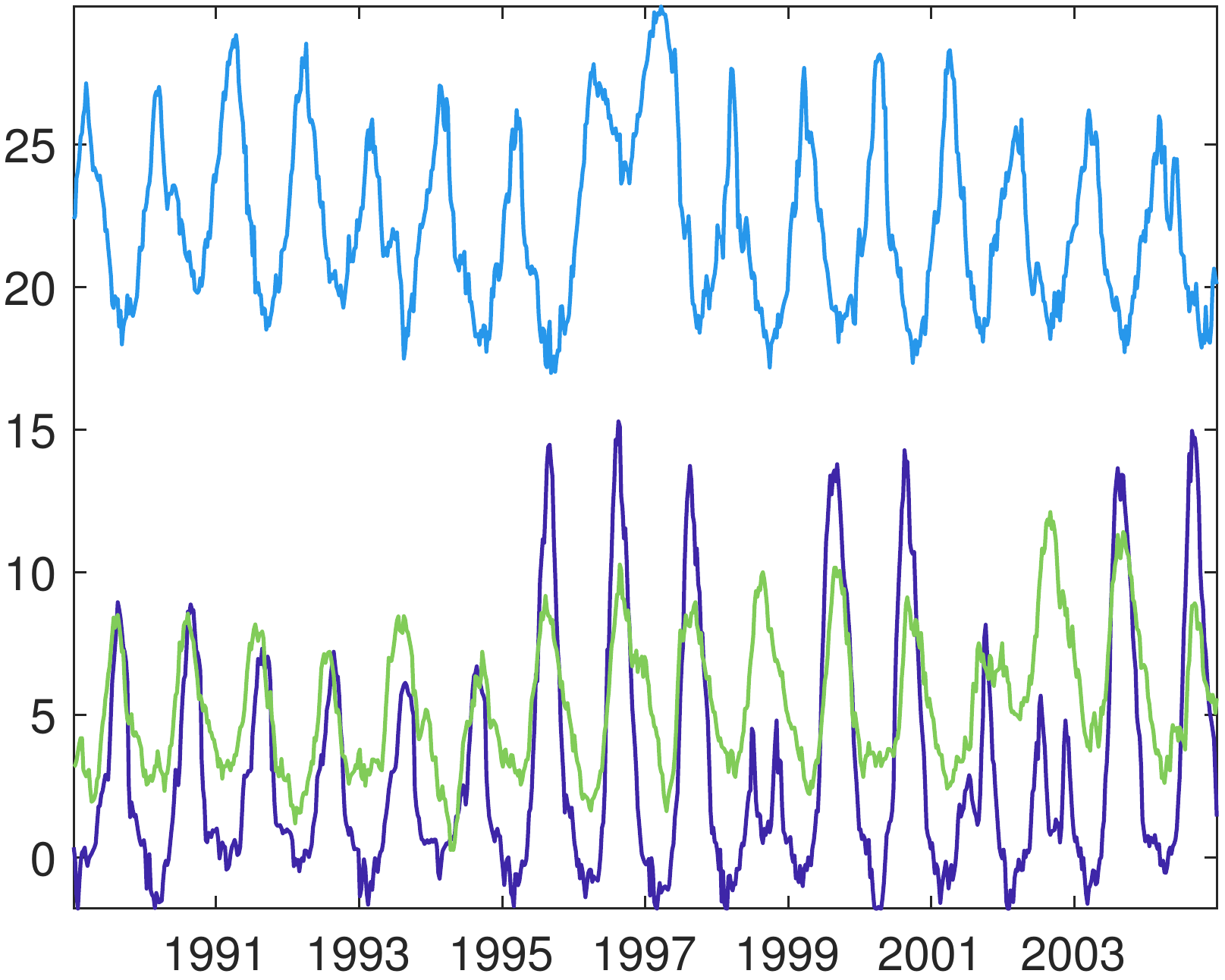}
		\label{fig:mssp_tseries_sst}}
	
		\subfloat[Singular values]{\includegraphics[width=.31\textwidth]{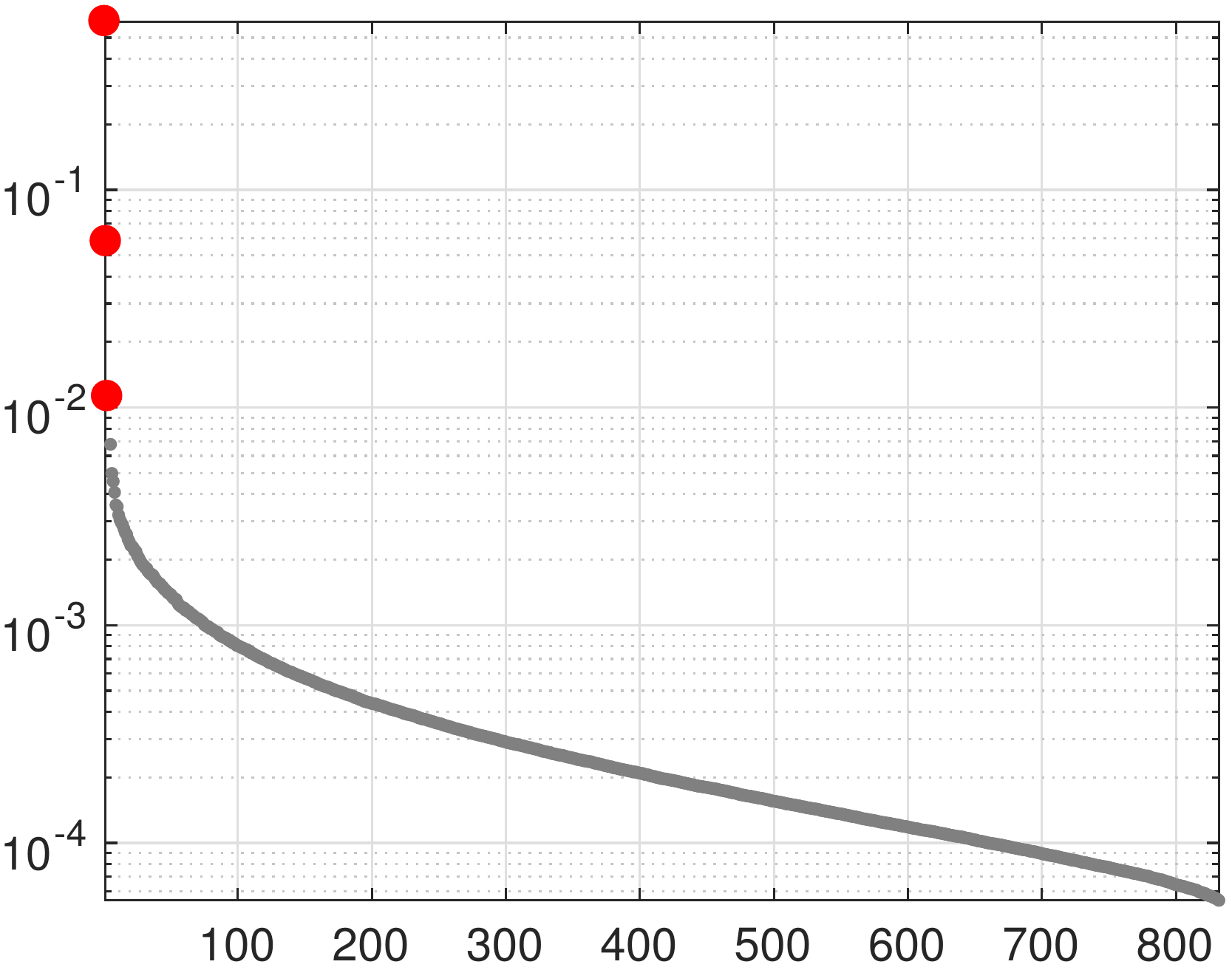}
			\label{fig:singval_sst}}
		~
		\subfloat[POD sensors from QR]{\fbox{\includegraphics[width=.315\textwidth]{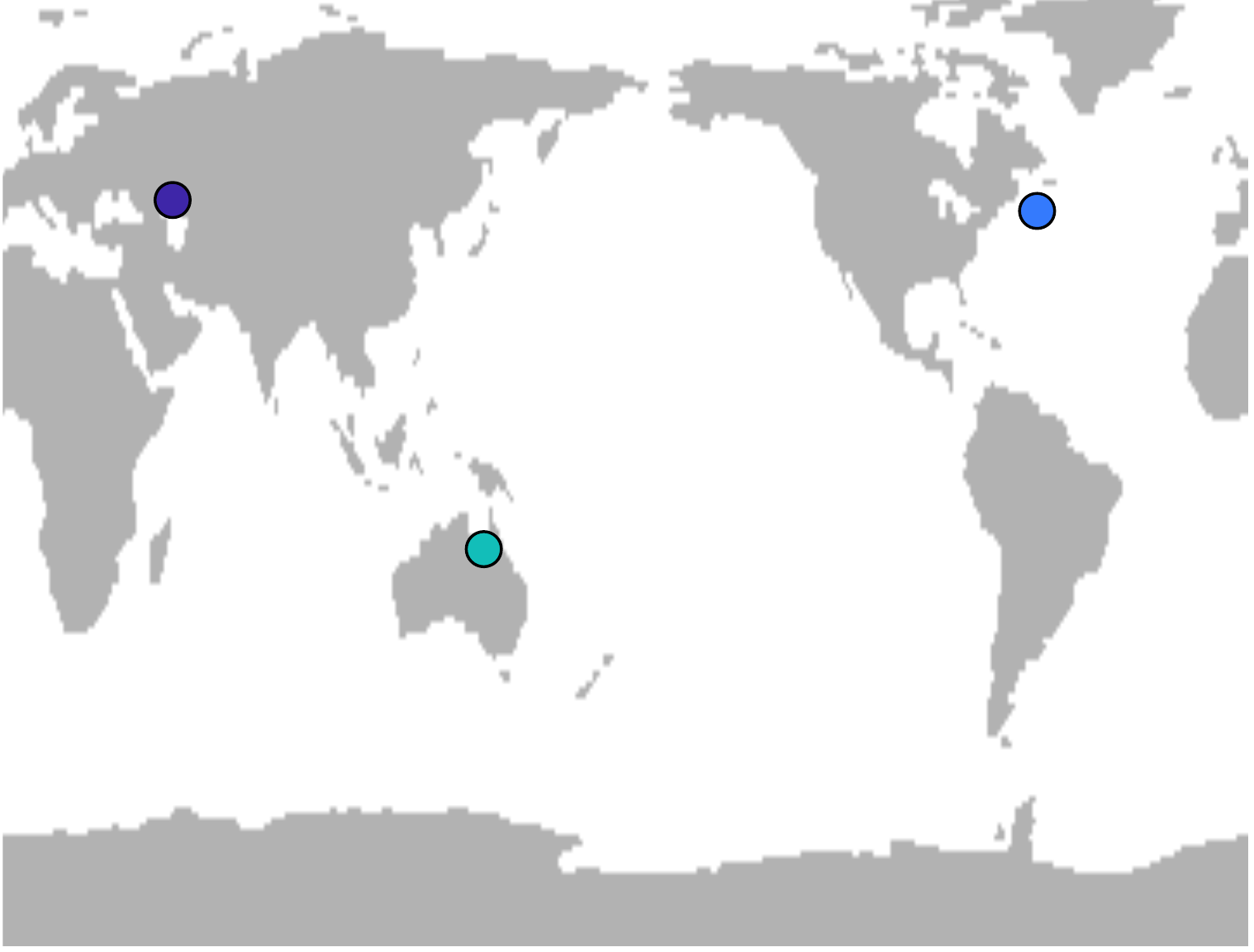}\label{fig:qdeim_sens_sst}}}
		~
		\subfloat[Sensor time series]{\includegraphics[width=.3\textwidth]{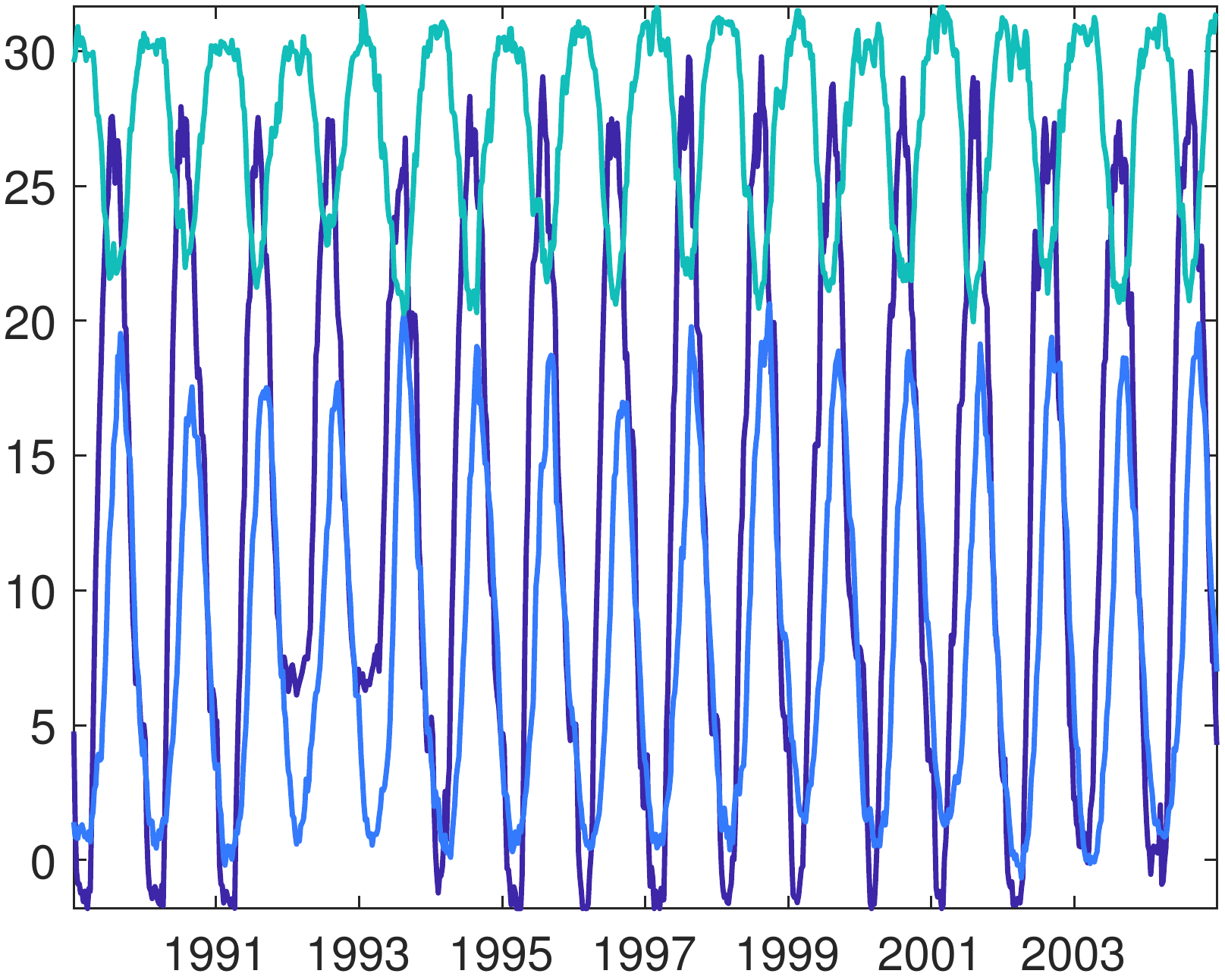}
			\label{fig:qdeim_tseries_sst}}
	\caption{{\bf mrDMD vs. POD sensors.} QR sensors sampled from the ENSO and BSO mrDMD modes (see ~\cref{fig:mrdmd_map_sst}) identify the relevant warming region in the South Pacific, unlike the POD sensors. The time series from mrDMD sensors also display richer temporal trends compared to the dominant annual trends seen in the POD time series.}
	\label{fig:enso_sensors}
\end{figure}
We study now the consequences of multiscale analysis on optimal sampling, which can be readily interpreted as physical sensors in this application. 
\Cref{fig:enso_sensors} compares mrDMD and POD in terms of spectral and spatiotemporal information relayed by mrDMD and POD-based sensors. 
DMD eigenvalues with nonzero imaginary part (blue, \cref{fig:mrdmd_freqs_sst}) identify active oscillations within the 1990-2006 training window. One of the conjugate frequency pairs confirms the true ENSO yearly frequency with $\mbox{Im}(\omega)=0.5$.
On the other hand, the POD spectrum characterizes data covariances and exhibits slow decay due to noisy, low-energy dynamics in the system. This prohibits a low-rank representation of this data using a few POD modes. 
We now compute sensors from both sets of modes. We form $\bPhi_{\Sub}$ from the oscillatory dynamics in the mrDMD spectrum, and the first three POD eigenmodes in the POD spectrum. QR-based sensors are computed for both sets of modes separately. The first observation is that, unlike the toy video example, mrDMD and POD sensors are different. One mrDMD sensor is located over the ENSO band, proving that multiscale characterization is crucial for this system. Meanwhile, POD sensors appear sensitive to temperature extremes occurring along coastlines and inland seas. POD sensor measurements also entirely miss the ENSO phenomenon in the South Pacific. This can be seen from the time series measured at all sensor locations -- mrDMD time series reflect a richer set of dynamics than those of POD, in which only a dominant annual scale is present.
The comparison with POD-based sampling is particularly valuable and demonstrates POD's inability to isolate low-energy intermittent ENSOs in the training window. POD is unable to obtain a low-rank representation of the dynamics because the ocean temperature system is driven by many more degrees of freedom than the toy example. By contrast, mrDMD curtails artificial rank inflation through localized time-frequency analysis.


%
\begin{figure}
	\centering
	\subfloat[Envelope of ENSO coefficient from training window 1990-2006]{\begin{overpic}[width=.8\textwidth]{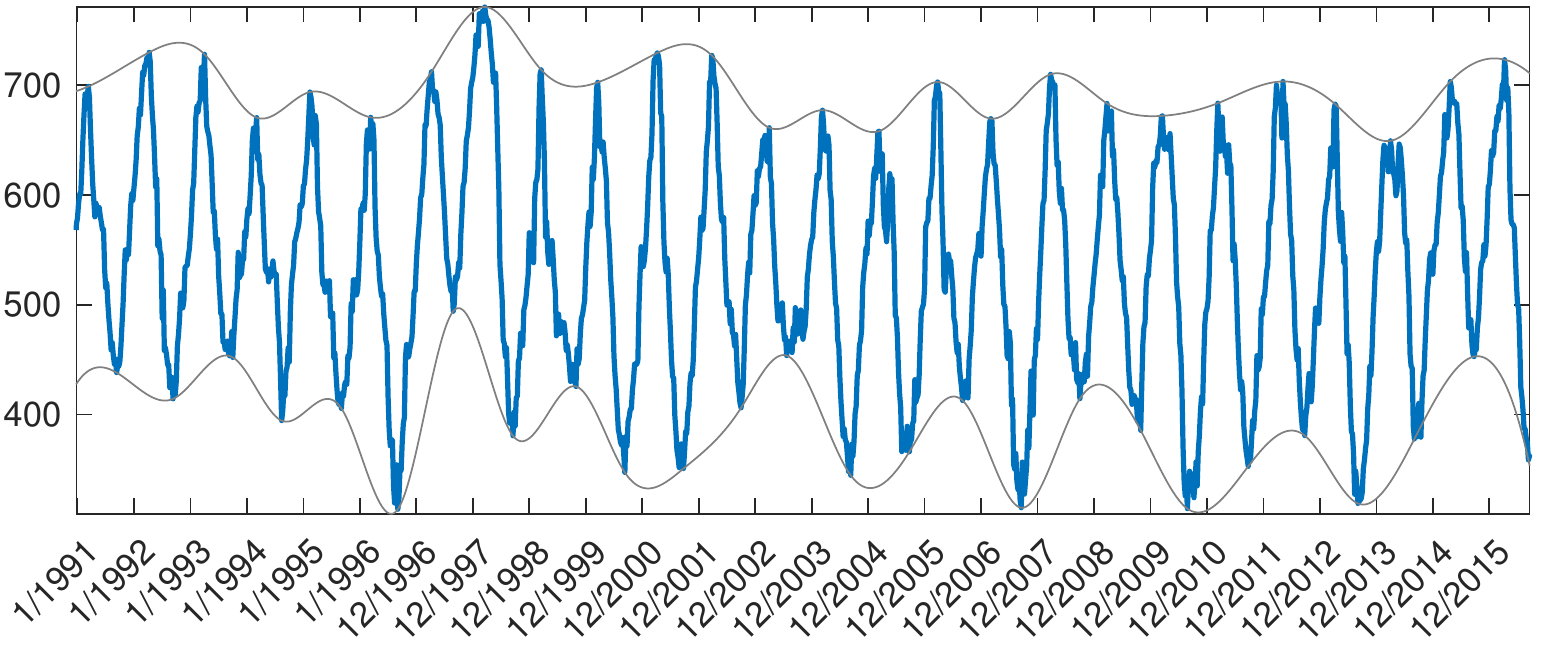}
	\thicklines
	\put(60,8.5){\line(0,1){33}}
	\put(60,37){\vector(-1,0){15}}
	\put(46,38){\small Training}
	\put(60,39){\vector(1,0){15}}
	\put(62,36.5){\small Testing}	
	\put(-3,20){\rotatebox{90}{$a_{ENSO}$}}
	\end{overpic}
		\label{fig:enso_env}}
		\\
	\subfloat[Prediction envelope mean]{\begin{overpic}[width=.8\textwidth]{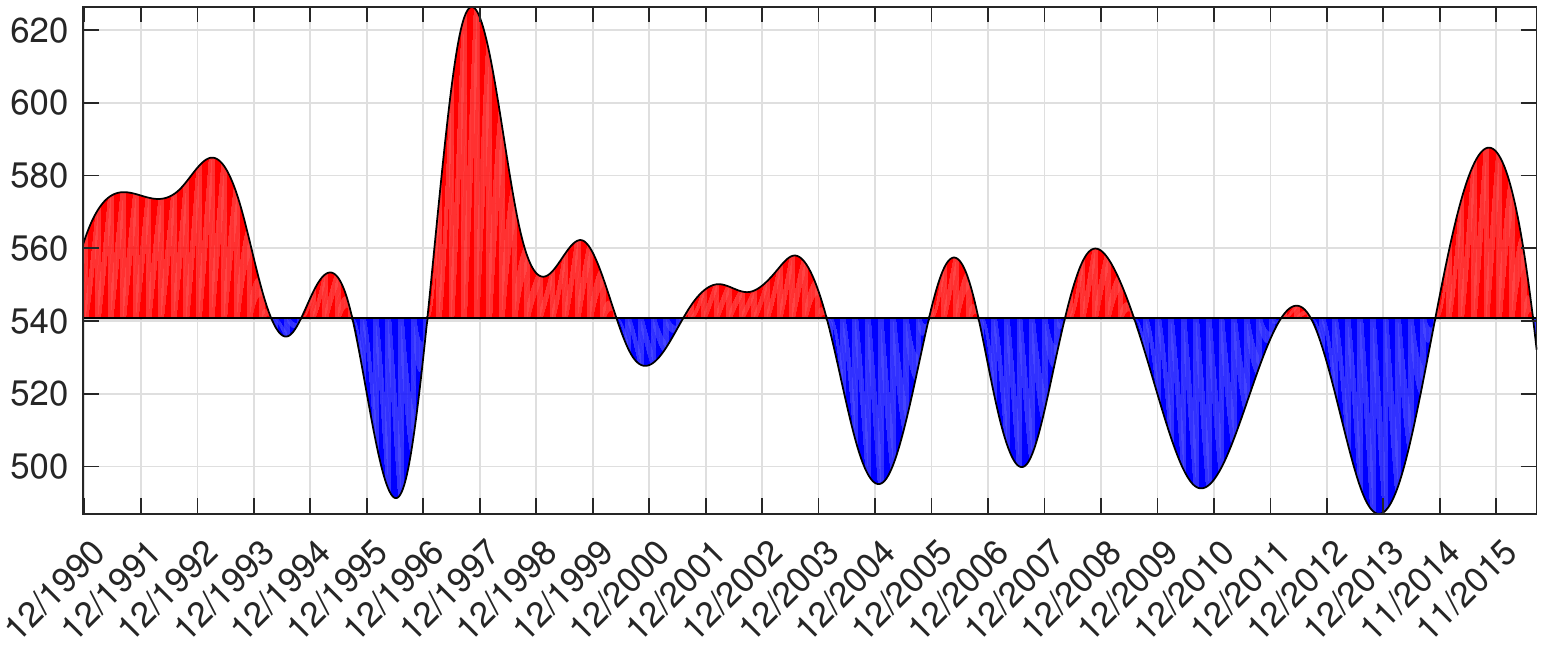}
	\thicklines
	\put(60,8.5){\line(0,1){33}}
	\put(60,37){\vector(-1,0){15}}
	\put(46,38){\small Training}
	\put(60,39){\vector(1,0){15}}
	\put(62,36.5){\small Testing}
	\put(-3,20){\rotatebox{90}{$a_{ENSO}$}}
	\end{overpic}
		\label{fig:enso_pred}}
\caption{{\bf ENSO coefficient prediction}. Coefficient estimation of ENSO mode from only 2 sensors accurately predicts El Ni\~no phenomena. Red peaks closely agree with well-documented strong ENSO events in 1997-99, 2014-2016, as well as with weak ENSO events in 2002-3, 2006-7 and 2009-10. Cooling events (blue) agree with documented La Ni\~na events in 2007-8, 2008-9 and 2010. Small discrepancies in the identified times can be explained by natural variations between the various warming events.}
\end{figure}
\begin{figure}
	\centering
	\includegraphics[width=.8\textwidth]{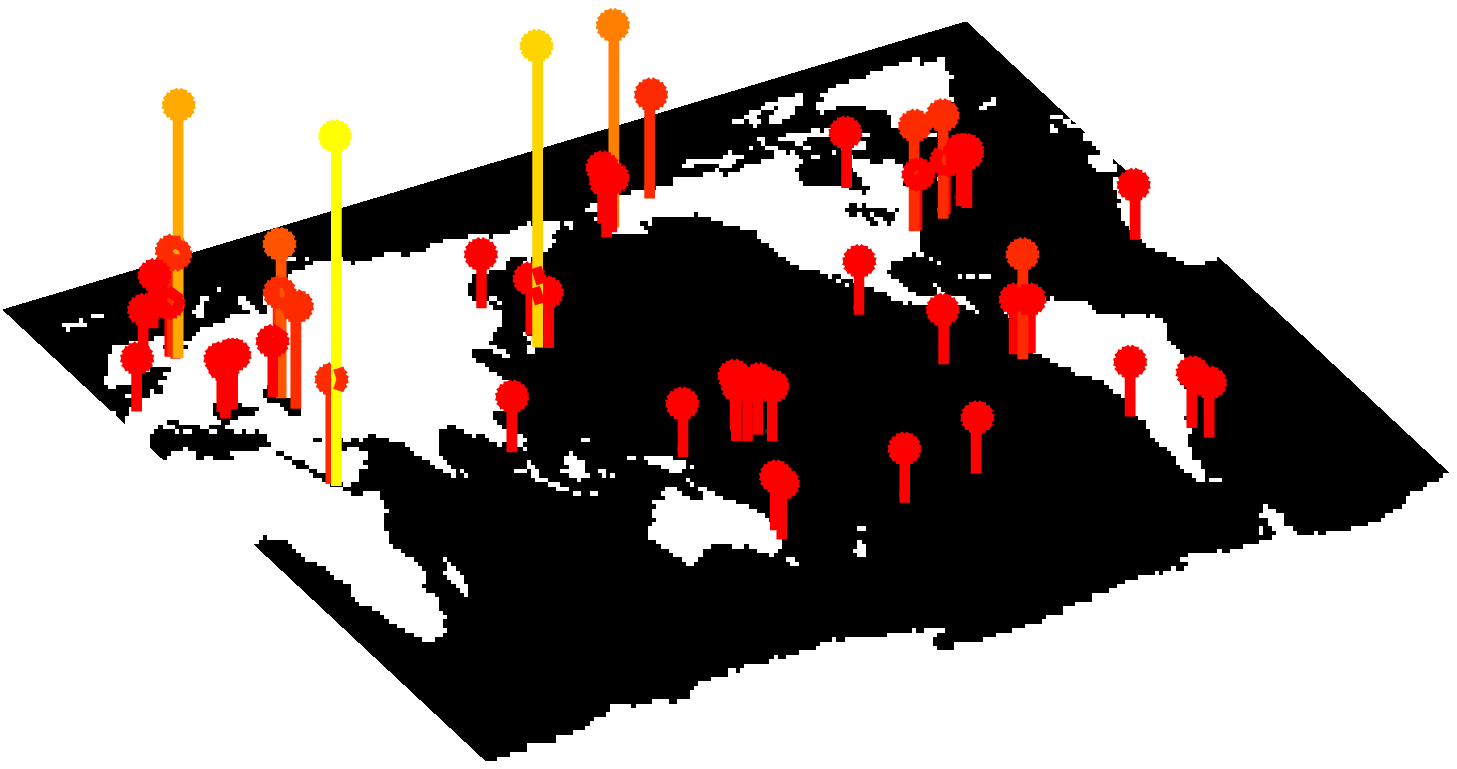}
		\label{fig:enso_sensor_ensemble}
\caption{{\bf Sensor Ensemble}. Time-dependent sensor aggregation pinpoints low-energy locations (\ie ~ENSO sensor along coastal Peru) that would have otherwise been ignored by a global analysis.}
\end{figure}

Optimal sampling can be exploited in a data assimilation-type framework in which prediction is accomplished by parameter estimation followed by full state reconstruction from sensors. 
We estimate ENSO intensity in a validation window using sparse estimation  of library ($\bPhi$) coefficients~\cref{eqn:sparse_pred}. 
\Cref{fig:enso_pred} plots the signal envelope in both the training and validation windows, then traces temperature anomalies about a baseline using the mean of upper and lower envelopes. The resulting measure, which demarcates events above and below baseline (red and blue respectively), closely agrees with well-documented El Ni\~no and La Ni\~na events.  

\begin{figure}[t]
	\centering
	\subfloat[mrDMD sensors from QR]{\frame{\includegraphics[width=.3\textwidth]{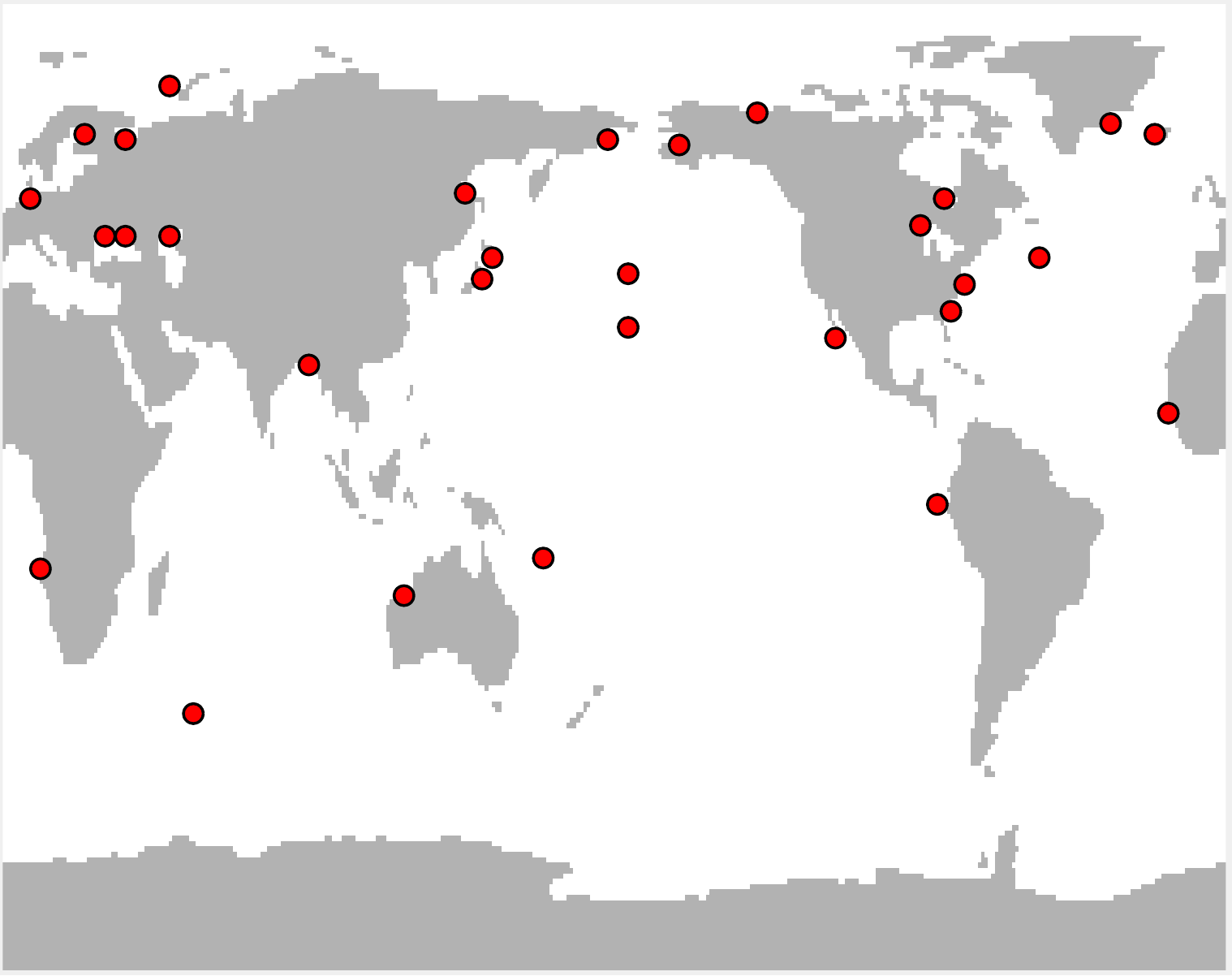}\label{fig:ms_sensors_recon}}}
	~
	\subfloat[$t_1$=6/1992]{\frame{\includegraphics[width=.32\textwidth]{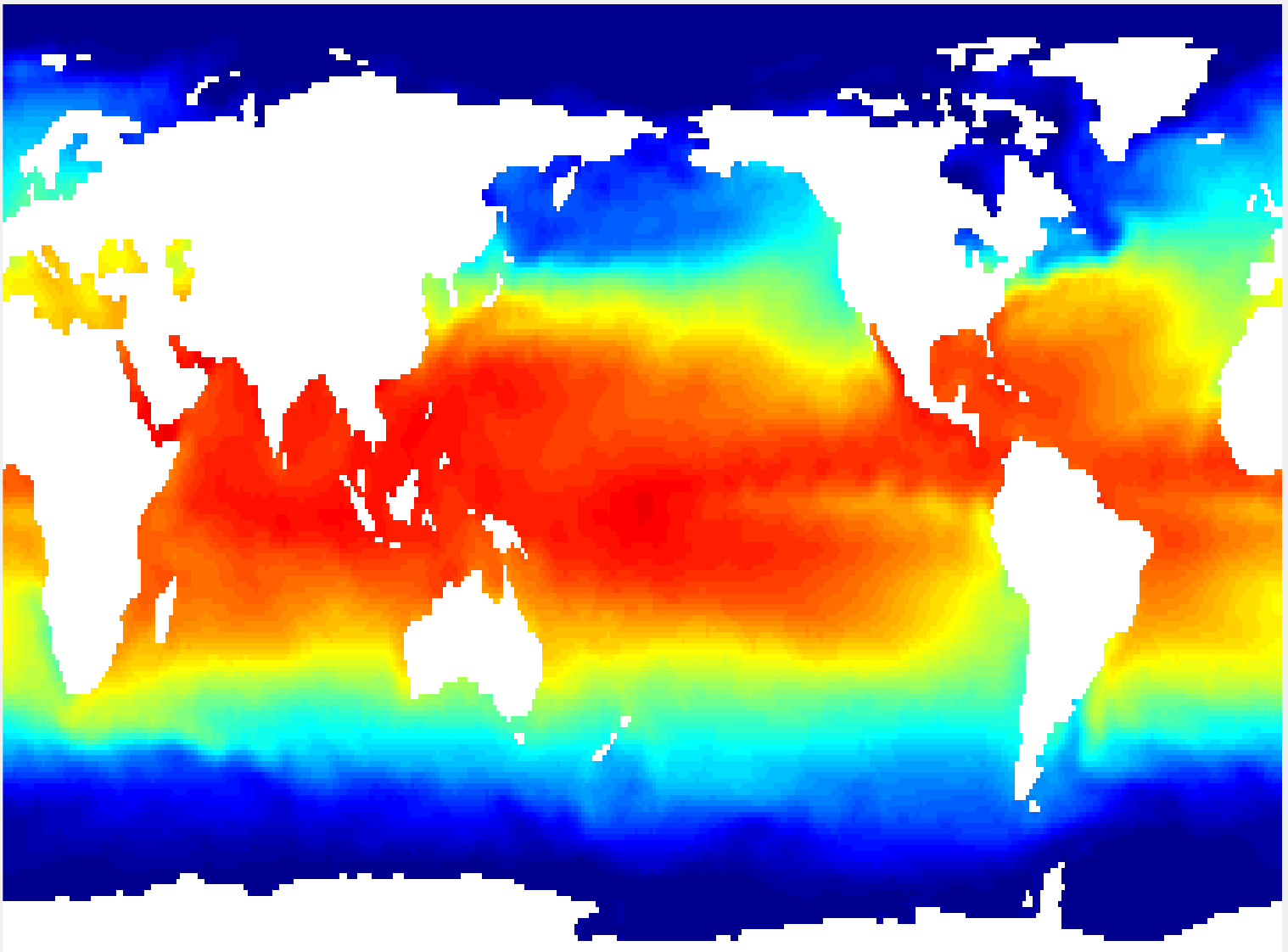}\label{fig:recon1}}}
	~
	\subfloat[$t_2$=6/1996]{\frame{\includegraphics[width=.32\textwidth]{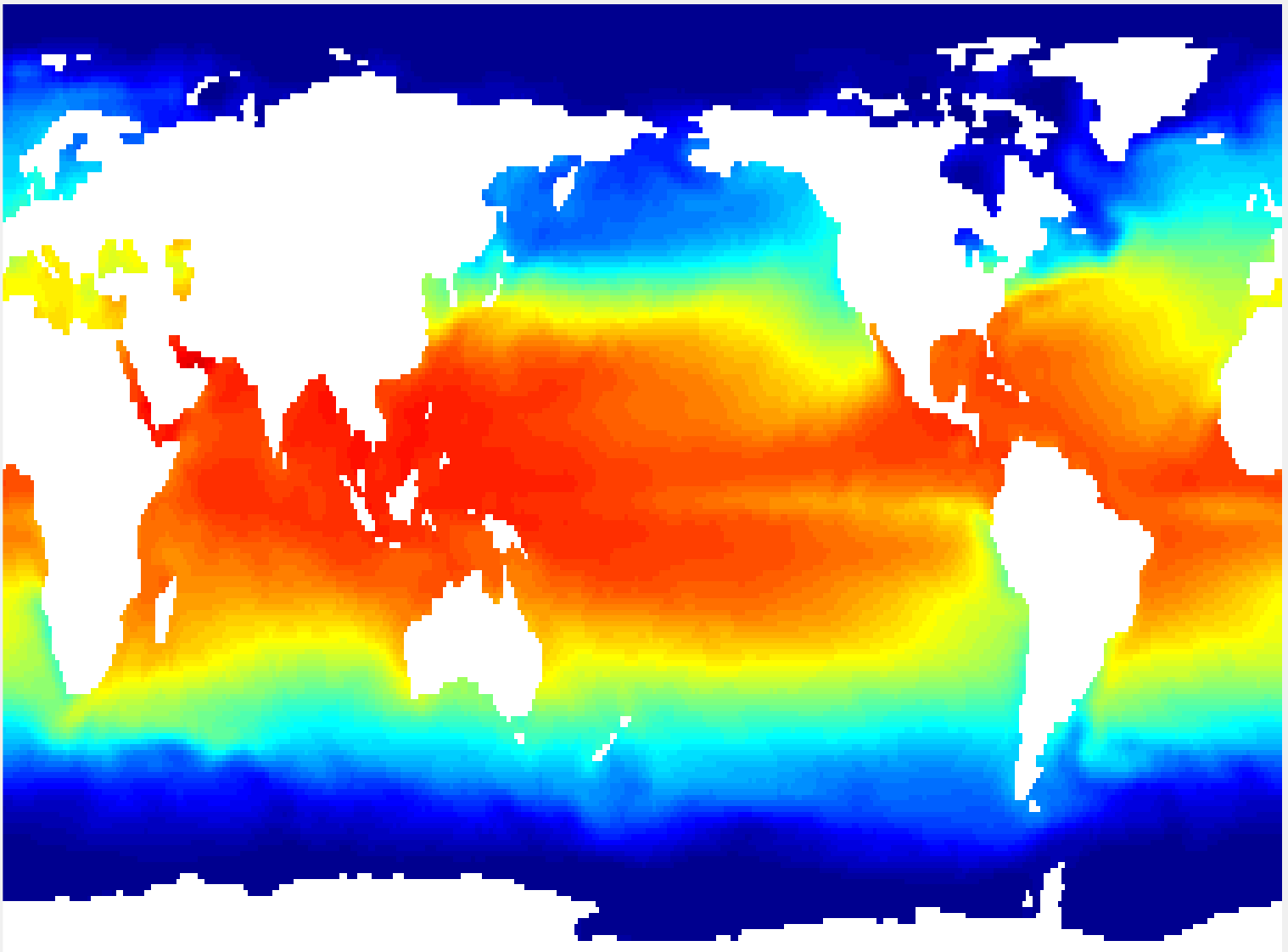}\label{fig:recon2}}}
		
		\subfloat[mrDMD map 1990-2006]{\begin{overpic}[width=.3\textwidth]{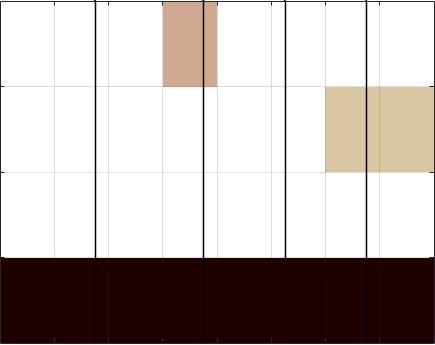}
		\put(14,71){$t_1$}\put(38,71){$t_2$}\put(57.5,71){$t_3$}\put(76,71){$t_4$}
		\end{overpic}\label{fig:mrdmd_map_recon}}
		~
		\subfloat[$t_3$=6/1999]{\frame{\includegraphics[width=.32\textwidth]{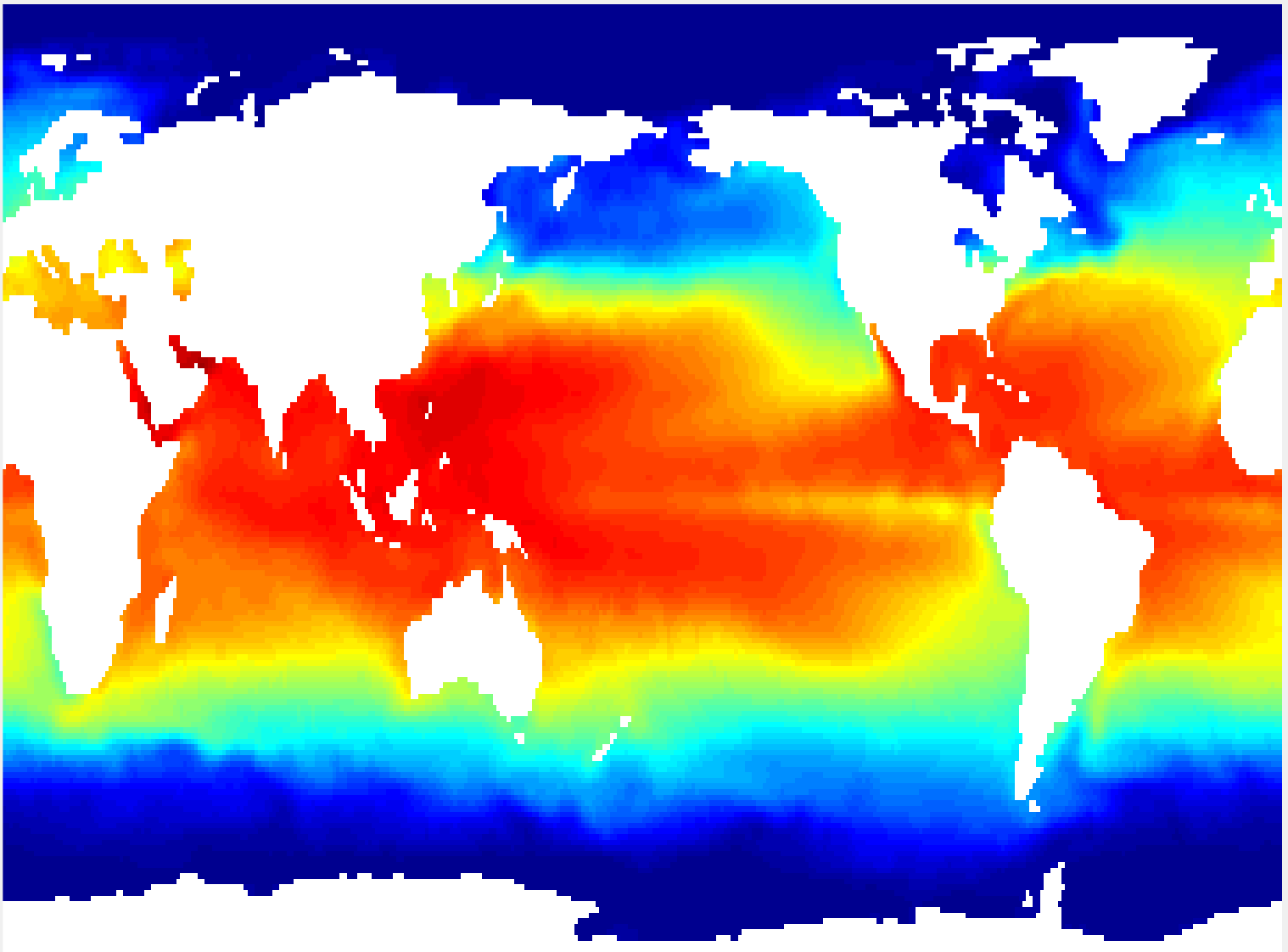}\label{fig:recon3}}}
		~
		\subfloat[$t_4$=6/2002]{\frame{\includegraphics[width=.32\textwidth]{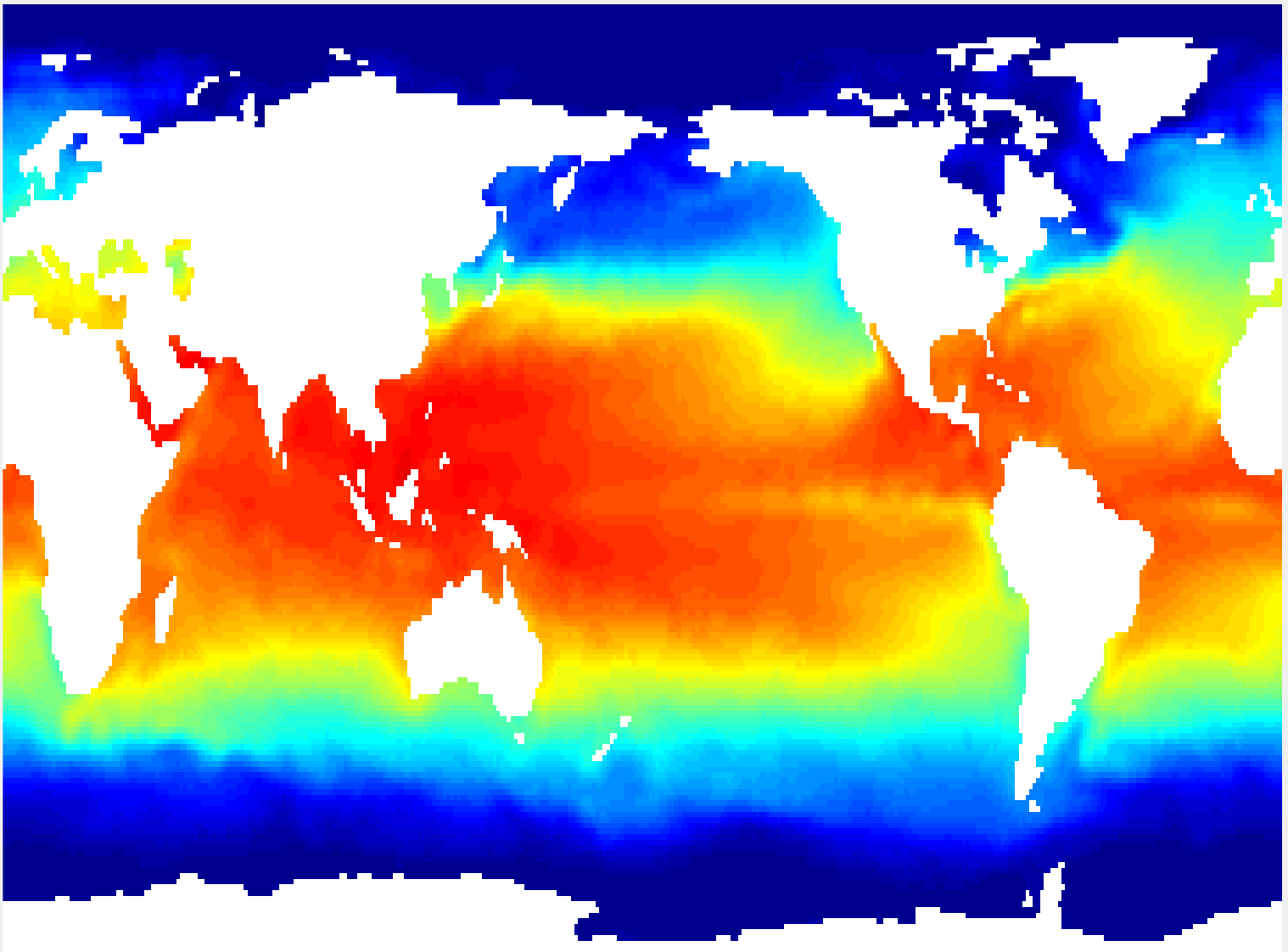}\label{fig:recon4}}}
	\caption{{\bf Reconstruction from 30 sensors in flow field.} MrDMD reconstruction from sensors is reasonably accurate across snapshots selected from different time windows of the decomposition. Each snapshot of 44219 gridpoints is approximated to less than 5\% relative error. }
	\label{fig:enso_reconstruction}
\end{figure}

Full state snapshots of ocean temperatures may be reconstructed using the identified coefficients as described in~\cref{Sec:prediction}. In addition, our sampling framework can be used as a form of multiscale model reduction. Traditionally, ROMs include a POD Galerkin projection to reduce dimensionality of the original system. However, in general DMD modes are not orthonormal so the projection cannot be easily inverted to lift back to the full state. We show here how a few active mrDMD modes can be used for global reconstructions of the state from 30 sensors in the temperature field. \Cref{fig:enso_sensor_ensemble} displays the number of times a spatial location is selected as a QR pivot in each temporal window of the decomposition. As expected, spatial locations corresponding to the ENSO mode are selected with less frequency, reflecting the intermittency of these dynamics. Finally, \cref{fig:enso_reconstruction} displays the results of state reconstruction at different time windows using optimal samples.

%

\section{Discussion}\label{Sec:Discussion}

Despite significant computational and algorithmic advances, the modeling and prediction of multiscale phenomena remains exceptionally challenging.  Much of the difficulty arises from our inability to disambiguate the underlying and diverse set of governing dynamics at various spatiotemporal scales.  Additionally, for large-scale systems, our limited number of sensor measurements can greatly restrict model discovery efforts, thus compromising our ability to produce accurate predictions.  In this manuscript, we have shown that the multiscale disambiguation and limited measurement problem can be simultaneously addressed with principled, algorithmic methods.   Specifically, combining the multiresolution dynamic mode decomposition with QR column pivots, a simple greedy algorithm is demonstrated, which gives accurate spatiotemporal reconstructions and predictions for intermittent multiscale phenomena.

In broader context, our algorithmic developments provide a principled framework for understanding the optimal placement of sensors in complex spatial environments and/or networked configurations.  The greedy QR column pivot selection used for selecting sensor locations leverages the dominant low-rank features of the multiscale physics in order to best predict and reconstruct the spatiotemporal dynamics.  By incorporating a multiresolution analysis tool, i.e. the mrDMD, respect is also given to the diverse temporal phenomena that are often observed in practice.  By systematically accounting for multiscale physics, this placement of sensors performs significantly better than the same limited number of sensors that do not account for these physics. 

The optimized sampling and prediction algorithms developed here have potential for a wide variety of technological applications.  Indeed, with the global increase in sensor networks and sentinel sites for monitoring, for instance, ocean and atmospheric dynamics, disease spread across countries, and/or chemical pollutants, new methods are needed that respect both limited budgets (i.e. limited sensors \& cost) and multiscale physics.  To our knowledge, our algorithms are the first to simultaneously take into account both of these critical components in an efficient, scalable algorithm.



\section*{Acknowledgments}
We are grateful to Travis Askham, Bingni W. Brunton, Joshua L. Proctor, Jonathan H. Tu, and N. Benjamin Erichson for valuable discussions on DMD and sparse sensing.

\bibliographystyle{siamplain}
\bibliography{references}
\end{document}